\documentclass{amsbook}
\usepackage{amssymb}
\usepackage{amsmath}
\newcommand{\Gdag}{G^{\dagger}}

\newtheorem{theorem}{Theorem}[chapter]
\newtheorem{lemma}[theorem]{Lemma}
\newtheorem{prop}[theorem]{Proposition}
\newtheorem{coro}[theorem]{Corollaray}
\newtheorem{main theorem}{MAIN THEOREM}
\theoremstyle{definition}
\newtheorem{definition}[theorem]{Definition}
\newtheorem{example}[theorem]{Example}

\theoremstyle{remark}
\newtheorem{remark}[theorem]{Remark}

\numberwithin{section}{chapter}
\numberwithin{equation}{chapter}

\newtheorem{conjecture}[theorem]{Conjecture}

\newcommand{\B}{\mathfrak{B}}
\newcommand{\Syl}{\mathrm{Syl}}
\newcommand{\skewdim}{\mathrm{skewdim}}
\newcommand{\Def}{\mathrm{Def}}
\newcommand{\rad}{\mathrm{rad}}
\newcommand{\Aut}{\mathrm{Aut}}
\newcommand{\End}{\mathrm{End}}
\newcommand{\im}{•\mathrm{im}}
\newcommand{\SL}{\mathrm{SL}}
\newcommand{\GL}{\mathrm{GL}}
\newcommand{\Sym}{\mathrm{Sym}}
\newcommand{\SU}{\mathrm{SU}}
\newcommand{\PSU}{\mathrm{PSU}}
\newcommand{\PSL}{\mathrm{PSL}}
\newcommand{\ZZ}{\mathbb{Z}}
\newcommand{\QQ}{\mathbb Q}
\newcommand{\NN}{\mathbb{N}}
\newcommand{\FF}{\mathbb{F}}
\newcommand{\HH}{\mathbb{H}}
\newcommand{\oa}{\overline{a}}
\newcommand{\oA}{\overline{A}}
\newcommand{\ob}{\overline{b}}
\newcommand{\opi}{\overline{\pi}}
\newcommand{\oX}{\overline{X}}
\newcommand{\oR}{\overline{R}}
\newcommand{\oJ}{\overline{J}}
\newcommand{\ow}{\overline{w}}
\newcommand{\oW}{\overline{W}}

\newcommand{\oV}{\overline{V}}
\newcommand{\oL}{\overline{L}}
\newcommand{\ov}{\overline{v}}
\newcommand{\Cent}{\mathrm{Cent}}

\newcommand{\characteristic}{\mathrm{char}}
\newcommand{\fata}{\mathbf{a}}
\newcommand{\fatb}{\mathbf{b}}
\newcommand{\F}{\mathfrak{F}}
\newcommand{\id}{\mathrm{id}}
\begin{document}

\frontmatter
\title{On Cubic Action of a Rank one Group}

\author{Matthias Gr\"uninger}
\address{Institut f\"ur Mathematik, Julius-Maximilians-Universit\"at W\"urzburg, Emil-Fischer-Stra\ss e 
30, 97074 W\"urzburg, Germany}
\email{matthias.grueninger@mathematik.uni-wuerzburg.de}
\thanks{The author was supported in part by DFG Grant GZ BA 200//3-1 and by 
ERC grant \# 278469. He likes to thank Bernhard M\"uhlherr and Franz Georg
Timmesfeld for their valuable remarks and suggestions.
}


\date{May 3, 2016}
\subjclass[2010]{20E42,51E24,20G15,17C50}
\keywords{Rank one groups, Moufang sets, quadratic pairs, quadratic action, cubic 
action, quadratic Jordan algebras, pseudo-quadratic forms, Freudenthal triple systems}
\begin{abstract} We consider a rank one group $G = \langle A,B \rangle $ which acts cubically on a 
module $V$, this means 
$[V,A,A,A] =0$ but 
$[V,G,G,G] \ne 0$. 
We have to distinguish whether the group  $A_0 :=C_A([V,A]) \cap C_A(V/C_V(A))$ is trivial or 
not. We show that if $A_0$ is trivial, $G$ is a rank one group associated to a quadratic 
Jordan division algebra. 
If $A_0$ is not trivial (which is always the case  
if $A$ is not abelian), then $A_0$ defines a subgroup $G_0$ of $G$ which acts quadratically on 
$V$. We will call $G_0$ the \textit{quadratic kernel} of $G$.  
By a result of Timmesfeld we have 
$G_0 \cong \SL_2(J,R)$ for a ring $R$ and a special quadratic Jordan division algebra 
$J \subseteq R$. We show that $J$ is either a Jordan algebra contained in a commutative field or 
a hermitian Jordan algebra. In the second case $G$ is the special unitary group of a 
pseudo-quadratic form $\pi$ of Witt index $1$, in the first case $G$ is the rank one 
group for a Freudenthal triple system. These results imply that if $(V,G)$ is a quadratic 
pair such that no two distinct root groups commute and $\characteristic V\ne 2,3$, then $G$ is 
a 
unitary group or an exceptional algebraic group.
\end{abstract}      
\maketitle
\setcounter{page}{4}
\tableofcontents
\mainmatter
\chapter{Introduction}
Abstract rank one groups were introduced by Franz Georg Timmesfeld in 
\cite{T0}. 
\begin{definition}\rm A group $G$ with two distinct, non-trivial, 
 nilpotent subgroups 
$A$ and $B$ is called an 
\textit{abstract rank 
one group with unipotent subgroups $A$ and $B$} if $G =\langle A,B 
\rangle$ and if for all $a 
\in A^*$ there 
is an element $\fatb(a) \in B^*$ with $B^a = A^{\fatb(a)}$. 
\end{definition}

The most common example for a rank one group is the group $\SL_2(K)$ with 
$K$ a (not necessarily 
commutative) field and  $A$ and $B$ the group of lower resp. upper 
unipotent matrices $2 \times 
2$ - matrices over $K$. Alternatively, one can take $G=\PSL_2(K)$. 
Moreover, every semisimple algebraic group $G=\mathbf{G}(K)$ 
defined over a field $K$ of 
relative $K$-rank $1$ is an abstract rank one group. Here, the 
unipotent subgroup $A$ is the unipotent residue of a maximal parabolic 
subgroup of $G$. Similar example are classical groups and groups of 
mixed type of rank $1$ respectively. These are all known examples 
apart from improper rank one groups, by which we mean rank one groups 
$G$ such that $G/Z(G)$ is a sharply $2$-transitive permutation group.

\begin{conjecture}\label{rank 1} If $G$ is an abstract rank one group, 
then $G$ is improper or $G$ is an algebraic or classical group or a 
group of mixed type. 
\end{conjecture}

At the moment, this conjecture seems to be far too difficult to prove.

The concept of an abstract rank one group is closely related to the theory of \textit{Moufang 
sets}. 
\begin{definition}\rm
A \textit{Moufang set} is a pair 
$(X,(U_x)_{x \in X})$, where $X$ is a set with at least $3$ elements and $(U_x)_{x\in X}$ is a 
family of subgroups of $\Sym X$ such 
that $U_x$ fixes $x$ and acts regularly on $X \setminus \{x\}$ and such that $U_x^g =U_{xg}$ 
for 
all $x,y \in X$ and all 
$g \in U_y$. The groups $U_x$ are called \textit{root groups} and the group $G^{\dagger} =\langle U_x| 
x\in X\rangle$ is called the 
\textit{little projective group} of the Moufang set. 
\end{definition}

If $(X, (U_x)_{x \in X})$ is a Moufang set 
with 
nilpotent root groups, then $G^{\dagger}$ is an 
abstract rank one group with unipotent subgroups $U_x$ and $U_y$ for all $x,y \in X$ with $x 
\ne y$. Conversely, if $G=\langle A,B \rangle$ 
is a rank one group, then $(X,(U_x)_{x \in X})$ with $X=\{A^g| g \in G\}$ and $U_{A^g} = A^g 
Z(G)/Z(G)$ is a Moufang set with 
$G^{\dagger} = G/Z(G)$. So a rank one group is a central extension of the 
little projective 
group of a Moufang set with nilpotent root 
groups. 
\newline 
One reason to examine abstract rank one groups is the connection 
between rank one groups and 
\textit{quadratic pairs}. 
A quadratic pair consists of a finite-dimensional $K$-vector space $V$ ($K$ a field of odd 
characteristic) and a subgroup 
$G$ of $\GL_K(V)$ generated by a set $\mathcal{Q}\subseteq \GL_K(V)^*$ of  
quadratic elements, i.e. elements with minimal polynomial $(X-1)^2$, such  
that $\id +\lambda (\sigma-\id)\in \mathcal{Q}$ 
for all $\lambda \in K^*$ and $\sigma\in \mathcal{Q}$. 
Note that if $\characteristic K =2$, then an 
element in $\GL_K(V)$ has minimal polynomial $(X-1)^2$ iff it is an 
involution. Finite quadratic pairs were introduced by Thompson in \cite{Th}, see also 
\cite{Ho}. In \cite{T0} the author proved that the quadratic elements $\sigma\in G$ for 
which $\dim C_V(\sigma)$ is maximal split into up to $1$ pairwise disjoint root groups, 
and two distinct root groups $A$ and $B$ either generate a nilpotent group of class at most 
$2$ or a rank one group with unipotent root groups $A$ and $B$.

Generalising the concept of quadratic pairs, Timmesfeld introduced the following notion.
\begin{definition}\rm(\cite{T2}) Let $G = \langle A,B \rangle $ be a rank one 
group. A $\ZZ G$-module $V$ 
is called 
\textit{quadratic} if $[V,A,A] =0$ but $[V,G,G] \ne 0$. In this case, $G$ 
is called a \textit{quadratic rank 
one group}.
\end{definition}
  
Note that by defintion $V$ is a module over $\ZZ$, but after passing 
over to a suitable submodule of $V$ one may assume that 
$V$ is either an elementary-abelian 
$p$-group or 
torsion free and uniquely divisible, so $V$ can be considered as 
a $kG$-module for $k={\mathbb F}_p$ or 
$k=\QQ$. Both $\characteristic k=2$ and $\dim_k V =\infty$ are allowed. 
\\
The most common example for a rank one group with a 
quadratic module is the rank one group $G=\SL_2(K)$ together with 
its standard module $V=K^2$ for a (not necessarily 
commutative) field 
$K$. One can generalise this as follows: Let $R$ be a ring with $1$. A 
subgroup $J$ of $(R,+)$ 
is called a 
\textit{special quadratic Jordan division algebra} if $1 \in J$ and if $a 
\in J^{\#} =J 
\setminus \{0\}$, then $a$ is invertible in $R$ and $a^{-1}$ is again in 
$J$. Note that the Hua 
identity implies that $bab \in J$ for all $a,b \in J$, see \cite{Hua}. If 
$R$ is a skew field, 
then $R$ itself can be regarded as a special quadratic Jordan division 
algebra; we will denote 
this special quadratic Jordan division algebra by $R^+$.

For $J\subseteq R$ a special quadratic Jordan division algebra set 
\begin{align*} 
A:= \left\{ \left( \begin{array}{cc} 1 & 0 \\ a & 1 \end{array} \right); 
a \in J \right\}, 
B:= \left\{ \left( \begin{array}{cc} 1 & a \\ 0 & 1 \end{array} \right);a 
\in J \right\}\end{align*} 
and $\SL_2(J,R) := \langle A, B \rangle$. 
Then $\SL_2(J,R)$ is a rank one group with unipotent subgroups $A$ and 
$B$ and $V=R^2 $ is a 
quadratic module for 
$\SL_2(J,R)$.
\newline Timmesfeld showed that for every quadratic rank one group $G$ 
there existsa ring $R$ and a special quadratic Jordan division algebra 
$J\subseteq R$ such that $G$ is isomorphic to $\SL_2(J,R)$.
\begin{theorem}[\cite{T2}, Theorem 1.1.]\label{Jordan} Let $G=\langle A,B 
\rangle$ be a rank one group 
which acts quadratically on a module $V$. 
Set $W:= [V,G]/(C_V(G) \cap [V,G])$, $X=[W,A]$ and $R=\End(X)$. Then 
there is a special 
quadratic Jordan division 
algebra $J \subseteq R$ such that $G/C_G(W) \cong \SL_2(J,R)$.
\end{theorem}
Quadratic Jordan division algebras were classified by McCrimmon and 
Zel'manov, see \cite[15.7]{McZ}: 
Every special quadratic Jordan division algebra is isomorphic to the 
quadratic Jordan algebra associated to a skew field, 
an ample hermitian Jordan algebra or a Jordan algebra of Clifford type.  
Note that these families are not disjoint.  
We therefore get a classification of all quadratic rank one groups and 
assert that these groups are all algebraic groups, classical groups or 
groups of mixed type. Therefore Conjecture \ref{rank 1} is true for this 
special case. 

One can generalise the notion of a quadratic action by introducing 
the following definition. If $G$ is a rank one group with unipotent 
subgroups $A$ and $B$ and $V$ is a $G$-module, then we say that 
$V$ has \textit{length} $n$ if $n\geq 2$ is the smallest positive integer 
such that $[V,\underbrace{A,\ldots,A}_{n-\text{times}}]=0$ but 
$[V,\underbrace{G,\ldots,G}_{n-\text{times}}]\ne 0$. So a quadratic 
module is a module of length $2$. (Modules of lenght $1$ do not exist 
since $G$ is the normal clousure of $A$). Every known proper rank one 
group has in fact a module of finite length. For this reason, a 
classification of all rank one groups having a module of "small" length 
might be an important step to prove Conjecture \ref{rank 1}.

In this paper we deal with \textit{cubic modules}, by which we mean 
modules of length $3$. As in the quadratic case we will later see that 
one can assume that a cubic module $V$ for a rank one group $G$  is a 
$KG$-module for a field $K$. Both $\characteristic K =2$ and 
$\dim_K V = \infty$ are allowed. 
Typical examples of cubic modules are pseudoquadratic spaces of Witt 
index $1$ or the endomorphism algebra of a quadratic rank one group. 
Furthermore, some exceptional algebraic groups, namely those 
 of type ${}^{3}\mathrm{D}^9_{4,1},
{}^{6}\mathrm{D}^9_{4,1} {}^{2}\mathrm{E}^{35}_6, 
\mathrm{E}^{66}_{7,1}$ or 
$\mathrm{E}_{7,1}^{78}$, possess a cubic module.

The class of all pairs $(V,G)$ with $G$ 
a rank one group with unipotent subgroups $A$ and $B$ and $V$ 
a cubic module for $G$ falls into three subclasses. To distinguish 
between these cases, we introduce a subgroup $A_0=C_A([V,A])
\cap C_A(V/C_V(A))$. Then $A' \leq A_0  \leq A$ and if $A_0\ne 1$, 
then $A_0$ is a \textit{root subgroup} of $A$. This means that there is 
a subgroup $G_0\leq G$ such that $A_0=A\cap G_0$ and that 
$G_0$ is a rank one group with unipotent subgroups $A_0=A\cap G_0$ and 
$B_0:=B\cap G_0$. Since $G_0$ acts quadratically on $V$, we will 
call $G_0$ the \textit{quadratric kernel} of $G$. 
 By Timmesfeld's result we have 
$G_0\/C_{G_0}(V) \cong\SL_2(J,R)$ for a quadratic Jordan division algebra inside a 
ring $R$. 

We will see that one of the following three cases holds:

\begin{itemize}
\item In the first case, $A_0=1$ is trivial and hence $A$ is abelian.
\item In the second case, $A_0\ne 1$ and $J$ is a commutative 
quadratic Jordan division algebra, hence there is a field $K$ such 
that $J=K$ and $G=\SL_2(K)$ or $\characteristic K =2$ and 
$K^2 \subseteq J \subseteq K$. Here, if $\characteristic K\ne 2$, 
then $A$ is not abelian. We will say that $G_0$ is of \textit{commutative 
type}.
\item In the third case, $A_0\ne 1$ and $J$ is a hermitian Jordan algebra 
inside a non-commutative skew field $K$ such that $J$ generates $K$ as 
a field. In this case, $A$ is not abelian. We will say that $G_0$ 
is of \textit{hermitian type}.
\end{itemize}

We can now state our three main results. 
\begin{main theorem} Let $G$ be a rank one group acting cubically 
on a module $V$ such that the quadratic kernel is trivial. Then 
there is a quadratic Jordan division algebra $J$ with $G/Z(G) \cong 
\PSL_2(J)$.
\end{main theorem}

\begin{main theorem} Let $G$ be a rank one group such that $G$ 
acts cubically on a module $V$. Suppose that the quadratic kernel 
$G_0$ is of hermitian type. Then $G$ is isomorphic to the special 
unitary group for a pseudoquadratic form of Witt index $1$.
\end{main theorem}

For the third case we only have a result for characteristic not $2$ or 
$3$,
\begin{main theorem} Let $G$ be a rank one group acting cubically one 
a module $V$. If the quadratic kernel $G_0$ is isomorphic to 
$\SL_2(K)$ with $K$ a commutative field of characteristic not $2$ or 
$3$, then $G$ is isomorphic to a rank one group for an anisotropic 
Freudenthal triple system.
\end{main theorem}

Freudenthal triple systems of finite dimension 
correspond to structurable division algebras of skewdimension $1$. 
In a current paper it is shown that semisimple 
algebraic groups of relative rank $1$ over a field $K$ with 
$\characteristic K \ne 2,3$ correlate to rank one groups defined via 
structurable division algebras, see \cite{BDMS}. 
Therefore these results strongly suggest that the following conjecture is 
true:

\begin{conjecture} Let $G$ be a cubic rank one group with 
unipotent subgroups $A$ and $B$. Then we have:
\begin{itemize}
\item[(i)] If $A$ is abelian, then there is a quadratic Jordan divsision 
algebra $J$ such that $G\cong \PSL_2(J)$.
\item[(ii)] If $A$ is not abelian, then one of the following holds:
\begin{itemize}
\item[(iia)] There is an involutory set $(K,K_0,*)$, a $K$-vector space 
$V$ and an anisotropic pseudo-quadratic form $\pi:V\to K/K_0$ such 
that $G\cong\SU(\pi)$.
\item[(iib)] There is a field $K$ such that $G$ is a semisimple algebraic 
group defined over $K$ of type ${}^{3}\mathrm{D}^9_{4,1},
{}^{6}\mathrm{D}^9_{4,1} 
{}^{2}\mathrm{E}^{35}_6, 
\mathrm{E}^{66}_{7,1}$ or $\mathrm{E}^{77}_{7,1}$.
\end{itemize}
\end{itemize}
\end{conjecture} 

The paper is organised as follows: In chapter 2 we repeat some facts 
about rank one groups and 
some facts concerning algebra and geometry that we will need later. 
This includes quadratic Jordan algebras, semi-prime rings, skew fields 
with involution, pseudo-quadratic forms, Moufang quadrangles, 
quadrangular algebras and Freudenthal triple systems. 
In chapter 3 we prove some elementary facts about 
cubic action. Especially we introduce the quadratic kernel of 
a cubic rank one group. In chapter 4 we introduce some examples. 
chapter 5 we establish a kind of normal form for cubic modules. 
In chapter 6 we will show how to obtain an irreducible cubic module 
by an arbitrary cubic module. chapter 7 and 8 deal with 
cubic rank one groups with trivial quadratic kernel. 
In chapter 7 we will show that such a group is always a rank one 
group for a quadratic Jordan division algebra. In chapter 8 we will 
give a characterisation of the adjoint module. chapter 9 deals 
with cubic rank one groups with non-trivial quadratic kernel. 
This case is divided into two subcases. In the first case the quadratic 
kernel is a hermitian rank one group; this case is dealt in chapter 10
and 11 where it is shown that the cubic rank one groups appearing in this 
case are unitary groups for a pseudo-quadratic form of Witt index $1$. 
The final chapter 12 deals with the other case where the quadratic kernel 
is a special linear group of dimension $2$ over a commutative case. 
Here we show that in this case we obtain rank one groups for a 
Freudenthal triple system provided the characteristic different from $2$ 
or $3$. 

\subsection*{Notation}
\begin{itemize}
\item For a multiplicative group $G=(G,\cdot,1)$, we set $G^* = G \setminus \{1\}$, while for an 
additive group $G=(G,+,0)$ we write $G^\# =G \setminus \{0\}$.
\item For a ring $R$, we denote the set of non-zero elements of $R$ by $R^\#$ and the set 
of units by $R^*$. 
We denote the centre of $R$ by $\Cent R :=\{a\in R|ba=ab \ \forall a\in R\}$ and $\Cent_R(a):=\{b\in R|ba=ab\}$ and for $a\in R$.
\item If a group $G$ acts on a set $\Omega$, we will also apply elements of $g$ on the right, 
i.e. we will write $\omega g$ instead of $g(\omega)$.
\item If a group $G$ acts on an abelian group $(V,+)$, then for subgroups $A,B \leq G$ we write
$[V,A,B]$ for $[[V,A],B]$.
\end{itemize}
\chapter{Preliminaries}
\section{Moufang sets}
We begin with some properties of Moufang sets. For further information and proofs we recommend 
\cite{DS} as an introduction to Moufang sets. \\
As mentioned in the introduction, a Moufang set $\mathbb{M}=(X,(U_x)_{x\in X})$ 
consists of a set $X$ with $|X| \geq 3$ and 
a family $(U_x)_{x\in X}$ of subgroups of $\Sym X$ such that for all $x,y \in X$ we have 
\begin{enumerate}
\item $U_x$ fixes $x$ and acts sharply transitively on $X\setminus \{x\}$.
\item For all $a\in U_x$ we have $U_x^a =U_{xa}$.
\end{enumerate}
The groups $U_x$ are called \textit{root groups}, the group $\Gdag :=\langle U_x| x\in X\rangle 
\leq \Sym X$ is called the \textit{little projective group} of $\mathbb{M}$. The group 
$\Gdag$ is acts $2$-transitively on $X$; the Moufang set $\mathbb{M}$ is called \textit{proper} 
if $G^{\dagger}$ is not sharply $2$-transitive.
\begin{example}\rm Let $K$ be a field and $G=\mathbf{G}(K)$ 
be a semisimple 
linear algebraic group of relative $K$-rank one and let $X$ be the set of 
all proper $K$-parabolic subgroups of $G$. For $P\in X$ let 
$U_P$ be the group of all $K$-rational points in the unipotent radical of 
$P$. Then $(X,(U_P)_{P\in X})$ is a Moufang set which will be denoted by 
$\mathbb{M}(G)$.
\end{example}

Every Moufang set can be constructed in the following manner: Let $(U,+)$ be a (not necessarily 
abelian) group, $X:=U \dot{\cup} \{\infty\}$ and $\tau \in \Sym X$ an element which 
interchanges $0$ and $\infty$. For $a\in U$ let $\alpha_a \in \Sym X$ be the permutation 
defined by $b\alpha_a =b+a$ for $b\in U$ and $\infty \alpha_a =\infty$. Set $U_{\infty}:=
\{\alpha_a|a\in U\}$. Then $U_{\infty}$ is a subgroup of $(\Sym X)_{\infty}$ isomorphic to $U$. 
Now define $U_0:=U_{\infty}^{\tau}$, $U_a :=U_0^{\alpha_a}$ for $a\in U$ and 
$\mathbb{M}(U,\tau):=(X,(U_x)_{x\in X})$. This is, however, not always a Moufang set. 
To detect whether $\mathbb{M}(U,\tau)$ one regards the following maps:
\begin{definition}\rm Let $a\in U^{\#} = U \setminus \{0\}$. Then $\mu_a:= \alpha_a 
\alpha_{-a\tau^{-1}}^{\tau} \alpha_{-(-a\tau^{-1})\tau}$ is called the $\mu$-map corresponding 
to $a$, and $h_a :=\tau \mu_a$ is called the Hua map corresponding to $a$.
\end{definition}
One easily sees that $\mu_a$ interchanges $0$ and $\infty$ and hence $h_a$ fixes these two 
elements.
\begin{theorem}
The following are equivalent:
\begin{enumerate}
\item $\mathbb{M}(U,\tau)$ is a Moufang set.
\item $U_0^{\mu_a} =U_{\infty} $ and $U_{\infty}^{\mu_a} =U_0$ for all $a\in U^{\#}$.
\item $U_{\infty}^{h_a} =U_{\infty}$ for all $a \in U^{\#}$.
\item For all $x,y,a\in U$ with $a\ne 0$ one has $(x+y)h_a =xh_a +yh_a$, i.e. $h_a$ induces 
an endomorphism on $U$.
\end{enumerate} 
\end{theorem}
\begin{proof}
See \cite[Theorem 3.1 and Theorem 3.2]{DW}.
\end{proof}
If $\mathbb{M}(U,\tau)$ is a Moufang set, then we have $\mathbb{M}(U,\tau)=\mathbb{M}(U,\mu_a)$ 
for all $a\in U^{\#}$. Hence we can assume that $\tau =\mu_a$ for some $a\in U^{\#}$.
\\
Set $H:= \langle \mu_a \mu(b)| a,b \in U^{\#}\rangle$. Then $H$ is called the 
\textit{Hua subgroup} of $\mathbb{M}(U,\tau)$. By \cite[Theorem 3.1(ii)]{DW} we have 
\begin{prop}\label{Hua Moufang}
$H=\Gdag_{0,\infty}$.
\end{prop}
Thus we get
\begin{theorem}\label{proper Moufang set}
For a Moufang set $M(U,\tau)$ the following is equivalent:
\begin{enumerate}
\item $M(U,\tau)$ is improper.
\item $H=1$.
\item $\mu_a=\mu_b$ for all $a,b\in U^{\#}$.
\end{enumerate}
\end{theorem}

We list some of the properties of the $\mu$-maps we will frequently use.
\begin{lemma}\label{properties mu map} 
Suppose that $\mathbb{M}(U,\tau)$ is a Moufang set. Let $a\in U^{\#}$.
\begin{enumerate}
\item $\mu_a= \alpha_{(-a)\tau^{-1}}^{\tau} \alpha_a \alpha_{-a\tau^{-1}}^{\tau}$.
\item $\mu_a $ is the unique element in the double coset $U_0 \alpha_a U_0$ which interchanges 
$0$ and $\infty$.
\item $\mu_a^{-1} =\mu_{-a}$.
\item $\mu_{ah} =\mu_a^h$ for all $h\in H$.
\item If $\mathbb{M}(U,\tau) =\mathbb{M}(U,\tau^{-1})$, then $\mu_{a\tau} =\mu_{-a}^\tau$;
especially we have $\mu_{a \mu_b} =\mu_a^{\mu_b}$ for all $b\in U^{\#}$.
\end{enumerate}
\end{lemma}
\begin{proof}
See \cite[4.3.1]{DS}.
\end{proof}
For $a\in U^{\#}$ we set $\sim a =(-a\tau^{-1})\tau$. 
Then we have
\begin{lemma}\label{formula sim}
\begin{enumerate}
\item $\sim a = -(-a)\mu_a$. Especially $\sim a$ does not depend on the 
choice of $\tau$.
\item $(-a)\tau = \sim (a\tau)$. 
\item $\sim (ah) =(\sim a) h$ for all $h\in H$.
\item $\mu_{\sim a} = \mu_{-a}$.
\item $a\mu_a = \sim - \sim a$ and $a\mu_{-a} =-\sim -a$. 
\item If $\mathbb{M}(U,\tau)$ is improper, then $\sim - \sim a =-\sim -a $ for all $a\in 
U^{\#}$.
\end{enumerate}
\end{lemma}
\begin{proof}
Part (a) is \cite[4.3.1 (6)]{DS}. Part (b) follows immediately from the definition of $\sim a$. 
Part (c) follows from \ref{properties mu map} and (a), (d) from \ref{properties mu map}(e). 
By (a) and (b) we have 
${(\sim a)\mu_{\sim a}} = (\sim a)\mu_a^{-1} = 
(-(-a)\mu_a)\mu_a^{-1} = \sim -a $. Replacing $a$ by $\sim a$ we get the first part of (e). 
The second part follows from applying the formula in (a) for $-a$ instead of $a$. Finally, 
if $\mathbb{M}(U,\tau)$ is improper, then $\mu_a=\mu_{-a}=\mu_a^{-1}$, so (f) follows.   
\end{proof}
\begin{definition}\rm Let $\mathbb{M}(U,\tau)$ be a Moufang set and $a\in U^{\#}$. Then 
$a$ is called \textit{special} if $(-a)\tau^{-1} = -a\tau^{-1}$. The Moufang set $\mathbb{M}(U,\tau)$ 
is called special if all non-trivial elements are special.
\end{definition}
\begin{lemma}\label{special elements in MS} The following are equivalent for $a\in U^{\#}$.
\begin{enumerate}
\item $a$ is special.
\item $-a =\sim a$.
\item $(-a)\mu_a =a$.
\item There is $b\in U^{\#}$ with $(-a)\mu_b =-a \mu_b$.
\item $(-a)\mu_b=-a \mu_b$ for all $b\in U^{\#}$.
\end{enumerate}
\end{lemma}
\begin{proof}
See \cite[3.7]{BG}.
\end{proof}

\bigskip \noindent
Special Moufang sets behave in many respects nicer than non-special ones. They have been 
studied intensively in the recent years and many deep results have been found. Probably 
the most interesting one is the following theorem by Segev (see \cite{S}).
\begin{theorem}\label{Segev} Let $\mathbb{M}(U,\tau)$ be a proper Moufang set with $U$ abelian.
Then $\mathbb{M}(U,\tau)$ is special.
\end{theorem}
It is conjectured that also the converse holds, i.e. that a special 
Moufang set has 
abelian root groups. This conjecture is still open.
\\
We will also make use of the following theorem of Segev and Weiss.
\begin{theorem}\label{invariant subgroups}
Let $M(U,\tau)$ be a special Moufang set, $H$ its Hua subgroup. Suppose 
that $V$ is 
an $H$-invariant subgroup of $U$. If $U$ is not an elementary-abelian $2$-
group, then $V=\{0\}$ or 
$V=U$.
\end{theorem}
\begin{proof}
This is \cite[Theorem 1.2]{SW}.
\end{proof}

\bigskip
For later use, we will need some properties of special elements.
\begin{lemma}\label{special elements}
Suppose that $a \in Z(U)^{\#}$ special and $b \in U^{\#} \setminus \{a,-
a\}$ with $\mu_a 
=\mu_{-a}=\mu_b$. Then
we have
\begin{enumerate}
\item $-a \cdot 3 = b-\sim b +\sim -b$.
\item $b$ is not special. 
\item $o(a)$ divides $6$.
\end{enumerate}
\end{lemma}
\begin{proof}
\begin{enumerate}
\item This is 3.12 (c) of \cite{BG}.
\item This is 3.12 (i) of \cite{BG}.
\item By (a) we have $-a\cdot 3 =b-\sim b +\sim -b=-(-a)\cdot 3 = a \cdot 
3$ and hence $a \cdot 
6 =0$.
Thus the claim follows.
\end{enumerate}
\end{proof}
\section{Rank one groups}
We will assume that $G$ is a rank one group with unipotent subgroups $A$ 
and $B$. As mentioned 
in the introduction, 
this means that $G=\langle A,B\rangle$ and there is a function 
$\mathbf{b}:A^* \to B^*: a \mapsto \fatb(a)$ with $A^{\fatb(a)} =B^a$ for 
all $a \in A^*$.

Like we have said before, if $\mathbb{M}(U,\tau)$ is a Moufang set with 
$U$ nilpotent, then 
$\Gdag$ is a rank one group with unipotent subgroups $U_{\infty}$ and 
$U_0$. For $A\in U^\#$ set $\fatb(\alpha_a) =\alpha_{a\tau^{-1}}^\tau$. 
Conversely, 
if $G=\langle A,B \rangle$ is a rank one group with unipotent subgroups 
$A$ and $B$, 
then $(X,(U_x)_{x\in X})$ is a Moufang set with little projective group 
$G/Z(G)$, where 
$X:=\{A^g|g\in G\}$ and $U_C =CZ(G)/Z(G)$ for $C\in X$.  

Note that not every central extension of $G^{\dagger}$ is a rank 
one group. If for 
example $G= A_5$ and 
$A,B \in \Syl_2(G)$, then $G$ is a rank one group with unipotent 
subgroups $A$ and $B$ and 
$Z(G)=1$, but ${\widehat G} = \langle {\widehat A}, 
{\widehat B} \rangle$ with ${\widehat G}=\SL_2(5)$ and with preimages ${\widehat 
A}, {\widehat B}$ of $A$ 
and 
$B$ in ${\widehat G}$ is not a rank one group with unipotent subgroups 
${\widehat A}$ and ${\widehat 
B}$. 
We refer to \cite{T2} where the notion of a rank one extension is 
introduced. 

Timmesfeld requires in his definition additionally that for all $b\in 
B^*$ there is an element $\fata(b)\in B^*$ with $A^b=B^{\fata(b)}$, but 
this is not necessary as the next lemma shows. In this lemma we introduce the $\mu$-maps for 
rank one groups. These maps will be of great importance for the following.
\begin{lemma}\label{symmetry} Let $G$ be a rank one group with unipotent subgroups 
$A$ and $B$. 
\begin{enumerate}
\item For $a\in A^*$ set $\mu(a)=\fatb(a^{-1})a\fatb(a)^{-1}$. Then we have
$A^{\mu(a)}=B$ and $B^{\mu(a)}A$.
\item $N_A(B)=1=N_B(A)$.
\item The function $\mathbf{b}$ is bijective. If $\mathbf{a}:B^* \to
A^*: b\mapsto \fata(b)$ is its inverse function, then we have
$B^{\fata(b)}=A^b$ for all $b\in B^*$.
\end{enumerate}
\end{lemma}
\begin{proof}
\begin{enumerate}
\item We have $A^{\mu(a)} =A^{\fatb(a^{-1})a\fatb(a)^{-1}}=(B^{a^{-1}})^{a\fatb(a)^{-1}} =
B^{\fatb(a)^{-1}}=B$ and $B^{\mu(a)} =B^{\fatb(a^{-1})a\fatb(a)^{-1}} =B^{a\fatb(a)^{-1}}=
(A^{\fatb(a)})^{\fatb(a)^{-1}}=A$, hence the claim follows.
\item The first equation is \cite[Lemma I (1.2)(3)]{T2}, the second follows immediately by (a).
\item For $e\in A^*$ fixed set $\mu=\mu(e)$. Then for $b\in B^*$ we have
$(A^b)^\mu= (A^\mu)^{b^\mu}=B^{b^\mu}=A^{\fatb(b^\mu)}$ and thus 
$A^b= 
(A^{\fatb(b^\mu)})^{\mu^{-1}}=(A^{\mu^{-1}})^{\fatb(b^\mu)^{\mu^{-1}}}=
B^{\fatb(b^\mu)^{\mu^{-1}}}$. Hence if we define $\fata: B^*\to A^*:
b\mapsto \fatb(b^\mu)^{\mu^{-1}}$, then we have 
$A^b=B^{\fata(b)}$ for all $b\in B^*$. For $b\in B^*$ we have 
$A^{\fatb(\fata(b))}=B^{\fata(b)} =A^b$ and thus $b=\fatb(\fata(b))$ by 
(b). Similarly we obtain $\fata(\fatb(a))=a$ for all $a\in A^*$, 
which shows that $\fatb$ is inverse to $\fata$.
\end{enumerate}
\end{proof}

\bigskip
Note that the formula for $\mu(a)$ 
equals to the formulas for $\mu_a$ in \ref{properties mu map}(a). 

As in the case of Moufang sets, one can define the Hua subgroup. 
\begin{definition}\rm 
 The group $H=\langle \mu(a) \mu(b)| 
a,b \in A^* \rangle$ is called the \textit{Hua subgroup} of $G$. 
\end{definition}
Since $G$ defines a Moufang set, we get 
using \ref{Hua Moufang}.

\begin{prop} $H=N_G(A)\cap N_G(B)$.
\end{prop}

We list some properties for 
the $\mu$-maps.

\begin{lemma}
Let $a \in A^*$. Then:
\begin{enumerate}
\item $\mu(a)$ is the unique element in the double coset $BaB$ satisfying 
$A^{\mu(a)} =B$ and $B^{\mu(a)} =A$.
\item $\mu(a^{-1}) =\mu(a)^{-1}$.
\item $\mu(a^h) =\mu(a)^h$ for all $h\in H$.
\end{enumerate}
\end{lemma}
\begin{proof}
\begin{enumerate}
\item This follows immediately by \ref{symmetry}(b).
\item This follows by (a) since $\mu(a)^{-1}$ interchanges $A$ and $B$ and is contained in the double 
coset $Ba^{-1}B$.
\item This follows again by (a) since $\mu(a)^h$ interchanges $A^h=A$ and $B^h=B$ and is contained in 
the double coset $(BaB)^h= B^h a^h B^h=Ba^hB$.
\end{enumerate}
\end{proof}

\bigskip
Translating the formula $\mu_a =\alpha_a \alpha_{-a\tau^{-1}} \alpha_{-(-a\tau^{-1})\tau}$ 
we get
\begin{lemma}\label{mu maps} Let $G=\langle A,B\rangle$ a rank one group. 
Then for $a\in A^*$ we have 
$\mu(a) = a \fatb(a^{-1}) \fata(\fatb(a^{-1})^{-1})$.
\end{lemma}
 \begin{proof}
 We have
 $$A^{a \fatb(a)^{-1}\fata(\fatb(a)^{-1})^{-1}} =A^{\fatb(a)^{-1} \fata(\fatb(a)^{-1})^{-1}} =B$$
 and
 $$B^{a \fatb(a)^{-1} \fata(\fatb(a)^{-1} )^{-1}} =A^{\fata(\fatb(a)^{-1})^{-1} }=A.$$
 Moreover,
 \begin{align*}
  a\fatb(a)^{-1} \fata(\fatb(a)^{-1})^{-1} =\fatb(a^{-1})^{-1} \mu(a) \fata(\fatb(a)^{-1})^{-1} =\\
  \fatb(a^{-1})^{-1} ( \fata(\fatb(a)^{-1})^{-1})^{\mu(a)^{-1}}\mu(a) \in B\mu(a) \subseteq BaB.
 \end{align*}
 Since $\mu(a)$ is the unique element in $BaB$ which interchanges $A$ and $B$, the claim 
 follows.
 \end{proof}
 
\bigskip 
For $a \in A^*$ and $\tau \in N_G(\{A,B\})$ 
let $a\diamond \tau$ be the unique element in $A^*$ such that $B^{a\diamond \tau} =
(B^a)^{\tau}$. Then $a\mapsto a\diamond \tau$ is a permutation of $A^*$. 
If $h\in N_G(A)\cap N_G(B)$, then we have $a^h =a\diamond  h$ for all $a\in A^*$. 
\begin{lemma}\label{b(a)} Let $a\in A^*$. Then we have 
$\fatb(a)=(a\diamond \tau^{-1})^{\tau}$ for all $\tau \in 
N_G(\{A,B\})^o:=N_G(\{A,B\}) \setminus N_G(A)\cap N_G(B)$.
\end{lemma}
\begin{proof}
We have $(a\diamond \tau^{-1})^\tau \in A^\tau =B$ and $A^{(a\diamond  
\tau^{-1})^\tau}=A^{\tau^{-1}(a\diamond \tau^{-1}) \tau}=
(B^{a\diamond \tau^{-1}})^\tau = (B^{a\tau^{-1}})^\tau =B^a$. 
Since $\fatb(a)$ is the unique element in $B$ with $B^a=A^{\fatb(a)}$, 
the claim follows.  
\end{proof}

For $a\in A^*$ we set $a^{\sim}: = \fata(\fatb(a)^{-1})$, $a^{-\sim}:=
(a^{-1})^\sim$ and $a^{\sim -} =(a^\sim)^{-1}$. Note that by \ref{b(a)}
this corresponds to the definition of $\sim a$ in Moufang sets and that 
$a^{\sim}  = (a\diamond \tau^{-1})\diamond \tau$ for 
all $\tau \in N_G(\{A,B\})^o$. Moreover, by \ref{mu maps} we have
$\mu(a) =a \fatb(a)^{-1} a^{\sim -}$. 
Furthermore, we have $a^{\sim \sim} =a$, $(a^{\sim -})^\sim =
(a^\sim)^{-\sim}$ and $(a^{-\sim})^{-1}=(a^{-1})^{\sim -}$ for all $a\in 
A^*$, so an expression in $-$ and $\sim$ is well-defined also for a length 
greater than two.

One easily sees that the formulas in \ref{formula sim} also hold for rank one groups.
\begin{lemma}\label{mu a tau} Let $a,b\in A^*$.
\begin{enumerate} 
\item 
$\mu(a\diamond  \mu(b)) =\mu(b)^{-1} \mu(a)^{-1} \mu(b)$ for all $a,b \in A^*$.
\item For a fixed $e\in A^*$ set $h(a):=\mu(e)^{-1}\mu(a)$. Then we have
$h(a\diamond  \mu(b)) =h(b^{-1})h(a)^{-1}h(b)$. Especially we have
$h(a\diamond  \mu(e))=\mu(e)^{-2} h(a)^{-1}$.
\end{enumerate}
\end{lemma}
\begin{proof}
\begin{enumerate}
\item
We have
\begin{align*}
\mu(a\diamond \mu(b)) =\fatb((a\diamond \mu(b))^{-1}) (a\diamond \mu(b)) \fatb(a\diamond \mu(b))^{-1} = \\
((a\diamond \mu(b))^{-1}\diamond \mu(b^{-1}))^{\mu(b)} 
(a\diamond \mu(b)) a^{-\mu(b)} =  ((a^{\sim}) (a\diamond  \mu(b))^{\mu(b^{-1})} a^{-1})^{\mu(b)} =\\
(a^{\sim} \fatb(a) a^{-1})^{\mu(b)} = (a \fatb(a)^{-1} a^{\sim-})^{-\mu(b)} =
\mu(a)^{-\mu(b)} = \mu(b)^{-1} \mu(a)^{-1} \mu(b).  
\end{align*}
\item 
We have by (a)
\begin{align*}
h(a\diamond \mu(b)) =\mu(e)^{-1} \mu(a\diamond \mu(b)) =\mu(e)^{-1} \mu(b)^{-1}
\mu(a)^{-1}\mu(b) =\\ \mu(e)^{-1}\mu(b^{-1}) \mu(a)^{-1}\mu(e)\mu(e)^{-1}
\mu(a)=h(b^{-1})h(a)^{-1} h(b).
\end{align*}
\end{enumerate}
\end{proof} 

 \begin{lemma}\label{formulas} Let $a,b \in A^*$ with $a \ne b$ and 
 $\tau \in N_G(\{A,B\})$. Then we have
 \begin{enumerate}
 \item $(a\diamond \tau^{-1}(b\diamond \tau^{-1})^{-1})\tau =(ab^{-1})\diamond  \mu(b) b^\sim$.
 \item $\mu((a\diamond \tau^{-1}(b\diamond \tau^{-1})^{-1})\diamond \tau) =\mu(b)^{-1}
 \mu(ba^{-1}) \mu(a)$.
 \end{enumerate}
 \end{lemma} 
 \begin{proof}
 This is Theorem 1.1 of \cite{DS} translated into the language of rank $1$ groups. 
 Since the unipotent groups of a rank one 
 group are isomorphic to the root groups of the corresponding Moufang set, part (a) follows 
 immediately. Part (b) holds a priori only in $G/Z(G)$, but one can adapt the proof of the 
 theorem above for rank one groups. Alternatively, one can apply 4.8 of \cite{L}.
\end{proof}
\begin{lemma}\label{h+} For $a,b\in A^*$ with $a\ne b$ and $\mu=\mu(e)$ we have
$h(b\diamond\mu(a\diamond \mu)^{-1})=h(b)^{-1} h(ab^{-1})h(a)^{\mu}$.
\end{lemma}
\begin{proof}
  We have by \ref{mu a tau}(a) and \ref{formulas} for $\tau=\mu^{-1}$:
$$\mu \mu((b\diamond \mu)(a\diamond \mu)^{-1})\mu^{-1} =
\mu(((a\diamond \mu)(b\diamond \mu)^{-1})\diamond \mu^{-1}) =\mu(b)^{-1}
 \mu(ba^{-1}) \mu(a)$$  
 and hence
\begin{align*}
h((b\diamond \mu) (a\diamond \mu)^{-1})=\mu(b)^{-1}\mu(ba^{-1})\mu(a)\mu=\\
\mu(b)^{-1} \mu^{-1} \mu \mu(ba^{-1}) (\mu\mu(a))^\mu =
h(b)^{-1} h(ab^{-1})h(a)^{\mu}.
\end{align*} 
\end{proof}
\bigskip  
  
We will frequently make use of the following observation:
\begin{prop}\label{faithful} Let $G$ be a rank one group with unipotent subgroups $A$ and $B$. 
\begin{enumerate}
\item If $N$ is a normal subgroup of $G$, then either $G =NA$ and thus $G/N$ is nilpotent, 
or $N \leq Z(G)$ and $A \cap N =1$.
\item $Z_2(G)=Z(G)$.
\item If $N$ is normal subgroup of $G$ with $G/N$ nilpotent, then $G =NA$.
\item If $V$ is a $G$-module, then either $A$ acts faithfully on $V$ or $G$ and $A$ have 
the same image in $\GL(V)$.
\end{enumerate}
\end{prop}
\begin{proof}
\begin{enumerate}
\item This is \cite[I (1.10)]{T1}.
\item This is \cite[I (2.1)]{T1}.
\item If $G \ne NA$, then $N \leq Z(G)$ and thus $G$ is nilpotent, which 
contradicts (b). 
\item By (a), either $A \cap C_G(V)\leq Z(G)$ or $G =C_G(V)A$. Thus the claim follows.
\end{enumerate}
\end{proof}

\bigskip
 As for Moufang set, one can define special elements in rank one groups. If $G=\langle A,B 
\rangle $ is a rank one group, then $a\in A^*$ is called special if $\fatb(a^{-1}) =\fatb(a)^{-1}$. 
This holds iff there is $\tau \in N_G(\{A,B\})^o$ with $(a^{-1})\diamond \tau=(a\diamond \tau)^{-1}$ 
for all $a\in A^*$; in this case, this formula holds for all $\rho \in N_G(\{A,B\})^o$.  
The rank one group $G$ is called special if all elements of $A^*$ special. This corresponds to the 
definition of special elements in Moufang sets. \\

In order to examine subgroups of rank one groups which are again rank one groups, we need
the following definition:
\begin{definition}\rm A subgroup $A_0$ of $A$ is called a \textit{root subgroup} if the set $B_0 :=\{1\} 
\cup \{\fatb(a)| a \in A_0^*\}$ 
is a subgroup of $B$.
\end{definition}
If $A_0$ is a root subgroup of $A$, then one easily sees that $G_0:=\langle A_0,B_0 \rangle $ 
is rank one group with unipotent 
subgroups $A_0$ and $B_0$. We say that $A_0$ is a special root subgroup if the rank one group 
$G_0$ is special.
\begin{lemma}\label{S3} Let $a\in A$ be a special involution. Then $\fatb(a)$ is also 
an involution and $\langle a,\fatb(a)\rangle \cong S_3$.
\end{lemma}
\begin{proof}
Since $a$ is special, we have $\fatb(a) =\fatb(a^{-1})=\fatb(a)^{-1}$, hence $b=\fatb(a)$ is 
an involution. We have $aba =\mu(a)=bab $ by \ref{mu maps}, which shows that $ab$ has order 
$3$. Thus the claim follows.
\end{proof}
\begin{theorem}\label{root subgroups}
Let $G$ be a special rank one group with unipotent subgroups $A$ and $B$. Suppose that 
$A$ is not an elementary-abelian $2$-group and that $A_0$ is an $H$-invariant subgroup of $A$.
Then $A_0=1$ or $A_0=A$.
\end{theorem}
\begin{proof} 
This  follows easily from \ref{invariant subgroups}.
\end{proof}
\section{Some ring theory}
A ring $R$ is called \textit{semi-prime} if $R$ has no nilpotent ideals. This means that if 
$I$ is 
an ideal 
of $R$ with $I^n =0$ for a natural number $n \geq 1$, then $I=0$. We set $\B(R):=\bigcap \{I| 
I 
\trianglelefteq 
R, R/I$ semi-prime$\}$. One easily sees that $R/\B(R)$ is semi-prime, so $\B(R)$ is 
the smallest ideal of $R$ such that the factor ring is semi-prime. $\B(R)$ is called the 
(lower) 
\textit{Baer radical} of $R$ (see \cite{B}). It is contained in the Jacobson radical of $R$ 
and 
every element of 
$\B(R)$ is nilpotent, so an element $x \in R$ is a unit iff $x+\B(R)$ is a unit in 
$R/\B(R)$. 
\newline The following lemma is a generalisation of \cite[Lemma 14.1.1]{Sc}, where both $S$ 
and $R$ 
are skew fields (and so $I_{x+1} =0$). 
\begin{lemma}\label{subrings} Let $R$ be a subring of ring $S$ such that $1\in R$.
Suppose there is a unit $x \in S$ such that $x+1$ is again a unit and such that 
$x^{-1} R x = (x+1)^{-1} R (x+1) = R$. Then $(x+1)^{-1} \in R$ or 
$I_{x+1}:=R \cap (x+1)R$ is an ideal of $R$ with $I_{x+1} \ne R$ and 
$x^{-1} u x - u \in I_{x+1}$ for all $u \in R$ and thus $x$ induces the trivial automorphism 
on $R/I_{x+1}$. 
\end{lemma}
\begin{proof} Let $u \in R$, set $v := x^{-1} u x$ and $w := (x+1)^{-1} u (x+1)$. Thus 
$xv = ux$ and $(x+1)w = u(x+1)$. We get $$(x+1)v - v+u = xv +u=ux+u =u(x+1) =(x+1)w$$ and hence
$v-u=(x+1)(v-w)  \in R \cap (x+1) R = I_{x+1}$. Since 
$$I_{x+1} =(x+1)R \cap R = (x+1) (x+1)^{-1} R (x+1) \cap R = R(x+1) \cap R,$$
it follows that $I_{x+1}$ is an ideal of $R$. Since $I_{x+1}=R$ iff 
 $1 \in (x+1)R$ iff $(x+1)^{-1} \in R$, the claim follows. \end{proof}
\section{Quadratic Jordan algebras}
\begin{definition}\rm Let $K$ be a field, $J$ a $K$-vector space, $e \in J^{\#}=J
\setminus \{0\}$ and let 
$Q: J \to \End_K(J):a\mapsto Q_a$ be a quadratic map, by which we mean that 
the map $Q_{.,.}: J \times J \to \End_K(J): (a,b)\mapsto Q_{a,b}:=Q_{a+b}-Q_a-Q_b$ is 
bilinear. 
For $a,b \in J$ define $V_{a,b}\in \End_K(J)$ by $cV_{a,b}:=bQ_{a,c}$. Then $(J,Q,e)$ is called a 
\textit{weak quadratic 
Jordan algebra} if the following identities hold for all $a,b \in J$.
\begin{itemize}
\item[(QJ1)] $Q_e=id_J$.
\item[(QJ2)] $Q_a V_{a,b}=V_{b,a}Q_a$.
\item[(QJ3)] $Q_{bQ_a}=Q_a Q_b Q_a$.
\end{itemize}
A weak quadratic Jordan algebra is called a \textit{quadratic Jordan algebra} if the identities (QJ1)-(QJ3) hold
strictly, that means that they continue to hold in $L \otimes_K J$ for every extension field $L$ of $K$.
\end{definition}
\begin{remark}\rm
\begin{enumerate}
\item (QJ1)-(QJ3) hold strictly iff their linearised versions holds. 
Thus (QJ2) and (QJ3) hold strictly iff
\begin{itemize}
\item[(QJ2*)]
\begin{enumerate}
\item $Q_{a_1} V_{a_1,b}=V_{b,a_1}Q_{a_1}$.
\item $Q_{a_1}V_{a_2,b} +Q_{a_1,a_2}V_{a_1,b}=V_{b,a_1}Q_{a_1,a_2}+V_{b,a_2}Q_{a_1}$. 
\end{enumerate}
\item[(QJ3*)]
\begin{enumerate}
\item $Q_{bQ_{a_1}}=Q_{a_1} Q_{b}Q_{a_1}$.
\item $Q_{bQ_{a_1},bQ_{a_1,a_2}}=Q_{a_1,a_2}Q_{b}Q_{a_1}+Q_{a_1}Q_{b}Q_{a_1,a_2}$.
\end{enumerate}
\end{itemize}
hold for all $a_1,a_2,b \in J$ hold.
\item The author does not know an example of a weak quadratic Jordan algebra which is not a 
quadratic 
Jordan algebra.
\item One can show that (QJ1-(QJ3) hold strictly iff they continue to hold in 
$K[t]\otimes_K J$.  This is automatically the case if $K$ has at least 
$4$ elements. 
\end{enumerate}
\end{remark}
\begin{example}\rm 
Let $R$ be a unital, associative algebra over $K$. For $a \in R$ define $Q_a: R \to R: b\mapsto 
bab$. Then
$R^+:=(R,Q,1)$ is a quadratic Jordan algebra. 
\end{example}
\begin{definition}\rm
\begin{enumerate}
\item If $(J,Q,e)$ and $(J^{\prime},Q^{\prime},e^{\prime})$ are weak quadratic Jordan algebras 
over $K$, 
then a \textit{Jordan homomorphism} between $J$ and $J^{\prime}$ is a homomorphism 
$f: J \to J^{\prime}$ such that $f(e)=e^{\prime}$ and 
$f(aQ_b)=f(a)Q^{\prime}_{f(b)}$ for all $a,b \in J$ holds.
\item Let $J^{\prime}$ be a subspace of a weak quadratic Jordan algebra $J$. Then $J^{\prime}$ 
is called
a \textit{Jordan subalgebra} of $J$ if $e \in J$ and if $J^{\prime} Q_a \subseteq J^{\prime}$ for all
$a \in J^{\prime}$ holds.
\item A quadratic Jordan algebra $J$ is called \textit{special} if there is an associative $K$-algebra 
$R$ such that $J$ is isomorphic to a Jordan subalgebra of $R^+$ and \textit{exceptional} else. 
\end{enumerate}
\end{definition}
\begin{definition}\rm
Let $J$ be a weak quadratic Jordan algebra.
\begin{enumerate}
\item An element $a \in J^{\#}$ is called \textit{invertible} if there is an element $b \in J$ with
$bQ_a =a$ and $e Q_b Q_a =e$, or equivalently, if $Q_a$ is invertible. The element 
$b=aQ_a^{-1}$ 
is uniquely determined and denoted by $a^{-1}$. 
\item $J$ is called a \textit{weak quadratic Jordan division algebra} if every $a \in J^{\#}$ is 
invertible.
\end{enumerate}
\end{definition}
If $J=R^+$ for a ring $R$, then one easily sees that $a\in J$ is invertible iff $a\in R^*$. In this 
case, the inverses coincide.\\
In \cite{G2} it was proved:
\begin{theorem}\label{QJDA} A weak quadratic Jordan division algebra is a Jordan division algebra.
\end{theorem}
\begin{definition}\rm Let $(J,Q,e)$ be a weak quadratic Jordan algebra. Then 
the subgroup of $\GL_K(J)$ generated by all $Q_a$ with $a\in J^{\#}$ invertible is called 
the \textit{inner structure group} of $J$.
\end{definition}
\begin{example}\rm
If $K$ is a field, then the inner structure group of $K^+$ is just $(K^*)^2$, the group of 
non-zero squares of $K^*$. 
\end{example}  
For special quadratic Jordan division algebras there is an alternative description.
\begin{lemma}
Let $R$ be a unital, associative $K$-algebra and $J \subseteq R$ with $1 \in J$. Suppose that
every $a \in J^{\#}$ is invertible in $R$ and that $a^{-1}$ is again in $J$. Then $J$ is a 
special 
quadratic Jordan division algebra.
\end{lemma}
\begin{proof}
By the Hua identity (\cite{H}) we have
$$aba =a-(a^{-1} -(a-b^{-1})^{-1})^{-1}$$
for all $a,b \in J$ with $b \ne a^{-1}$, 
which shows that $J$ is a Jordan subalgebra of $R^+$. Since $a=aba$ and $1=ab^2a$ holds for 
$b=a^{-1}$ for all $a \in J^{\#}$, every $a\in J^\#$ is invertible with inverse $a^{-1}$. 
\end{proof}

\begin{theorem}
Let $J$ be a quadratic Jordan division algebra and $X=J \dot{\cup} \{\infty\}$. 
Define $\tau \in \Sym X$ by $\infty \tau =0$, $0\tau =\infty$ and  
$a \tau =-a^{-1}$. Then $\mathbb{M}(J) :=\mathbb{M}(J,\tau)$ is a Moufang set.
\end{theorem}
\begin{proof}
See \cite[Theorem 4.2]{DW}.
\end{proof}
\bigskip

For a special quadratic Jordan division algebra there is a more direct way to construct a rank 
$1$ group.
Let $R$ be a unital, associative $K$-algebra and $J \subseteq R^+$ a quadratic Jordan division 
algebra. 
For $a \in J$ set 
$$\alpha_a:= \left( \begin{array}{cc} 1 & 0 \\ a & 1 \end{array} \right) \
\text{and} \ \beta_a:= \left( \begin{array}{cc} 1 & a \\ 0 & 1 \end{array} \right).$$
Define $A:=\{\alpha_a| a\in J\}$ and $B:=\{\beta_a| a\in J\}$. Then $A$ and $B$ are subgroups 
of 
$\SL_2(R)$ and $B^{\alpha_a}=A^{\beta_{a^{-1}}}$, so $\SL_2(J,R):=\langle A, B \rangle $ is a rank 
$1$ group 
with unipotent subgroups $A$ and $B$. Note that $\SL_2(J,R)$ may depend on $R$, but 
$\SL_2(J,R)/Z(\SL_2(J,R))$ is the little projective group for $\mathbb{M}(J)$ and thus 
independent of 
the choice of $R$.
\\
It is conjectured that every special Moufang set with abelian root groups is isomorphic to
$\mathbb{M}(J)$ for a quadratic Jordan division algebra $J$. More generally, let 
$M(U,\tau)$ be a special Moufang set with $U$ abelian and $\tau =\mu_e$ for some $e \in U^{\#}$. 
Since $U$ is either an elementary-abelian $p$-group for some prime $p$ or torsion-free and 
uniquely 
divisible, $U$ can be considered as a $\mathbb{F}$-vector space for $\mathbb{F}=\mathbb{F}_p$ in 
the first 
case and $\mathbb{F}=\mathbb{Q}$ in the second case. Let $\mathfrak{h}: U \to \End_{\mathbb{F}}
(U): a 
\mapsto h_a$ with the convention $h_0=0$. Then we have
\begin{conjecture}
$(U,\mathfrak{h},e)$ is a quadratic Jordan division algebra over $\mathbb{F}$.
\end{conjecture}
There has been partial progress to prove this conjecture (see \cite{DMS}). One sees easily that
(QJ1) and (QJ3) hold. Also, one sees that if $U$ is a (weak) 
quadratic Jordan algebra, then every
$a \in U^{\#}$ is invertible with inverse $-a\tau$. 
It remains to show that $\mathfrak{h}$ is bilinear and that (QJ2) holds. 
\begin{theorem}\label{QJ2}
If (QJ2) holds for $\mathbb{M}(U,\tau)$, then $\mathfrak{h}$ is bilinear, therefore 
$\mathbb{M}(U,\tau)$ is the Moufang set for a quadratic Jordan division algebra. 
\end{theorem}
\begin{proof}
See \cite[Theorem 5.4]{DS} and also \cite{G2} for the case $\characteristic \FF \in \{2,3\}$.
\end{proof}

\bigskip
One can show that for $\characteristic \FF
\ne 2,3$ the bilinearity of $\mathfrak{h}$ and (QJ2) are equivalent (see \cite{DS}).\\
(QJ2) holds by definition exactly iff $$ah_{b,c}h_a = bh_{a,ch_a}$$ 
holds for all $a,b,c \in U$. Actually, it is enough to prove a weaker version of this formula.
(QJ1)-(QJ3) depend on the choice of $e$. If $e^{\prime} \in U^{\#}$, 
$h^{\prime}_a:=\mu_{e^{\prime}}
\mu_a$ for $a \in U^{\#}$ and $\mathfrak{h}^{\prime}: U \to \End_{\FF}(U): a \mapsto 
h^{\prime}_a$, 
then $(U, \mathfrak{h}^{\prime},e^{\prime})$ is called an \textit{isotope} of $(U,\mathfrak{h},e)$. 
One easily sees that $(U,\mathfrak{h},e)$ is a (weak) quadratic Jordan division algebra iff
$(U, \mathfrak{h}^{\prime},e^{\prime})$ is.
\begin{lemma}\label{weak QJ2}
Suppose that
$$\hbox{(QJ2)'} \hspace{1cm} eh_{b,c} =bh_{c,e}$$
holds for all $b,c \in U$ in every isotope of $U$, then 
(QJ2) holds in $U$.
\end{lemma}

We will prove some identities for quadratic Jordan algebras which we will need in chapter 4.

\begin{lemma}\label{identities} Let $J$ be a quadratic Jordan algebra and $a\in J^\#$ invertible. 
Then we have $V_{y,a^{-1}} =V_{a,yQ_a^{-1}}$ and $V_{a^{-1},y} =V_{yQ_a^{-1},a}$ for all $y\in J$.
\end{lemma}
\begin{proof}
We have again by (QJ2)
$$ aQ_{y,x}Q_a = xV_{y,a}Q_a = xQ_a V_{a,y} = yQ_{xQ_a,a}$$ for all $x\in J$.  
Replacing $x$ by $xQ_a^{-1}$, $y$ by $yQ_a^{-1}$ and applying (QJ3) yields
\begin{equation*}
x V_{y,a^{-1}} = a^{-1} Q_{y,x} = a Q_{yQ_a^{-1},xQ_a^{-1}} Q_a =yQ_a^{-1} Q_{a,x} =xV_{a,yQ_a^{-1}},
\end{equation*}
which yields $V_{y,a^{-1}} =V_{a,yQ_a^{-1}}$. 
Moreover, we have again by (QJ2)
$$xV_{yQ_a^{-1},a} = yQ_a^{-1} V_{x,a}=yV_{a,x}Q_a^{-1} = xQ_{a,y}Q_a^{-1}$$
and thus \begin{equation}\label{equ1} V_{yQ_a^{-1},a} = Q_{a,y}Q_a^{-1}. \end{equation}
We also have \begin{equation}
aQ_{y,x}Q_a =xV_{y,a}Q_a =xQ_a V_{a,y} = yQ_{xQ_a,a} \end{equation} 
Since the left side of (2) is symmetric in $x$ and $y$, the right side of (2) is also symmetric in 
$x$ and $y$. Therefore we have
\begin{equation}
yQ_{xQ_a,a} =xQ_{yQ_a,a}.
\end{equation}
Replacing $y$ by $yQ_a^{-1}$ and applying (QJ3) yields
$$xV_{a^{-1},y} = yQ_{a^{-1},x} = yQ_a^{-1} Q_{xQ_a,a}Q_a^{-1}=xQ_{y,a}Q_a^{-1}.$$
Thus by (\ref{equ1}) we get $V_{a^{-1},y} = Q_{y,a}Q_a^{-1} =V_{yQ_a^{-1},a}$, as desired.
\end{proof}
\begin{lemma}\label{commutator}
Let $J$ be a special quadratic Jordan algebra inside a ring $R$. 
Then for all $x,y,z\in R$ we have 
$[x,y]^2\in J$ and $[x,y,z]=[[x,y],z]\in J$. 
\end{lemma}
\begin{proof}
We have $$[x,y]^2 =(xy-yx)^2 = x(xy+yx)y+y(xy+yx)x-x^2y^2-y^2 x^2 -xy^2 x-yx^2y\in J$$ and 
$$[x,y,z] =(xy-yx)z-z(xy-yx) =xyz+zyx-yxz-zxy\in J.$$
\end{proof}
\section{Envelopes of special quadratic Jordan algebras}
If $J$ is a special quadratic Jordan algebra over $K$,  
$R$ is a $K$-algebra and $f:J \to R^+$ is an injective Jordan homomorphism such that $f(J)$ 
generates $R$ as 
a ring, then $(R,f)$ is called an \textit{envelope} for $J$. An envelope $(R,f)$ of $J$ is called 
\textit{universal} if for any other envelope 
$(S,g)$ of $J$ there is a $K$-algebra 
homomorphism $\varphi:R \to S$ with $ \varphi \circ f = g$. One can construct a universal 
envelope as follows:
Let $T(J) = K \oplus J \oplus (J \otimes_K J) \oplus (J \otimes_K J \otimes_K J) \ldots$ be the 
tensor algebra over 
$J$ and let $I(J)$ be the ideal generated by all elements of the form $bQ_a- a \otimes b 
\otimes a$ and 
by $1_K -1_J$. Then $T(J)/I(J)$ together with the canonical embedding $a \mapsto a + I(J)$ is a 
universal envelope 
of $J$. This construction depends on the choice of $K$. If $F$ is a subfield of $K$, then $J$ 
is also a Jordan 
algebra over $F$, and it is not clear if the universal envelope over $F$ equals the universal 
envelope over $F$. 
\\
If $U$ is the universal envelope of $J$, then we call $U/\B(U)$ the universal semi-prime 
envelope of $J$. 
If $J$ is a Jordan division algebra, then this definition does not depend on the ground field 
$K$. If $F$ is a subfield of $K$, $U_F$ the 
universal envelope over $F$ and $U_K$ the universal envelope over $K$, then there is an 
epimorphism from $U_F$ to $U_K$, and one can 
easily see that the kernel of this map is contained in $\B(U_F)$.
\begin{prop}\label{semi-prime} Let $J$ be a special quadratic Jordan division algebra. Then the 
following statements are equivalent.
\begin{enumerate}
\item There is a commutative field $F$ such that $J =F^+$ or $\characteristic J =2$ and 
$J$ is a $F^2$-vector space with $F^2 \subseteq J \subseteq F$.
\item 
The universal semi-prime envelope of $J$ is a commutative field.
\item There is a commutative envelope $R$ for $J$.
\item The inner structure group of $J$ is abelian.
\end{enumerate}
\end{prop}
\begin{proof}
If $J = F^+$ for a commutative field $F$, then $\dim_F J = 1$ and so one easily sees that $F$ 
is the universal $F$-envelope of $J$. 
Suppose $F^2 \subseteq J \subseteq F$ for a commutative field $F$ with $\characteristic F =2$ 
and let $U=T(J)/I(J)$ the universal $F^2$-envelope of $J$.
Then for all $a,b \in J$ we have $a^2 - a \otimes a \in I(J)$ and thus   
\begin{align*} (a+b)^2 -(a+b) \otimes (a+b) -(a^2 -a \otimes a) -(b^2 -b \otimes b) = \\
a^2 +b^2 +a\otimes a +b \otimes b +a\otimes b +b\otimes a +a^2 
+a\otimes a +b^2 +b\otimes b\\=a \otimes b + b \otimes a \in I(J).\end{align*}
This shows that $U$ is commutative 
Since $F$ is an envelope for $J$, there is a homomorphism $\varphi:U \to F$ with $\varphi(a) 
=a$ for all $a \in J$. The kernel 
of $\varphi$ is generated by all elements $a b - a \otimes b+I(J)$ with $a,b,ab \in J$. Since 
$a \otimes b +b \otimes a \in I(J)$, we get \begin{align*}
(ab-a\otimes b)^2 +I(J) = a^2 b^2 + a \otimes b \otimes a \otimes b + I(J) =\\ a^2 b^2 + 
a\otimes a \otimes b \otimes b + I(J) =
a^2 b^2 +a^2 b^2 =0.\end{align*} This shows that $\ker \varphi$ is a nilpotent ideal and so the 
semiprime envelope of $J$ is just $F$. Thus we have shown 
that (a) implies (b). 
\newline
The implications (b) to (c) and (c) to (d) are obvious. For the remaining implication (d) to 
(a) we use the fact that every Jordan division algebra $J$ defines a special Moufang set 
$\mathbb{M}(J)$ and the Hua 
group of $\mathbb{M}(J)$ is just the inner structure group of $J$ (\cite{DW}, 4.1 and 4.2). 
Thus (d) to (a) is a consequence of 6.1 in \cite{DW} and the main theorem of \cite{G}.
   \end{proof}
   
\bigskip \noindent   
Alternatively, this proposition follows by the classification of quadratic Jordan division 
algebras (\cite[15.7]{McZ}).
\section{Involutory sets and pseudo-quadratic forms}   
If $R$ is a ring with involution $*$, then we set $\mathcal{H}(R,*):=\{x \in R| x^* =x\}$. 
\begin{definition}\rm
An additive subgroup $\mathcal{H}_0(R,*)$ of $\mathcal{H}(R,*)$ 
is called an \textit{ample Hermitian Jordan algebra} or an \textit{ample subspace} of 
$\mathcal{H}(R,*)$ if $1 \in 
\mathcal{H}_0(R,*)$ and $x^* \mathcal{H}_0(R,*) x \subseteq \mathcal{H}_0(R,*)$ for all $x \in K$. 
\end{definition}
If every non-zero element of $\mathcal{H}_0(R,*)$ is invertible, then 
$\mathcal{H}_0(R,*)$ is a special quadratic Jordan division algebra since $x^{-1} =x^{-*} = x^{-*} x 
x^{-1} \in 
\mathcal{H}_0(R,*)$ for all $x \in \mathcal{H}_0(R,*)^{\#}$. Note that $x+x^* = (x+1)^* (x+1)-x^* x -1 \in 
\mathcal{H}_0(R,*)$ for all $x \in R$. 
Set $R_* :=\langle x^* x| x \in R\rangle$. Then $R_*$ is the smallest Hermitian Jordan algebra relative 
to $*$. 
If $2$ is invertible in $R$, 
then $x = \frac{x}{2} + \frac{x}{2}^* \in \mathcal{H}_0(R,*)$ for all $x 
\in \mathcal{H}(R,*)$ and so 
$K_* = \mathcal{H}_0(R,*) =\mathcal{H}(R,*)$. If $2$ is not invertible in $R$, 
then $R_* \subsetneq \mathcal{H}_0(R,*) \subsetneq \mathcal{H}(R,*)$ 
is possible. 
\begin{lemma}\label{amplifier}
Let $R$ be a ring with involution $*$ and $J$ an additive subgroup of $
\mathcal{H}(R,\star)$ such that $J$ generates $R$ as a ring,
$r+r^\star\in J$ for all $r\in R$ and $aba\in J$ for all $a,b\in J$. Then $J$ is an ample 
subspace of $\mathcal{H}(R,*)$. 
\end{lemma}
\begin{proof}
Set $N_R(J):=\{x\in R|x^*ax \in J\ \forall a\in J\}$. For $x,y\in N_R(J)$ and $a\in J$, we 
have $$(xy)^* a xy = y^* \underbrace{x^* ax}_{\in J}y \in J$$ and 
$$(x+y)^* a (x+y)= x^*ax +y^* ay +\underbrace{x^* ay +(x^*ay)^*}_{\in J} \in J.$$ 
Since $J\subseteq N_R(J)$ and $J$ generates $R$ as a ring, we conclude $N_R(J)=R$, thus the 
claim follows.
\end{proof}

If $K$ is a skew field, then one has either $\mathcal{H}_0(K,*) \subseteq \mathrm{Cent} K$ or $\langle \mathcal{H}_0(K,*) \rangle =K$ 
as a ring (\cite{TW}, (23.23)).  
If $K$ is a skew field with involution $*$ and $K_0$ an ample hermitian Jordan algebra, then one 
calls 
$(K,K_0,*)$ an \textit{involutory set}. 
\begin{definition}\rm
\begin{enumerate}
\item
 If $L_0$ is a right 
 $K$-vector space, then a map $\pi:L_0 \to K/K_0$ is called a \textit{pseudo-quadratic form} 
relative to $(K,K_0,*)$ 
if $\pi(a \lambda ) = \lambda^ * \pi(a) \lambda $ for all $a \in L_0, \lambda \in K$ and if 
there is 
a skew-hermitian form $f: L_0 \times L_0 \to K$ relative to $*$ with $\pi(a+b) \equiv \pi(a) + 
\pi(b) + 
f(a,b) \ \mod \ K_0$ for all $a,b \in L_0$. 
\item A subspace $X_0$ of $L_0$ is called \textit{isotropic} if 
$\pi(a) =0$ for all $a \in X_0$.
\item The maximal dimension of an isotropic subspace is called the \textit{Witt index} of 
$\pi$.
\item If the Witt index is $0$, then $\pi$ is called \textit{anisotropic}. 
\end{enumerate}\end{definition}
\section{Quadrangular algebras and Moufang quadrangles}
Quadrangular algebras were introduced in \cite{W}. They 
are the algebraic structures that describe certain Moufang quadrangles. 
The prototype of a quadrangular algebra is an anisotropic pseudo-quadratic space over a 
field or a quaternion algebra. By omitting the multiplication of the field resp. quaternion 
algebra and keeping all other axioms one gets an algebraic structures that enables a unique 
treatment of pseudo-quadratic Moufang quadrangles over a field or quaternion algebra and 
exceptional Moufang quadrangles of type $\mathrm{E}_6$, $\mathrm{E}_7$ and $\mathrm{E}_8$.   
\begin{definition}\rm A \textit{pointed quadratic space} is a $4$-tuple $(K,L,q,1)$ such that 
$K$ is a field, $L$ is a $K$-vector space, $q:L \to K$ is a quadratic form and 
$1 \in L$ is an element with $q(1)=1$.
\end{definition}
\begin{definition}\rm Let $(K,L,q,1)$ be a pointed quadratic space and let $f$ be the 
bilinear form associated to $q$. Then the linear map
$$\sigma: L \to L: x \mapsto f(x,1)1 -x$$
is called the \textit{standard involution} associated to $(K,L,q,1)$.
\end{definition}
\begin{definition}\rm
\begin{enumerate}
\item
Let $(K,L,q,1)$ be a pointed quadratic space with associated bilinear 
form $f$ and standard involution $\sigma $, let $X$ be a $K$-vector space and let 
\begin{align*} \cdot: X \times L \to L, (a,v)\mapsto av,\\
h: X\times X \to L \hbox{ and }
\theta: X \times L \to L\\
\end{align*}
be maps such that the following holds:
\begin{itemize}
\item[(A1)] The map $\cdot$ is $K$-bilinear
\item[(A2)] $a\cdot 1 =a$ for all $a \in X$.
\item[(A3)] $(av)v^{\sigma} =q(v) a $ for all $a \in X$ and all $v \in L$.
\item[]
\item[(B1)] $h$ is $K$-bilinear.
\item[(B2)] $h(a,bv) =h(b,av)+f(h(a,b),1)v$ for all $a,b \in X$ and all $v \in L$.
\item[(B3)] $f(h(av,b),1)=f(h(a,b),v)$ for all $a,b \in X$ and all $v \in L$.
\item[]
\item[(C1)] For all $a \in X$, the map $v \mapsto \theta(a,v)$ is $K$-linear.
\item[(C2)] $\theta(ta,v)=t^2 \theta(a,v)$ for all $a \in X, v\in V$ and $t \in K$.
\item[(C3)] There is a function $g:X\times X \to K$ such that 
$$\theta(a+b,v)=\theta(a,v)+\theta(b,v)+ h(a,bv)-g(a,b)v$$ 
for all $a,b \in X, v \in L$.
\item[(C4)] There is a function $\phi: X \times L \to K$ such that 
\begin{align*} \theta(av,w)+\theta(a,w^{\sigma})^{\sigma} q(v) -f(w,v^{\sigma}) 
\theta(a,v)^{\sigma} \\
+f(\theta(a,v),w^{\sigma}) v^{\sigma}+\phi(a,v)w^{\sigma} \end{align*}
for all $a,b \in X, v,w \in L$.
\item[]
\item[(D1)] Let $\pi:X \to L: a \mapsto \theta(a.1)$. Then
$$a\theta(a,v)=(a\pi(a))v$$
for all $a \in X$ and all $v \in L$.
\item[(D2)] $\pi(a) \in \langle 1 \rangle$ if and only if $a=0$.
\end{itemize}
Then $(K,L,q,1,X,\cdot,h,\theta)$ is called a \textit{quadrangular algebra}.
\item 
A quadrangular algebra is called \textit{proper} if $\sigma \ne 1$, \textit{improper} if it is not isotopic 
to a proper quadrangular algebra (see 8.6 and 8.7 of \cite{W} for the definition of an isotope 
of a quadrangular algebra) 
\item A quadrangular algebra is called \textit{regular} if $f$ is non-degenerate 
and \textit{defective} if it is not regular. 
\item Let $\delta \in L$. A quadrangular algebra is called 
$\delta$-\textit{standard} if f$(\pi(a),\delta)=0$ for all $a\in X$ and 
\begin{itemize}
\item[\rm (i)] $\characteristic K \ne 2$ and $\delta=\frac{1}{2}(\in K \subset L)$. 
\item[\rm (ii)] $\characteristic K=2$ and $f(1,\delta)=0$.
\end{itemize}
\end{enumerate}
\end{definition}
Note that a regular quadrangular algebra is always 
proper and that defective quadrangular algebras only exist for fields with characteristic $2$. 
\begin{example}\rm\label{example}
\begin{enumerate}
\item Let $K$ be a field,
$L$ be a separable quadratic extension of $K$ or a quaternion division algebra over 
$K$, $q: L \to K$ the (reduced) norm, $\sigma$ the generator of $\mathrm{Gal}(L|K)$ or 
the standard involution of $L$, $(L,\sigma, X,\pi,h)$ an anisotropic pseudo-quadratic space. 
Let $\cdot: X\times L \to X$ be the scalar multiplication and set $\theta(a,u) =\pi(a)u$ 
for $a\in X$ in $u\in L$. Then $\Xi =(K,L,q,1,X,\cdot, h.\pi)$ is a quadrangular with 
$\phi=0$ in (C4) (see 1.18 of \cite{W}). A quadrangular algebra isomorphic to $\Xi$ is 
called \textit{special}. 
(For the definition of an isomorphism between quadrangular algebras see 
1.25 of \cite{W}). 
\item Let $(K,L,q)$ be a pointed quadratic space of type $\mathrm{E}_6,\mathrm{E}_7,\mathrm{E}_8$ or $F_4$ (see 2.13 
and 2.15 of \cite{W}) and $u\in L^{\#}$. Then there exist $X,\cdot,h,\pi$ such that 
$\Xi =(K,L, \frac{1}{q(u)}q, u, X,\cdot , h,\pi)$ is a quadrangular algebra (see 10.1 of 
\cite{W}). This quadrangular algebra is not special. If $(K,L,q)$ is of type $\mathrm{E}_6,\mathrm{E}_7$ or $\mathrm{E}_8$, 
then $\Xi$ is regular, if $(K,L,q)$ is of type $F_4$, then $K$ is a non-perfect field of 
characteristic $2$ and $\Xi$ is proper but defective.
\item Let $(L,X,\pi)$ be a quadratic space with bilinear form $h$ such that $\characteristic 
L=2$. Let $K$ be a subfield of $L$ with $L^2 \subseteq K \subseteq L$. Let $q:L \to K: u
\mapsto u^2$ and let $\theta(a,u) =\pi(a)u$ for $a\in X$ and $u\in L$. We regard $X$ and 
$L$ as vector spaces over $K$ and $K$ as a vector space over $L$ via the scalar multiplication 
$\ast: K \times L \to K: (a,u) \mapsto a \ast u = au^2$. Let $q_K$ be the embedding of $K$ into $L$. 
Then $q_K$ is a quadratic form since $q_K(a \ast u) = q_K(au^2) =au^2 =q_K(a)u^2$. Suppose that 
the quadratic space 
$$ (L,X,\pi) \perp (L,K,q_K) $$
is anisotropic. Then 
$$ \Xi=(K,L,q,1,X,\cdot, h,\theta) $$  
is an improper quadrangular algebra.         
\end{enumerate}
\end{example}
\begin{theorem}
Let $\Xi = (K,L,q,1,X,\cdot, h,\theta)$ be a quadrangular algebra with $h$ not identical zero.
\begin{enumerate}
\item If $\Xi$ is regular but not special, then $(K,L,q)$ is of type $\mathrm{E}_6,\mathrm{E}_7$ and $\mathrm{E}_8$ and 
$\Xi$ is uniquely determined by the pointed quadratic space $(K,L,q,1)$.
\item If $\Xi$ is proper but defective, then $(K,L,q)$ is a quadratic space of type $F_4$ and 
$\Xi$ is uniquely determined by $(K,L,q,1)$.
\item If $\Xi $ is improper, then $\Xi$ is isomorphic to a quadrangular algebra constructed in 
\ref{example}.  
\end{enumerate}
\end{theorem}
For a proof see 3.2,3.3 and 9.26 of \cite{W}. Note that also quadrangular algebras with 
$h$ identically zero can be classified (9.33 of \cite{W}), but they are not interesting for 
us.\bigskip
\\
We now turn our attention to Moufang quadrangles. If $\Gamma$ is a Moufang quadrangle, 
$\Sigma = (x_i)_{i \in \mathbb{Z}/8\mathbb{Z}}$ an ordered apartment, $\alpha_i=(
x_i,x_{i+1},
x_{i+2},x_{i+3},x_{i+4})$ 
and $U_i =U_{\alpha_i}$ the root group corresponding to $\alpha_i$ and 
$U_+ =(U_1, \ldots, U_4)$, then the \textit{root group sequence} $(U_+,U_1,U_2,U_3,U_4)$ 
determines the quadrangle up to isomorphism (see \cite[(7.5)]{TW}).  
It is well known (\cite[(5.5) and (6.1)]{TW}) that \begin{enumerate}
\item $[U_i,U_{i+j}] \leq U_{i+1} \ldots U_{i+j-1}$ for all $i,j$ with $0 < j < 4$.
\item $G_i:=\langle U_i, U_{i+4} \rangle$ is a rank $1$ group with unipotent subgroups $U_i$ 
and 
$U_{i+4}$.
\end{enumerate} 

Suppose that $\Xi =(K,L,q,1,X,\cdot, h,\theta)$ is a quadrangular algebra. Let $S=X \times K$, 
set 
$$ (a,t) \cdot (b,s) :=(a+b,s+t+g(b,a)),$$
where $a,b \in X, s,t, \in K$ and $g$ is the function in (C4). Then $(S,\cdot)$ is a group 
with $(a,s)^{-1} =(-a,-s+g(a,a))$ for $(a,s) \in S$. 
\begin{theorem}\label{quadrangular algebras} 
\begin{enumerate}
\item If $\Xi=(K,L,q,1,X,\cdot,h,\theta)$ is a proper quadrangular algebra 
which is $\delta$-standard for some $\delta\in L$ and $\phi$ as in 
(C4), then there is 
a Moufang quadrangle $\Omega_\Xi$ with root group sequence $(U_+,U_0,
\ldots, U_4)$ such that
\begin{itemize}
\item[\rm (i)] $U_1$ and $U_3$ are isomorphic to the group $S$ defined 
above via isomorphisms $x_1$ and $x_3$.
\item[\rm (ii)] $U_2$ and $U_4$ are isomorphic to $(L,+)$ va isomorphisms 
$x_2$ and $x_4$.
\item[\rm (iii)] $[U_1,U_2]=[U_2,U_3]=[U_3,U_4]=1$.
\item[\rm (iv)] For all $a,b \in X, s,t\in K$ and $u,v\in L$ we have
\begin{align*}
[(x_1(a,t),x_3(b,s)^{-1}] &=& x_2(h(a,b)),\hspace{4.4cm} \\
[x_2(u),x_4(v)^{-1}] &=& x_3(0,f(u,v))\text{ and} \hspace{3.4cm}\\
[x_1(a,t),x_4(v)^{-1}] &=& x_2(\theta(a,v)+tv)\cdot x_3(av,tq(v)+\phi(a,v)).
\end{align*}
\end{itemize}
\item If $\Xi$ and $\hat{\Xi}$ are quadrangular algebras, then 
$\Omega_\Xi $ and $\Omega_{\hat{\Xi}}$ are isomorphic iff $\Xi$ and 
$\hat{\Xi}$ are isotopic as defined in \cite[(8.7)]{W}.
\end{enumerate}
\end{theorem}
\begin{proof}
(a) follows from \cite[(11.12)]{W}, (b) is \cite[(11.13) and (11.14)]{W}
\end{proof}
\begin{definition}\rm Let $\Omega$ be a Moufang quadrangle with 
root group sequence $(U_+,U_1,\ldots,U_4)$. 
\begin{enumerate} \item
Set $V_i:=C_{U_i}(U_{i+2})$. Then $\Omega$ is called
\textit{indifferent} if $U_i=V_i$ for all $i$, \textit{reduced} if $U_i=V_i$ for some but not 
for all $i$ 
and \textit{wide} if $U_i \ne V_i$ for all $i$.
\item Suppose that $\Gamma$ is a wide Moufang quadrangle with root group sequence $(U_+, U_1, \ldots, U_4)$ 
such that $1 \ne V_i$ for all odd $i$. 
Then by \cite[(21.4)]{TW} there is Moufang subpolygon $\Omega$ of $\Gamma$ with root group sequence
$(V_+,V_1,U_2,V_3, U_4)$ with $V_+=V_1 U_2 V_3 U_4$. 
We call $\Gamma$ an \textit{extension} of $\Omega$.
\end{enumerate}
\end{definition}

We briefly recapitulate the classification of Moufang quadrangles. 
\begin{theorem}
Let $\Gamma$ be a Moufang quadrangle.
\begin{enumerate}
\item If $\Gamma$ is indifferent, then $$\Gamma \cong \mathcal{Q}_{\mathcal{D}}(K,K_0,L_0)$$ 
for some 
indifferent set $(K,K_0,L_0)$ with $\mathcal{Q}_{\mathcal{D}}(K,K_0,L_0)$ defined as in\cite[(16.4)]{TW}.
\item If $\Gamma$ is reduced, then either 
$$\Gamma \cong \mathcal{Q}_{\mathcal{I}}(K,K_0,\sigma)$$ for an involutory set
$(K,K_0,\sigma)$ such that $\sigma \ne 1$ and $K $ is generated by $K_0$ 
as a ring or 
$$\Gamma \cong \mathcal{Q}_{\mathcal{Q}}(K,L_0,q)$$
for an anisotropic quadratic space $(K,L_0,q)$. Here 
$\mathcal{Q}_{\mathcal{I}}(K,K_0,\sigma)$ 
and
$ \mathcal{Q}_{\mathcal{Q}}(K,L_0,q)$ are defined as in \cite[(16.2) 
resp. (16.3)]{TW}.
\item If $\Gamma$ is wide and the extension of the quadrangle 
$\mathcal{Q}_{\mathcal{I}}(K,K_0,
\sigma)$
for an involutory set $(K,K_0,\sigma)$ with $\sigma \ne 1$ and $K =\langle 
K_0 \rangle$, then 
$$\Gamma \cong \mathcal{Q}_{\mathcal{P}}(K,K_0,\sigma, L_0,q)$$ for an 
anisotropic pseudo-quadratic space $(K,K_0,\sigma, L_0,q)$, where the 
Moufang quadrangle 
$\mathcal{Q}_{\mathcal{P}}(K,K_0,\sigma, L_0,q)$ is defined as in 
\cite[(16.5)]{TW}.
\end{enumerate}
\end{theorem}

\section{Freudenthal triple systems and structurable algebras}
Freudenthal triple systems were introduced in order 
to have an algebraic approach to exceptional algebraic groups of type 
$\mathrm{E}_7$, see \cite{Br}, \cite{Fr} and \cite{Mey}. 
\begin{definition}\rm Let $\F$ be a $K$-vector space for a field $K$ with $\characteristic K \ne 
2,3$, $\langle.,.\rangle: \F \times \F \to K$ an alternating
$K$-linear form and $\{.,.,.\}:\F \times \F \times \F \to \F: (x,y,z) \mapsto \{x,y,z\}$ a 
symmetric, $K$-trilinear map. Then
$(\F,\langle.,.\rangle, \{.,.,.\})$ is called a \textit{Freudenthal triple system} if the following holds:
\begin{enumerate}
\item The map 
$q:\F \times \F \times \F \times \F \to K: (w,x,y,z) \mapsto \langle w, \{x,y,z\} \rangle$ is  
symmetric.
\item $\{x, \{x,x,x\}, y\} = \langle y,x\rangle \{x,x,x\} +\langle y, \{x,x,x\} \rangle x$ for all $x,y \in 
\F$.
\end{enumerate} 
A Freudenthal triple system $\F$ is called \textit{anisotropic} if $Q(x):= \langle x, \{x,x,x\} \rangle \ne 0$ 
for all $x\in \F \setminus \{0\}$.
\end{definition}
\begin{remark}\rm
Usually, the definition of a Freudenthal triple system also includes the condition that $\dim\F$ is 
finite, but we also want to regard infinite-dimensional Freudenthal triple systems. 
\end{remark}
\begin{example}
\rm 
\begin{enumerate}
\item
(see \cite{Fer}) 
Let $J$ be a quadratic Jordan algebra with base point 
$1$ and $\mathfrak{N}:J \to K$ be an admissible quadratic or 
cubic norm function. Let $a\mapsto a^\#$ 
be the adjoint map and $T: J\times J \to K$ the associated trace form. 
Then we have $\mathfrak{N}(1)=1$ and 
\begin{align*} 
\mathfrak{N}(a+b)=\mathfrak{N}(a)+T(a^\#,b)+T(a,b^\#) +\mathfrak{N}(b) 
\text{ and } (a^\#)^\#=\mathfrak{N}(a)a
\end{align*} for $\mathfrak{N}$ cubic resp.
\begin{align*}
a^\# =0 \text{ and } \mathfrak{N}(a+b)=\mathfrak{N}(a)+\mathfrak{N}
(b)+T(a,b)
\end{align*}
for $\mathfrak{N}$ quadratic, where $a,b\in J$.
Define $\mathfrak{M}(J,\mathfrak{N},1):=K\oplus K \oplus J \oplus 
J$, $$\langle (\alpha_1,\alpha_2,a_1,a_2). (\beta_1,\beta_2,b_1,b_2)\rangle 
:= \alpha_1 \beta_2-\alpha_2\beta_1+T(a_1,b_2)-T(a_2,b_1)$$ and 
\begin{align*}
Q(\alpha_1,\alpha_2,a_1,a_2):=\\
12(4\alpha_1\mathfrak{N}(a_1)+4\alpha_2
\mathfrak{N}(a_2)-4T(a_1^\#,a_2^\#)+(T(a_1,a_2)-\alpha_1\alpha_2))^2
\end{align*}
for $\alpha_1,\alpha_2,\beta_1,\beta_2\in K$ and $a_1,a_2,b_1,b_2\in J$.
Let $Q(.,.,.,.)$ be the linearisation of $Q$. Then for every $x,y,z \in 
\mathfrak{M}(J,\mathfrak{N},1)$ there is a unique element $\{x,y,z\}\in 
\mathfrak{M}
(J,\mathfrak{N},1)$ such that 
$$Q(x,y,z,w) =\langle w,\{x,y,z\}\rangle$$ 
holds for all $w\in \mathfrak{M}(J,\mathfrak{N},1)$, and 
$(\mathfrak{M}(J,\mathfrak{N},1),\langle.,.\rangle, \{.,.,.\})$ is a Freudenthal 
triple system. We say that a Freudenthal triple system $\F$ \textit{splits} 
if it isomorphic to $\mathfrak{M}(J,\mathfrak{N},1)$ for $(J,
\mathfrak{N},1)$ as above.
\item (see \cite[Theorem 4.1]{BDM}) Let $\Xi =(K,L,q,1,X,\cdot,h,\theta)$ 
be a quadrangular algebra with $\characteristic K \ne 2,3$. Let $g$ as 
in (C3) and $\pi$ as in (D1). Then $(X,g,\{.,.,.\})$ with 
$$\{x,y,z\} = \frac{1}{2}
(x(h(y,z)+h(z,y))+y(h(x,z)+h(z,x))+z(h(x,y)+h(y,x)))$$ is an anisotropic 
Freudenthal triple system.
\end{enumerate}
\end{example}

\begin{definition}\rm Let $\F$ be a Freudenthal triple system.
\begin{enumerate}
\item An \textit{ideal} of $\F$ is a subspace $\mathfrak{I}$ of $\F$ such 
that $\{x,y,z\}\in \mathfrak{I}$ holds for all $x,y\in \F$ and $z \in 
\mathfrak{I}$. 
\item
The Freudenthal triple system $\F$ is called \textit{simple} 
if $\F$ has no proper ideal.
\item An element $0\ne u\in \F$ is called \textit{strictly regular} if 
$\{u,u,x\}\in \langle u \rangle$ for all $x\in \F$.
\item $\F$ is called \textit{reduced} if there is a pair of strictly 
regular elements $\{u,v\}$ with $\langle u,v\rangle =1$.
\end{enumerate}
\end{definition}

It can be shown that $\F$ is simple iff $\langle.,.\rangle$ is non-degenerate 
(\cite[Theorem 2.1]{Fer}). Hence an anisotropic Freudenthal triple system 
is always simple. 
Note that if $u\in \F$ is strictly regular, then $\{u,u,u\} =\lambda u$ 
for some $\lambda\in K$, thus $Q(u)=\langle u,\{u,u,u\}\rangle =\lambda \langle u,u\rangle 
=0$. Thus an anisotropic Freudenthal triple system contains no strictly 
regular element.

\begin{theorem}\label{Ferrar}
Let $\mathfrak{F}$ be a Freudenthal triple system. Then
\begin{enumerate}
\item $\F$ is reduced iff it contains an element $u\ne 0$ such 
that $12 Q(u)$ is a square.
\item If $\mathfrak{F}$ is a simple, reduced Freudenthal triple system 
with $\dim_K \F <\infty$, then $\F\cong \mathfrak{M}(J,\mathfrak{N},1)$ 
for a quadratic Jordan division $J$ and an admissible quadratic or cubic 
norm function $\mathfrak{N}:J\to K$.
\end{enumerate}
  \end{theorem}
\begin{proof}
See \cite[3.12 and 5.1]{Fer}.
\end{proof}

\bigskip 
It follows that a simple, finite-dimensional 
Freudenthal triple system over a quadratically closed 
field always splits.
Moreover, if $K$ is arbitrary, then 
there is a quadratic extension field $L$ of $K$ such that 
$L\otimes_K \F$ splits.

Closely related to the theory of Freudenthal triple systems is the concept of a structurable 
algebra.
\begin{definition}\rm
Let $K$ be a field of characteristic different from $2$ or $3$ and $\mathcal{A}$ a (not 
necessarily associative) $K$-algebra with 
involution $\bar{}$. 
For $x,y \in \mathcal{A}$ define $V_{x,y} \in \End(\mathcal{A})$ by $V_{x,y} z :=(x\bar{y})z 
+(z\bar{y})x -(z\bar{x})y$. Then 
$(\mathcal{A},\bar{.} )$ is called a \textit{structurable algebra} iff $$[V_{x,y}, V_{z,w}] 
=V_{V_{x,y},w}-V_{z,V_{y,x}w}$$
holds for all $w,x,y,z \in \mathcal{A}$. 
\end{definition}
\begin{remark}\rm If $\mathcal{A}$ is a $K$-algebra for a field $K$ with $\characteristic K \ne 
2,3$,
then $(\mathcal{A},\mathrm{id})$ is a structurable algebra iff $\mathcal{A}$ is a (linear) 
Jordan 
algebra.
\end{remark} 
\begin{definition}\rm
Let $(\mathcal{A},\bar{.})$ be a structurable algebra. Set $\mathcal{S}=\{a \in \mathcal{A}| 
\bar{a}=-a\}$
Then $\skewdim \mathcal{A}:=\dim \mathcal{S}$ is called
the \textit{skew-dimension} of $\mathcal{A}$. 
\end{definition}
\begin{prop} Let $(\mathcal{A},  \bar{.}  )$ be a structurable algebra of skew-dimension $1$ 
and $s_0 \in 
\mathcal{S} \setminus 
\{0\}$. Let $\langle.,.\rangle :\mathcal{A} \times \mathcal{A} \to K$ be defined by $x\bar{y}-
y\bar{x}=\langle x,y\rangle 
s_0$ and
define $\{.,.,.\}: \mathcal{A} \times \mathcal{A} \times \mathcal{A} \to \mathcal{A}$ by 
$$\{x,y,z\}:=
2 V_{x,s_0 y} z-\langle y,z\rangle x -\langle y,x \rangle z-\langle x,z \rangle y .$$
Then $(\mathcal{A}, \langle.,.\rangle, \{.,.,.\})$ is a Freudenthal triple system.
\end{prop}
\begin{proof}
Set \cite[Proposition 2.18]{AF}.
\end{proof}
\begin{definition}\rm Let $(\mathcal{A},\bar{.})$ be a structurable algebra. An element $x \in 
\mathcal{A}$ is called \textit{conjugate invertible} iff there is an element $\hat{x}\in \mathcal{A}$ 
such 
that $V_{\hat{x},x} =V_{x,\hat{x}} =\mathrm{id}$. The algebra $(\mathcal{A},\bar{.})$ is called 
a 
\textit{structurable division algebra} if every $x \in \mathcal{A}^{\#}$ is invertible. 
\end{definition}
\begin{remark}\rm \begin{enumerate}
\item 
If $x\in \mathcal{A}$ is conjugate invertible, then $\hat{x}$ is uniquely determined, 
and we have $\hat{\hat{x}}=x$ (see \cite[Lemma 6.1]{AH}.
\item If $\F$ is a Freudenthal triple system and $(\mathcal{A},\bar{•})$ its corresponding structurable structurable algebra, 
then $x\in \F$ is conjugate invertible iff $Q(x) \ne 0$. Thus $\F$ is anisotropic iff 
$\mathcal{A}$ 
is a structurable division algebra.
\end{enumerate}
\end{remark}
In \cite[Theorem 5.1.6]{BDMS} the authors show that every structurable division algebra 
$\mathcal{A}$ gives rise to a Moufang set $\mathbb{M}(\mathcal{A})$, see also \cite{Boe}. The authors give
an explicit formula for the root group and the map $\tau$. 
Furthermore, they showed that if $\characteristic K \ne 2,3$, then 
for every semisimple algebraic group defined over $K$ of relative 
$K$-rank one the Moufang set $\mathbb{M}(G)$ can be obtained in this 
fashion, see \cite[Theorem 5.2.1]{BDMS}.

Translating their results into the language of Freudenthal triple 
systems, we get 
\begin{theorem}\label{Freudenthal Moufang}
Let $K$ be a field with $\characteristic K \ne 2,3$.
\begin{enumerate} 
\item Let $\F$ be an anisotropic Freudenthal triple system over $K$. Set 
$U =\F \times K$ 
with 
addition $(x,s) +(y,t) =(x+y,s+t+\langle x,y\rangle)$. Define $\tau: 
U^{\#} \to U^{\#}$ by 
$$(x,t)\tau = \left(-\frac{1}{Q(x)-12t^2}(2\{x,x,x\}+12tx),\frac{12t}
{Q(x)-12t^2}\right).$$
Then $M(U,\tau)$ is a Moufang set which we will denote by $\mathbb{M}
(\F)$.
\item Let $G^\dagger$ be the little projective group for $\mathbb{M}
(\F)$ with $\F$ an anisotropic Freudenthal triple system over $K$. If 
$\dim_K \F <\infty$, then $G^\dagger$ is an algebraic group defined over 
$K$ of type ${}^{3}\mathrm{D}^9_{4,1},{}^{6}\mathrm{D}^9_{4,1} 
{}^{2}\mathrm{E}^{35}_6, 
\mathrm{E}^{66}_{7,1}$ or $\mathrm{E}^{133}_{8,1}$.
\item Conversely, suppose that $G$ is a semisimple algebraic group 
defined over $K$ of type ${}^{3}\mathrm{D}^9_{4,1},
{}^{6}\mathrm{D}^9_{4,1} {}^{2}\mathrm{E}^{35}_6, 
\mathrm{E}^{66}_{7,1}$ or $\mathrm{E}^{133}_{8,1}$. Then there is an 
anisotropic Freudenthal triple system $\F$ over $K$ with 
$\mathbb{M}(G)\cong \mathbb{M}(\F)$.
\end{enumerate}
\end{theorem} 
\begin{proof}
Part (a) follows from formula (6.9) of \cite{Boe} resp. Theorem 5.1.6 of \cite{BDMS}. Since $\mathcal{S}=K s_0$, we may replace 
$\alpha \mu s_0 \in k s_0$ by $t\in k$ and then apply 3.27 and 3.28 of 
\cite{Boe}.
Part (b) and (c) are \cite[Theorem 6.4.3 and 6.4.4]{BDMS}. 
\end{proof}


  
\chapter{Cubic Action} 
From now  on, we assume that $G =\langle A,B \rangle$ is a rank one group. 
We repeat the definition for cubic action:
\begin{definition}\rm
Let $V$ be a $\mathbb{Z}G$-module. Then $G$ acts \textit{cubically} on 
$V$ 
if $[V,A,A,A]$ $ =0$ 
but $[V,G,G,G] \ne 0\ne [V,A,A]$. In this case $V$ is called a 
\textit{cubic module}
 for $G$, and $G$ is called 
a \textit{cubic rank one group}.
\end{definition} 
For the rest of this chapter $V$ is always a cubic module for $G$ on which $G$ 
acts faithfully. 
\begin{lemma}\label{faithful2} 
$A$ acts faithfully on $V$.
\end{lemma}
\begin{proof}
Since we have $[V,A,A,A]=0$ but $[V,G,G,G] \ne 0$, the image of 
$A$ in $\GL(V)$ is properly contained in the image of $G$. 
Thus $A$ acts faithfully on $V$ by \ref{faithful}(d).\end{proof}

\bigskip
We now define an important subgroup of $A$. We set
$$A_0:=C_A(V/C_V(A)) \cap C_A([V,A]).$$
Then we have 
\begin{prop}\label{quadratic kernel} \begin{enumerate}
\item $A^ {\prime} \leq A_0 \leq Z(A)$.
\item If $A_0 \ne 1$, then $A_0$ is a special root subgroup of $A$ which acts quadratically on 
$V$.
\item $A_0$ is $H$-invariant.
\end{enumerate}
\end{prop}
\begin{proof} \begin{enumerate}
\item 
This follows easily with the $3$-subgroup lemma: We have $[V,A,A_0] =0$ and $[V,A_0,A] \leq 
[C_V(A),A] =0$, 
thus $[V,[A,A_0]]=0$. Since $A$ acts faithfully on $V$, this implies $A_0 \leq Z(A)$. Moreover,
$[[V,A],A,A] =[V,A,A,A]=0$ implies $ [[V,A],A^{\prime}] =0$, and $[V/C_V(A),A,A] =0$ implies 
$[V,C_V(A),A^{\prime}] =0$. Thus $A^{\prime} \leq A_0$.
\item Let $a \in A_0^*$, $\mu =\mu(a)$,  $B_0 = A_0^{\mu} =C_B([V,B]) \cap C_B(V/C_V(B))$ 
and $G_0=\langle A_0,B_0 \rangle$. 
We have $[V,A_0] \subseteq C_V(A)$ and thus 
 $[V,A_0,A] \subseteq [C_V(A),A] = 0$. Suppose that $[V,G,G_0] = 0$. Then $G_0 \leq N:=
 C_G([V,G])$. We have $1 \ne A_0 \leq G_0 \leq N$, so $A\cap N\ne 1$ and thus $G=NA$ by
 \ref{faithful}. So $[V,G,G,G,G] =[V,G,A,A,A] = 0$. This implies that $G$ is nilpotent, a 
 contradiction. 
Thus $A_0$ is a special root subgroup of $A$ by \cite[I (2.4)]{T1}. Since $[V,A_0,A] = 0$, we 
get $[V,A_0,A_0] =0$. Thus $A_0$ acts quadratically on $V$.
\item Since $H$ normalises $A$, $H$ normalises both $C_V(A)$ and $[V,A]$ and 
thus $H$ 
normalises $C_A([V,A])$ and 
$C_A([V,A])$ and also the intersection of these two groups which is $A_0$.  
\end{enumerate}
\end{proof}
\begin{definition}\rm
The group $G_0:=\langle A_0,B_0\rangle $ with $B_0 = \{1\} \cup \{\fatb(a)|a\in 
A_0^*\}$ is 
called the \textit{quadratic kernel} of $G$.
\end{definition}
As a consequence of \ref{quadratic kernel} and \cite[Theorem 1.1]{T2} we 
get
\begin{coro}
Suppose that $G_0 \ne 1$. Then we have
\begin{enumerate}
\item $G_0$ is a rank one group with unipotent subgroups $A_0$ and $B_0$.
\item There is a ring $R$ and a special quadratic Jordan division algebra $J 
\subseteq R$ such 
that $G_0/C_{G_0}(V) \cong \SL_2(J,R)$.
\end{enumerate}
\end{coro}
\chapter{Examples for cubic modules}
\section{Pseudoquadratic spaces}\label{pqf}
Let $(K,K_0,\ast)$ be an involutory set, $\oV$ a $K$-vector space and ${\opi:
\oV \to K/K_0}$
 an anisotropic pseudoquadratic form with associated skew-hermitian form 
 $f$. Define 
 $V = K \oplus \oV \oplus K$ and $\pi:V \to K/K_0$ by 
$$\pi(r,x,s) = s^\ast r + \opi(x) +K_0.$$ Then $\pi$ is a pseudo-quadratic form 
with associated 
skew-hermitian form $g$ 
given by
$$g((r,x,s),(t,y,u)) = r^\ast u - s^\ast t +f(x,y)$$
for $x,y \in \oV$ and $r,s,t,u \in K$. 
We have $\pi(r,x,0) \in K_0$ iff $x=0$ and $\pi(r,x,1) =0$ iff $r +\opi(x) 
\in K_0$.
Thus $\{(1,0,0)K\} \cup \{(r,x,1)K| r+\opi(x) \in K_0\}$ is the set of 
anisotropic $1$-dimensional spaces. Therefore the Witt index of $\pi$ is 
$1$. 
For $v \in \oV, t \in K$ with $\pi(v) -t \in K_0$ let 
$\alpha_{(v,t)} \in \End(V)$ be defined by $$(r,x,s)\alpha_{(v,t)} = 
(r-f(v,x) + t^\ast s, x +vs, s)$$ for $x \in \oV, r,s \in K$. Then 
$\alpha_{(v,t)} \in \GL_K(V) $, 
$\alpha_{(v,t)} \alpha_{(w,u)} = \alpha_{(v+w, t+u + f(v,w))}$ and
\begin{eqnarray*} 
\pi( (r,x,s)\alpha_{(v,t)}) = \pi ( (r-f(v,x) + t^\ast s,x+vs,s))  = \hfill\\
 s^\ast (r-f(v,x) +t^\ast s) +\opi(v+xs) + K_0  = \\
 s^\ast r +s^\ast f(x,v)^\ast + s^\ast t^\ast s +\opi(v) +f(x,vs) + s^\ast \opi(x) s + K_0  = 
 \\
s^\ast r+ \opi(v) +\underbrace{s^\ast (t^\ast + t -(t -\opi(v)) s + f(x,v) s + (f(x,v)s)^\ast}_{\in K_0}+ K_0 
 = 
\pi(r,x,s).
\end{eqnarray*}
Similarly, define $\beta_{(t,v)}$ by $$(r,x,s) \beta_{(t,v)} = (r,x-vr,s-f(v,x)-t^\ast r).$$
Then again $\pi((r,x,v)\beta_{(t,v)}) = \pi(r,x,v)$. Set $A=\{\alpha_{(v,t)}| 
\opi(v) -t \in 
K_0\}$, $B =\{\beta_{(v,t)}| \opi(v) -t \in K_0\}$ and $G:=\SU(V,\pi):=\langle A,B \rangle$. 
Then $G$ is an abstract rank one group and $V$ is a cubic module for $G$. 
(This can be seen either by a direct computation or by the 
fact that this example lives inside the Moufang quadrangle 
$\mathcal{Q}_{\mathcal{P}}(K,K_0,\ast,\oV,\opi)$, see 
section \ref{Moufang quadrangle}.)
It is 
$[V,A] = \{(r,v,0)| r \in K, v \in V\}$ and $C_V(A) = \{(r,v,0)| r \in K, v 
\in \rad(f)\}$, 
where $\rad(f):=\{v\in \oV| f(v,w)=0 \ \forall w\in V\}$.
Therefore $$A_0 = C_A([V,A]) =\{\alpha_{(v,t)}| v \in \rad(f)\}=Z(A).$$
Now let $V'=\{(r,x,s)|r,s\in K,x\in \rad(f)\}$ and $V''=\{(0,x,0)|x\in 
\rad(f)\}$. Then $G_0$ acts on $V'/V''\cong K^2$, and it can 
be easily seen that $G_0/C_{G_0}(V'/V'')\cong \SL_2(J,K)$ for $J=K_0+\{\pi(x)| x\in \rad(f)\}$. 
\section{Adjoint action}
Suppose that $G =\langle A,B \rangle$ is a rank one group and that $V$ is a quadratic module for $G$. 
If $F$ is a subring of $\End_G(V)$, 
then $G$ acts on $\End_F(V)$ by conjugation. 
We claim that this action is cubic. For $\varphi \in \End_F(V), a \in A$ and $v \in [V,A]$ we 
have 
$$v(-\varphi + a^{-1} \varphi a ) = -v\varphi +v\varphi a =[v\varphi,a] \in [V,A], $$
so $[\End_F(V),A] \subseteq \{\varphi \in \End_F(V)| [V,A]\varphi \subseteq [V,A]\}$.
For $\varphi \in [\End_F(V),A], v \in V$ we get
$$v(-\varphi + a^{-1} \varphi a) = -v\varphi +v\varphi a + [v,a^{-1}]\varphi a =
[v\varphi,a] + [v,a^{-1}]\varphi \in [V,A].$$
If $v \in [V,A]$, then $v(-\varphi + a^{-1} \varphi a) =0$. So we get 
$[\End_F(V),A,A] \subseteq \{\varphi \in \End_F(V)| V\varphi \subseteq [V,A], [V,A]\varphi 
=0\}$. 
With the formula above we obtain that ${[\End_F(V),A,A,A] =0}$. 

By \ref{Jordan} $G\cong \SL_2(J,R)$ for a quadratic Jordan division algebra inside a ring $R$, so 
either $G$ is perfect or $G\cong \SL_2(p)$ and $V\cong \mathbb{F}_p^2$ for $p\in \{2,3\}$.
So one sees that $[\End_F(V),G,G,G] \ne 0$. 
If $G=\SL_2(J,R)$ with $R=\langle J \rangle$, $V=R^2$ and $A$ the group of lower unipotent triangle 
matrices with entries in $J$, 
then $C_A([\End_F(V),A]) =A$ iff $R$ is commutative and $\characteristic R =2$ and 
$C_A([\End_F(V),A]) =1$ in all other cases. 
In all cases $A_0 =1$.
\section{The Tits-Kantor-Koecher module}\label{TKK}
Let $J$ be a quadratic Jordan division algebra over a field $K$. For $x,y \in J$ let 
$D(x,y)=(V_{x,y},-V_{y,x}) \in \End_K(J)^2$ and let $\mathcal{D}$ be the subspace of 
$\End_K(J)^2$ 
generated by $\{D(x,y)|x,y \in J\}$. Note that if $(\delta_1,\delta_2) \in \mathcal{D}$, then 
also $(\delta_2,\delta_1)\in \mathcal{D}$, so $\bar{.}: \mathcal{D}\to \mathcal{D}:
(\delta_1,\delta_2) \mapsto (\delta_2,\delta_1)$ is a well-defined involution of 
$\mathcal{D}$. 
Set $\mathcal{K}=J \oplus \mathcal{D} \oplus J$.
Let $[.,.]: \mathcal{K} \times \mathcal{K} \to \mathcal{K}$ be bilinear and anti-symmetric, 
such that for $x,y\in J $ and $(\delta_1,\delta_2)\in \mathcal{D}$ we have
\begin{align*}
[(x,0,0,0),(0,0,0,y)] = (0,V_{x,y},-V_{y,x},0), \\
[(0,\delta_1,\delta_2,0),(x,0,0,0)] =(x\delta_1,0,0,0), \\
[(0,\delta_1,\delta_2,0),(0,0,0,y)] =(0,0,0,y\delta_2) \end{align*}
and $[.,.]$ is zero in all other cases. Then $(\mathcal{K},[.,.])$ is a $\{-1,0,1\}$-graded 
Lie algebra (see \cite[Section VIII.5]{J}).  
For $a\in J$ let $\alpha_a \in \End_K(\mathcal{K})$ defined by 
$$(x,\delta_1,\delta_2,y)\alpha_a =(x+a\delta_1+yQ_a,\delta_1 +V_{a,y}, \delta_2-V_{y,a},
y)$$
for $x,y \in J$ and $(\delta_1,\delta_2) \in \mathcal{D}$. Then a straightforward but tedious computation 
shows that $\alpha_a$ is a Lie automorphism of $\mathcal{K}$.
If $\characteristic K \ne 2$, 
then we have $\alpha_a =\exp (\mathrm{ad} (a,0,0,0)) =\mathrm{id} +\mathrm{ad} (a,0,0,0) +
\frac{1}{2} (\mathrm{ad} (a,0,0,0))^2$, where $\mathrm{ad}: \mathcal{K} \to \End(\mathcal{K}):
x  \mapsto [x,.]$. 
One easily computes $\alpha_a \alpha_b=\alpha_{a+b}$, hence $A:=\{\alpha_a| a\in J\}$ is a 
subgroup of $\GL_K(\mathcal{K})$ isomorphic to $(J,+)$. Let $\tau \in \End_K(\mathcal{K})$ 
be defined as 
$$(x,\delta,y)\tau =(y,\bar{\delta},x)$$
for $x,y\in J$ and $\delta \in \mathcal{D}$. Again, $\tau$ is a Lie automorphism of $\mathcal{K}$.
For $b\in J$ we set $\beta_b=\alpha_b^\tau$ and 
let $B:=A^{\tau}$. 
Then we have
\begin{align*}
(x,\delta_1,\delta_2,y)\beta_b=(y,\delta_2,\delta_1,x)\alpha_b\tau = 
(y+b\delta_2+xQ_b,\delta_2 +V_{b,x},\delta_1-V_{x,b},x)\tau =\\
(x,\delta_1-V_{x,b},\delta_2+V_{b,x},y+b\delta_2+xQ_b). 
\end{align*}
Set $L_{\infty} :=\{(x,0,0,0)|x \in J\} \leq \mathcal{K}, 
M_\infty:=\{a\in \mathcal{K}| [L_\infty,a] \subseteq L_\infty\} =\{(x,\Delta,0)|x\in J, 
\Delta \in \mathcal{D}\},
L_0:=\{(0,0,0,y)|y \in J\} \leq 
\mathcal{K}$ and $L_a =L_0 \alpha_a$ for $a\in J$. Then for $a\in J^{\#}$ we have by \ref{identities}
\begin{align*}
L_a =\{(xQ_a,V_{a,x},-V_{x,a}, x)| x\in J\} =\{(y,V_{a,yQ_a^{-1}},-V_{yQ_a^{-1},a},yQ_a^{-1})| 
y\in J\} = \\
\{(y,V_{y,a^{-1}},-V_{a^{-1},y},yQ_{-a^{-1}})|y\in J\} =\\
\{(y,-V_{y,-a^{-1}},V_{-a^{-1},y}, 
y Q_{-a^{-1}})| y\in J\} = L_{\infty}\beta_{-a^{-1}} 
\end{align*} 
Moreover, for $a\in J^\#$ we have again by \ref{identities}
\begin{align*}
L_a\tau = \{(xQ_a,V_{a,x},-V_{x,a},x)\tau|x\in J\} = \{(x,-V_{x,a},V_{a,x},xQ_a)|a\in J\} = \\
\{xQ_a^{-1},-V_{xQ_a^{-1},a},V_{a,xQ_{a^{-1}}},x)|x\in J\}=
\{(xQ_{-a^{-1}},V_{-a^{-1},x},-V_{x,-a^{-1}},x)|x\in J\}\\
 =L_{-a^{-1}}.
\end{align*}
Set $G=\langle A,\tau \rangle$, 
$\Omega:=\{L_{\infty}\} \cup \{L_a|a\in J\} =\{L_0\} \cup \{L_{\infty}\beta_{a}|a\in J\}$, 
$Z =C_G(\Omega)$, $\bar{G}:=G/Z$ and $\bar{A}:=AZ/Z$. Then $(\Omega, (\bar{A}^g)_{g\in \bar{G}})$ is a 
Moufang set isomorphic to $\mathbb{M}(J)$. Moreover, $[Z,A]$ fixes all elements in $\Omega$ and 
centralises $L_\infty$, $M_\infty/L_\infty$ and $\mathcal{K}/M_\infty$. Since $L_0$ is a $[Z,A]$-invariant 
complement for $M_\infty$ in $\mathcal{K}$, we get that $[Z,A]$ also centralises $L_0$. Therefore 
$Z_A$ centralises $\mathcal{K}=L_0 +L_\infty + [L_0,L_\infty]$. Therefore $[Z,A]=1$, and by the same 
argument we get $[Z,B]=1$. Thus $Z =Z(G)$ and $G$ is a rank one subgroup with unipotent subgroups $A$ and 
$B$ acting cubically on $\mathcal{K}$.

\section{Quadratic pairs without commuting root groups}
\label{algebraic groups}

As Timmesfeld noted in \cite{T3}, every quadratic pair without commuting root subgroups gives 
rise 
to a cubic module. Recall from the introduction that a quadratic pair consists of 
finite-dimensional 
$K$-vector space $M$ with $K$ a field of characteristic not $2$, and a group $G$ generated 
by a set $\mathcal{Q}\leq \GL_K(M)^*$ consisting of quadratic elements, hence $[M,\sigma]\leq 
C_M(\sigma)$ 
for all $\sigma\in \mathcal{Q}$. We demand moreover 
that for $\lambda\in K$ and $\sigma\in \mathcal{Q}$ also the 
element $\lambda\circ \sigma:=\id +\lambda(\sigma-\id)$ is contained in $\mathcal{Q}$.  

We set $d(\sigma)=\dim [M,\sigma]$ for $\sigma\in \mathcal{Q}$, 
$d=\min\{d(\sigma)|\sigma\in \mathcal{Q}\}$ and $\mathcal{Q}_d=\{\sigma\in \mathcal{Q}|
d(\sigma)=d\}$. 
For $\sigma\in \mathcal{Q}_d$ we set $E_\sigma=\{\tau\in \mathcal{Q}| [M,\sigma]=[M,\tau]\}$ 
and 
$E(\sigma) =\langle E_\sigma\rangle$. It was shown in the introduction of \cite{T0} that $E(\sigma) 
=E_\sigma\cup 
\{\id\}$ holds. We set $\Sigma=\{E(\sigma)|\sigma\in \mathcal{Q}_d\}$. The elements of 
$\Sigma$ are called \textit{root groups}. 

\begin{theorem}
Let $A,B\in E(\Sigma)$ be two distinct root groups. Then one of the following cases holds
\begin{itemize}
\item[\rm (i)] $A$ and $B$ commute.
\item[\rm (ii)] $\langle A,B\rangle$ is nilpotent of class $2$, and $[A,B]=Z(\langle A,B\rangle) \in 
\Sigma$.
\item[\rm (iii)] $\langle A,B\rangle$ is a rank one group with unipotent subgroups $A$ and $B$.
\end{itemize}
\end{theorem}
\begin{proof}
This follows from \cite[Theorem 2]{T0}, see also \cite[Chapter V, Theorem (1.5)]{T1}. 
\end{proof}

A group $G$ generated by a family of subgroups $\Sigma$ that satisfy the conditions above are 
called \textit{groups generated by abstract root groups} and were examined in \cite{T0} and 
\cite{T1}. To determine the possible structure of $G$ one may, by applying \cite[Theorem 1]
{T0}, assume that $G$ is quasi-simple. If $\Sigma$ contains two distinct root groups that 
commute, or equivalently, if the geometry associated to $\Sigma$ has rank at least $2$, then 
the possible structure of $G$ was determined by Timmesfeld, compare \cite{T0} and \cite{T1}. 
This leaves the case that two distinct root groups always generate a rank one group. 
In this case we have:
  
\begin{theorem}\label{quadratic pairs}
Suppose that the following conditions are satisfied:
\begin{enumerate}
\item $M =[M,G]$.
\item $G=\langle \Sigma\rangle$, but $G\ne \langle A,B\rangle$ for all $A,B\in \Sigma$.
\item For $A,B\in \Sigma$ with $A\ne B$ the group $X=\langle A,B\rangle$ is an abstract rank one 
group with unipotent subgroups $A$ and $B$.
\end{enumerate}
For $A \in \Sigma$ we set $$U(A):= \{\varphi \in G| [M,\varphi] \subseteq C_M(A), [C_M(A),
\varphi] \subseteq [M,A] \hbox{ and } 
[M,A,\varphi] =0\}.$$
Then for all $A,B \in \Sigma$ distinct the group $G$ is an abstract rank one group with unipotent subgroups 
$U(A)$ and $U(B)$, and $M$ is a cubic module for $G$. 
\end{theorem}
\begin{proof}
This follows by \cite[Proposition 1]{T3}.
\end{proof}

\bigskip
Note that this statement is not true for $\characteristic K = 2$. For example, let $K$ be a commutative 
field of characteristic $2$ with $|K|>2$, let $L/K$ be a quadratic Galois extension and consider 
$G=\PSU_3(L)$ acting on the vector space $M$ of $3\times 3$-matrices of trace $0$ with entries in $L$.
Furthermore, set $\Sigma:=\{Z(A)|A\leq G \text{ unipotent}\}$. Then $G$, $\Sigma$ and $M$ satisfy 
the hypothesis of the hypothesis, but we have $U(A)=A$, and though $G$ is a rank one group, the action of 
$G$ on $M$ is 
not cubic.

For the case that $G$ can be generated by two but not by three elements of $\Sigma$ Timmesfeld proved 
the following result: 
\begin{theorem}
Suppose that in the situation of \ref{quadratic pairs} the group $G$ is generated by three distinct 
elements of $\Sigma$. Then one of the following cases holds:
\begin{enumerate}
\item There are extension fields $E \supset F \supset K$ such that $E/F$ is a quadratic extension 
with non-trivial Galois automorphism $\sigma$ and a 
$\sigma$-skew-hermitian form $f:E^3\times E^3\to E$ of Witt index $1$ 
such that $G\cong \SU(f)$. In this case we have $\langle A ,B\rangle\cong \SL_2(F)$ for all $A,B\in \Sigma$ with 
$A\ne B$. 
\item There is a division algebra $L$ containing $K$, an involution $*$ of $L$ such that $L$ is 
generated by $L_*$ and a $*$-hermitian form 
$f:L^3\times L^3\to L$ of Witt index $1$ such that $G\cong \SU(f)$. In this case we have 
$\langle A,B\rangle \cong \SL_2(K_*,K)$ for all $A,B\in \Sigma$ with $A\ne B$.
\end{enumerate}
\end{theorem}
Now let $G=\mathbf{G}(K)$ be a universal algebraic group of type ${}^3 \mathrm{D}_{4,1}, {}^6 \mathrm{D}_{4,1}, {}^2 \mathrm{E}_{6,1}^{35}$ 
or $\mathrm{E}_{7,1}^{66}$ defined over a field $K$ with $\characteristic K \ne 2$ and let 
$\bar{K}$ be the algebraic closure of $K$, Then $\mathbf{G}({\bar{K}})$ 
has an irreducible module $M$ such that $G_{\bar{K}}$ 
is generated by quadratic elements (see for example \cite{PS}). Examining 
the Tits diagram 
we see that also $G \leq \mathbf{G}({\bar{K}})$ is generated by quadratic elements. Thus 
$M$ satisfies Hypothesis (H) and is therefore a cubic module for $G$.  
\section{Quadratic spaces}\label{orth}
We present an example with abelian root groups such that $A_0 \ne 1$ is possible.
Let $K$ be a commutative field and let $(L_0,q)$ be an anisotropic quadratic space over $K$ 
with associated 
symmetric bilinear form $f$. Set $\oL_0 = L_0/\Def(q)$.
Here, $\Def(q) =L_0^{\perp} =\{ v \in L_0| f(v,w) =0$ for all $w \in L_0\}$. 
For $v \in L_0$, $\ov$ denotes the image of $v$ in $\oL_0$. 
Set $V = K \oplus \oL_0 \oplus K$. For all $v \in L_0$ we define $\alpha_v, \beta_v \in 
\GL_K(V)$ by 
$$(x,\ow,y) \alpha_v = (x,\ow +\ov x, y+f(\ow,\ov) + xq(v))$$ and 
$$(x,\ow,y)\beta_v =(x+f(\ow,\ov) +q(v)y,\ow +\ov y,y).$$
Then $\alpha_v \alpha_w = \alpha_{v+w} $ and $\beta_v \beta_w = \beta_{v+w}$ for all $v,w \in 
L_0$. Thus 
$A:=\{\alpha_v| v \in L_0\}$ and $B:=\{\beta_v| v\in L_0\}$ are groups isomorphic to $L_0$. If 
$\Def(q) =0$ 
(which is always case if $\characteristic K \ne 2$), then $Q:V \to K: Q(x,\ow,y) = xy - q(w)$ 
is a well-defined quadratic form of
 Witt index $1$. 
 
 Suppose $\Def(q) \ne 0$. Then we may assume that there is an element $e \in 
 \Def(q)$ with $q(e) =1$ 
(if not, we replace $q$ by $q(e)^{-1} q$ for some $e\in L_0^\#$). Set $K_0 := q(\Def(q))$. Then 
$(K,K_0,\mathrm{id})$ is an 
involutory set. 
If we define $Q:V \to K/K_0$ by $Q(x,\ow,y) = xy + q(w) +K_0$ (here $w$ is a preimage of $\ow$ 
in $L_0$), then $Q$ is a well-defined 
pseudo-quadratic form. 

In both cases, $Q$ has Witt index $1$ and the two groups $A$ and $B$ are 
contained in $O(Q):=\{
\varphi \in \GL(V,K)| Q(z\varphi) =Q(z)$ for all $z \in V\}$. 
If $X =\{(0,0,1)K\} \cup \{(1,\ow,q(w))K| w \in L_0\}$, then $X$ is the set of isotropic $1$-
dimensional subspaces of $V$.
We write $\infty$ for $(1,0,0)K$ and $v$ for $(q(v),\ov,1)K$. One can see that $A$ is the 
centraliser in $O(Q)$ of 
$\infty$, $V/\infty^{\perp}$ and $\infty^{\perp}/\infty$. Similarly, $B$ is the centraliser of 
$0$, $V/0^{\perp}$ and 
$0^{\perp}/0$.
 
Moreover, since $\alpha_v$ maps $w$ to $w+v$, $A$ acts regularly on $X \setminus \{\infty\}$.
Similarly, $B$ acts regularly on $X \setminus \{0\}$. 
Thus $G =\langle A,B \rangle$ is a rank one group with unipotent subgroups $A$ and $B$. Since $[V,A] 
\leq
\infty^{\perp}, [\infty^{\perp},A] \leq \infty$ and $[\infty,A] =0$, $V$ is a quadratic or 
cubic module for $G$. It is 
quadratic iff $V= \Def(q)$. One can easily see that $A_0 =\{\alpha_v| v \in \Def(q)\}$, so $A_0 
\ne 1$ is only possible if $\characteristic K =2$.
\section{Connection to Moufang Quadrangles}\label{Moufang quadrangle}
Let $\Omega$ be a Moufang quadrangle and $\Sigma=(x_i)_{i\in \ZZ/8\ZZ}$ an ordered apartment in $\Omega$. 
For $i\in\ZZ/8\ZZ$ let $U_i$ the root group belonging to the root $(x_i,\ldots,x_{i+4})$, i.e. 
the subgroup of $\Aut \Omega$ that stablises every vertex adjacent to $x_{i+1},x_{i+2}$ or $x_{i+3}$, 
and let $V_i$ be the pointwise stabiliser of the $2$-neighbourhood of $x_{i+2}$. Then $V_i\leq U_i$, 
and by \cite[(21.3)]{TW} we have for $i \in \ZZ/8\ZZ$
\begin{enumerate}
\item $V_i \leq Z(U_i)$.
\item $[U_i,U_{i+k}] \leq U_{i+1} \ldots U_{i+k-1}$ for $1 \leq k \leq 3$.
\end{enumerate}
Now let $W_i:=U_{i+1} U_{i+2} U_{i+3}$ and $M_i:= W_i/[W_i,W_i]$. Let $\bar{U}_j$ be the image 
of 
$U_j$ in $M_i$ for $j\in \{i+1,i+2,i+3\}$. Then we have $M_i=\bar{U}_{i+1} \oplus \bar{U}_{i+2} 
\oplus \bar{U}_{i+3}$ and 
 (b) implies $[\bar{U}_{i+3},U_i] \leq \bar{U}_{i+2}
\bar{U}_{i+2}, [\bar{U}{i+2},U_i] \leq \bar{U}_{i+1}$ and $[\bar{U}_{i+1},\bar{U}_i]=\{1\}$.
 Thus $M_i$ is a quadratic or cubic module for $G_i$. 
 One sees immediately that $M_i$ is quadratic if $V_i =U_i$. 
Therefore we only have to regard reduced and wide quadrangles.
\begin{itemize}
\item
\textbf{Reduced quadrangles:}\\
We may assume that $U_i =V_i$ for $i$ odd and $U_i \ne V_i$ for $i$ even. By \cite[(21.3)]{TW} 
$U_i$ is abelian for all $i$, so $M_0 = U_1 \times U_2 \times U_3$. By \cite[(21.8) and (21.9)]{TW}
 either $\Gamma \cong \mathcal{Q}_{\mathcal{I}} (K,K_0,\sigma)$ for an involutory set 
$(K,K_0,\sigma)$ with 
$\langle K_0 \rangle =K$ or $\Gamma \cong \mathcal{Q}_{\mathcal{Q}}(K,L_0,q)$ for an anistropic 
quadratic 
space $(K,L_0,q)$. For the definition of $\mathcal{Q}_{\mathcal{I}} (K,K_0,\sigma)$ and 
$\mathcal{Q}_{\mathcal{Q}}(K,L_0,q)$ compare \cite[(16.2) and (16.3)]{TW}.
\begin{enumerate}
\item If $\Gamma \cong \mathcal{Q}_{\mathcal{I}} (K,K_0,\sigma)$, then by \cite[(16.2)]{TW} we have 
$U_0 \cong U_2 \cong U_4 \cong (K,+)$ and $U_1 \cong U_3 \cong (K_0,+)$. Thus $G_0 \cong 
\PSL_2(K)$ 
and $M_0 \cong K_0 \oplus 
K \oplus K_0$, and the action of $U_4=\{\alpha(a)| a \in K\}$ is given by 
$$(s,b,t)\alpha(a)= (s,b-sa,t-a^{\sigma} s a-a^{\sigma} b -b^{\sigma} a)$$
for $a,b \in K$ and $s,t \in K_0$.
By \cite[(32.9)]{TW} we have 
$$(r,b,s)\mu = (s,-b^{\sigma},r)$$
for $\mu =\mu(\alpha(1))$. Since $G_0 =\langle A_4,\mu\rangle$, this determines the action of 
$G_0$ 
on $M_0$. The module $M_0$ can also described as follows. Set 
$$\mathcal{H}:=\left\{ \left( \begin{array}{cc} r & b \\ b^{\sigma} & s \end{array}\right)| 
r,s\in K_0, b\in K\right\}.$$
Then $G_0=\SL_2(K)$ acts on $\mathcal{H}$ by $x^y:=(y^{\sigma})^t x y$ for $x \in \mathcal{H}$ and 
$y\in G_0$ 
and $\mathcal{H}$ and $M_0$ are isomorphic as $G_0$-modules.
\item If $\Gamma \cong \mathcal{Q}_{\mathcal{Q}}(K,L_0,q)$ for an anisotropic quadratic space 
$(K,L_0,q)$ with bilinear form $f$, 
then by \cite[(16.3)]{TW} 
we have $U_0 \cong U_2 \cong U_4 \cong (L_0,+)$ and 
$U_1 \cong U_3 \cong (K,+)$. We may assume that there is an element $e \in L_0^*$ with $q(e) 
=1$
(otherwise, replace $q$ by $\hat{q} = {1 \over {q(e)}} q$ for an arbitrary $e\in L_0^*$). 
Thus $G_0 \cong \PSL_2(L_0,C(q,e))$ with $C(q,e)$ the Clifford algebra of $q$ with basepoint 
$e$ 
and $M_0 \cong K \oplus L_0 \oplus K$, 
and the action of $U_4=\{\alpha(a)| a \in L_0\}$ is given by 
$$(s,b,t)\alpha(a)= (s,b-as,t-f(a,b)-sq(a))$$
for $a,b \in L_0$ and $s,t \in K$.
If $\mu =\mu(\alpha(e))$, then by \cite[(32.7)]{TW} we have
$$(s,b,t) \mu = (t,b-q(a)^{-1} f(a,b)b,s)$$
for $a,b \in L_0$ and $s,t \in K$. Therefore the action of $G_0$ on $M_0$ is determined. \\
The module can be described as follows. Let $\sigma: L_0 \to L_0: a \mapsto f(a,e) e -a$   
and set
$$\mathcal{H} :=\left\{ \left( \begin{array}{cc} r & b \\ b^{\sigma} & s \end{array}\right)| b\in 
L_0, r,s 
\in K \right\} \subseteq \mathrm{Mat}_2(C(q,e)).$$
Then $G_0 = \SL_2(L_0,C(q,e))$ acts on $\mathcal{H}$ by $x^y = (y^{\sigma})^t  x y$ for $x\in 
\mathcal{H}$
 and $y\in G_0$, and $\mathcal{H}$ and $M_0$ are isomorphic $G_0$-modules. This example is 
 isomorphic to the one in \ref{orth}.
\end{enumerate}
\item
\textbf{Wide quadrangles:}\\
By \cite[(21.6)]{TW} every wide Moufang quadrangle is an extension of a reduced Moufang 
quadrangle 
in the sense of \cite[(21.5)]{TW}. 
\begin{enumerate}
\item If $\Gamma$ is an extension of $\mathcal{Q}_{\mathcal{I}} (K,K_0,\sigma)$ for an involutory set 
$(K,K_0,\sigma)$ with $\sigma \ne 1$ and $K =\langle K_0 \rangle$, then by \cite[(21.11)]{TW}
there is a $K$-vector space $L_0$ and an anisotropic pseudoquadratic form $q:L_0 \to K/K_0$ 
relative 
to $(K,K_0,\sigma)$ with skew-hermitian form $f$ 
such that $\Gamma \cong \mathcal{Q}_{\mathcal{P}}(K,K_0,\sigma, L_0,q)$ (see \cite[16.5]{TW} for the definition of 
this polygon). We have $U_2 \cong U_4 \cong K$ and 
$$U_1 \cong U_3 \cong U_5 \cong T:=\{(a,t) \in L_0 \times K| q(a)-t \in K_0\}$$ with 
$$(a,t) \cdot (b,u) = (a+b, t+u+f(b,a)$$ for $a,b \in L_0$ and $t,u \in K$.  
By the commutator formulas given in \cite[(16.5)]{TW} we have $[W_0,W_0] =U_2$ and $[W_1,W_1] = 
[U_3,U_3] \cong u_3(\{(0,a+a^{\sigma})|a \in K\})$ with $u_3: T\to U_3$ an isomorphism. Thus 
$M_0 $ is a quadratic module for $G_0$ but $M_1$ is a cubic module for $G_1$. Moreover, one has 
that 
$C_{G_1}(M_1) $ is the image of $u_3(\{(a,t)| a \in L_0^{\perp},q(a)-t \in K_0\})$ in $M_1$. 
The action of $U_1$ and $\mu=\mu((0,1))$ is given by \cite[(16.5) and (32.9)]{TW}  
One easily sees that $(G_1,M_1/C_{G_1}(M_1))$ is isomorphic to $(G,V)$ in \ref{pqf}.  
\item If $\Gamma$ is an extension of $\mathcal{Q}_{\mathcal{Q}}(K,L_0,q)$, then one 
of the following holds:
\begin{enumerate}
\item There is a multiplication on $L_0$ such that either $L_0$ is a separable quadratic 
extension of $K$ or a quaternion division algebra with centre $K$. In both cases $q$ is 
the norm of $L_0$ over $K$. In this case $\Gamma \cong \mathcal{Q}_{\mathcal{P}}(L_0,K,
\sigma, M_0,\pi)$ for an anisotropic pseuso-quadratic space $(L_0,K,\sigma,M_0,\pi)$, 
where $\sigma$ is the generator of $\mathrm{Gal}(L_0|K)$ in the first case and the standard 
involution of $L_0$ in the second case. Here we get a cubic module which is a pseudo-quadratic 
space over $(L_0,K,\sigma)$ and already appears in \ref{pqf}.
\item $(L_0,K,q)$ is of type $\mathrm{E}_6,\mathrm{E}_7,\mathrm{E}_8$ or $F_4$ and there is an 
element $1\in L_0$ with $q(1)=1$, a $K$-vector space $X$ and maps $\cdot: X \times L_0 \to X, 
h: X\times X \to L_0$ and $\theta: X\times L_0 \to L_0$ such that $\Xi =(K,L_0,q,1,X,\cdot,h,\theta)$ is 
a quadrangular algebra and $\Gamma\cong \Omega_\Xi$ as defined in \ref{quadrangular algebras}. 
Hence $U_2 \cong (L_0,+) \cong U_4$ and $U_1,U_3$ and $U_5$ are isomorphic to the group $S=X\times K$ 
with $(a,s)\cdot (b,t) =(a+b,s+t+g(b,a))$ for $a,b\in X,s,t\in K$. The commutator relations show that we 
may identify $M_0$ with $L_0 \oplus X\oplus L_0$. For $i=1,5$ let $x_i:S\to U_i$ be an 
isomorphism. Then for $u,v\in L_0, 
a,x \in X$ and $t\in K$ we have
$$(u,x,v)x_1(a,t) =(u+\theta(v,a)+h(x,a)+tv, x+a\cdot v, v)$$ and 
$$(u,x,v)x_5(a,t) =(u,x+a\cdot u,v+\theta(u,a)+h(x,a)+tu).$$
If $(K,L_0,q)$ is of type $\mathrm{E}_8$, then $G$ is an algebraic 
group of type $\mathrm{E}_{7,1}^{66}$ and $M_0$ is its $56$-dimensional 
standard module. Thus if $\characteristic K \ne 2$, then this example 
already appears in \ref{algebraic groups}.
\end{enumerate}
\end{enumerate}
\end{itemize}


\section{Suzuki and Ree groups}
 If $G$ is a Suzuki group or a Ree group, then there is no cubic module for $G$. 
In case of the Suzuki groups, this follows by the fact that every Suzuki group contains a 
Frobenius group of order 
$20$ which has only one non-trivial irreducible module in characteristic $2$ which is not 
cubic.
If $G$ is a Ree group, then there is an element $g \in A$ with $o(g) =9$. So the minimal 
polynom of $g$ can 
not be $(X-1)^3 $ and $g$ cannot act cubically. 
\chapter{On the structure of a cubic module}
We continue to assume that $G =\langle A,B \rangle$ acts cubically on a 
$\ZZ G$-module $V$. 
We start with an easy but crucial observation.
\begin{lemma}\label{trilinear}
The map $V \times A \times A \to V: (v,a,b) \mapsto [v,a,b]$ is additive in 
every component. 
Moreover, for $a \in A$ we have $[v,a,b] = 0 =[v,b,a]$ for all $b\in A$ 
if and only if $a \in 
A_0$.\end{lemma}
\begin{proof}
Let $a,b,c \in A$ and $v \in V$. Then we have
\begin{align*}
0=[v,a,b,c] = v(a-1)(b-1)(c-1)=\\ v(abc-ab-ac-bc+a+b+c-1) =\\
v(abc-bc-a+1)-v(ab-b-a+1)-v(ac-c-a+1)=\\
[v,a,bc] -[v,a,b]-[v,a,c],\end{align*}
hence $[v,a,bc] =[v,a,b]+[v,a,c]$. Similarly,
\begin{align*} 0=[v,a,b,c] = v(abc-ab-ac-bc+a+b+c-1) = \\ 
v(abc-ab-c+1)-v(ac-c-a+1)+v(bc-c-b+1) =\\
[v,ab,c]-[v,a,c]-[v,b,c],\end{align*}
hence $[v,ab,c] =[v,a,c]+[v,b,c]$. It is clear that this map is additive in 
the first component. 
For $a \in A$ we have $[v,a,b] =0$ for all $b\in A$ iff $a \in 
C_A(V/C_V(A))$ and 
$[v,b,a]=0$ for all $b\in A$ iff $a\in C_A([V,A])$. Thus the last claim 
follows.
\end{proof}

\bigskip
In \cite{T2} it was shown that for a quadratic module one may assume that 
$C_V(G)=0$ and $[V,G]=0$. We will now explain why this is also true for 
cubic modules. 
\begin{lemma}\label{standard form} 
\begin{enumerate}
\item Set $W = [V,G,G,G]$. Then $W$ is a cubic or quadratic module for $G$ with $[W,G] =W$.
\item Suppose that $[V,G] =V$. Set $Z_0 =0$, let $Z_{i+1}$ be the preimage of $C_{V/Z_i}(G)$ in 
$V$.
 Then $V/Z_3$ is a cubic or quadratic module for $G$ with $C_{V/Z}(G) =0$. 
\end{enumerate}
\end{lemma}
\begin{proof}
\begin{enumerate}
\item By definition $W \ne 0$. Set $X = V/[W,G]$ and $N:=C_G(X)$. Then $G/N$ is nilpotent.
Thus $G=NA$ by \ref{faithful}(c), and we get $0 = [X,A,A,A] = [X,G,G,G] =[V,G,G,G]/[W,G] 
=W/[W,G]$. Thus the claim follows.
\item Since $[V,G] =V$, we have $V \ne Z_n$ for all $n$.
Set $M :=C_G(Z_4)$. Then again, we get $G=MA$ and $0=[Z_4,A,A,A] =[Z_4,G,G,G]$. Thus $[Z_4,G,G] \leq 
C_V(G)=Z_1$, $[Z_4,G] \leq C_{V/[Z_4,G,G]}(G) \leq C_{V/Z_1}(G) =Z_2$ and so $Z_4=Z_3$. Therefore 
$C_{V/Z_3}(G) =0$ and the claim follows. 
\end{enumerate}
\end{proof}

As a corollary we get
\begin{coro}\label{normal form} 
Suppose that $G$ is a rank one group acting cubically on a module $V$. Then there 
are $G$-submodules $V_1 < V_2\leq V$ such that $W:=V_2/V_1$ is  
a cubic or quadratic $G$-module with $[W,G]=G$ and $C_W(G)=0$.
\end{coro}
\begin{lemma}\label{vmu}
Let $G=\langle A,B\rangle $ a rank one group acting cubically on $V$.
Let $v \in C_V(B)$ and $a \in A^*$. Then we have
$$v\mu(a) = v+[v,a] +[v,a,\fatb(a)^{-1}].$$
\end{lemma}
\begin{proof}
We have \begin{align*}
v\mu(a) = v\fatb(a^{-1}) a \fatb(a)^{-1} = va\fatb(a)^{-1} = \\ (v+[v,a])\fatb(a)^{-1} = v+[v,a] 
+[v,a,\fatb(a)^{-1}].\end{align*}
\end{proof}
\begin{definition}\rm
Let $X$ be a group acting on a module $V$. We define $Z_2(V,X):=\{v \in V| [v,x,y]=0 \ \forall 
x,y \in X \}$.
\end{definition}
By definition $Z_2(V,X)/C_V(X)= C_{V/C_V(X)}(X)$ and $[V,X,X,X]=0$ iff $[V,X] \leq Z_2(V,X)$.
\begin{lemma}\label{qusubmodule}
Let $V$ be a cubic module for a rank one group $G$ with $C_V(G)=0$. 
Set $W:=(C_V(A) \cap Z_2(V,B)) +(C_V(B) \cap Z_2(V,A))$.
\begin{enumerate}
\item $W$ is a $G$-submodule of $V$.
\item If $W \ne 0$, then $G$ operates quadratically on $W$. 
\end{enumerate}
 \end{lemma}
\begin{proof} 
Set $W_0:= C_V(A) \cap Z_2(V,B)$ and $W_1 =C_V(B) \cap Z_2(V,A)$. 
Note that $W_0 \cap W_1 \leq C_V(\langle A,B \rangle) =C_V(G) =0$. 
 Let $v \in W_1$ and $a\in A^*$. Then $[v,a] \in 
C_V(A)$ because $v\in Z_2(V,A)$, $[v,a,\fatb(a)^{-1}] \in [V,B] \subseteq 
Z_2(V,B), v\in C_V(B) \subseteq Z_2(V,B)$ and $v\mu(a) \in W_1\mu(a)=
W_0\subseteq Z_2(V,B)$. Thus with \ref{vmu} we have
$$[v,a] = v\mu(a) -v-[v,a,\fatb(a)^{-1}] \in C_V(A) \cap Z_2(V,B)=W_0.$$
Therefore we have $[W_1,A] \leq W_0$. Now $G=\langle A,\mu \rangle$, $[W_0,A]=0$, 
$W_0\mu_a = W_1$ and $W_1\mu_a =W_0$, thus $W$ is $G$-invariant. Furthermore, we have 
$[W,A] =[W_0,A] + [W_1,A] \leq W_0 =C_W(A)$, therefore $W$ is a quadratic module for $G$ 
if $W\ne 0$.
\end{proof}

\bigskip
If $G$ is a rank one group 
with unipotent subgroups $A$ and $B$ which acts quadratically on a module $V$, 
then $[V,G] =V$ and 
$C_V(G)=0$ implies $V =[V,A] \oplus [V,B]=C_V(A) \oplus C_V(B)$ 
(see \cite[2.1(8)]{T2}). In the cubic case, there is a similar 
decomposition of $V$.
\begin{definition}\rm
Let $V$ be a cubic module for a rank one group $G$. 
\begin{enumerate}
\item We say that $V$ is 
in \textit{standard form} if $[V,G]=V$ and $C_V(G)=0$. 
\item $V$ is called \textit{squarefree} if 
there is no non-trivial submodule $W$ of $G$ such that $G$ acts 
quadratically on $W$ or on $V/W$.
\end{enumerate} 
\end{definition}
\begin{prop}\label{decompV}
Suppose that $V$ is a squarefree cubic module in standard form for a rank one group $G$. Then we have
$V =C_V(A) \oplus (Z_2(V,A) \cap Z_2(V,B)) \oplus C_V(B)$.

\end{prop}
\begin{proof}
As in \ref{qusubmodule}, set $W:= (C_V(A) \cap Z_2(V,B)) +(C_V(B) \cap 
Z_2(V,A))$. Since $V$ is squarefree, there is 
no quadratic $G$-submodule of $V$, hence $W=0$. Now we set
$V_0 =C_V(A)+ (Z_2(V,A) \cap Z_2(V,B)) + C_V(B)$. We show that this sum is 
direct. 
Let $v_1 \in C_V(A), v_2 \in Z_2(V,A) \cap Z_2(V,B), v_3 \in C_V(B) $ such 
that
$v_1 + v_2 +v_3 =0$. Then for all $a \in A$ we have
$$0=[v_1+v_2+v_3,a]=[v_1,a]+[v_2,a]+[v_3,a]=[v_2,a]+[v_3,a],$$
hence $[v_3,a]=-[v_2,a]\in C_V(A)$ and so $v_3 \in Z_2(V,A) \cap C_V(B)=0$. 
Thus 
$[v_2,a]=-[v_3,a]=0$ 
and so $v_2 \in C_V(A)\cap Z_2(V,B)=0$. It follows $v_1=0$ and so we get
$V_0 =C_V(A)\oplus (Z_2(V,A) \cap Z_2(V,B)) \oplus C_V(B)$.
\\ Let $v \in C_V(B)$ and $a \in A^*$. Then we 
have $v\mu(a) \in C_V(A)$ and $[v,a,\fatb(a)^{-1}] \in [V,B] \leq Z_2(V,B)$, 
therefore
$v\mu(a) -[v,a] =v+[v,a,\fatb(a)^{-1}] \in Z_2(V,A) \cap Z_2(V,B)$ by \ref{vmu} 
and thus 
$[v,a]=v\mu(a) -v-[v,a,\fatb(a)^{-1}] \in C_V(A)+(Z_2(V,A) \cap Z_2(V,B))$. 
Since $[Z_2(V,A),A] 
\leq C_V(A)$, 
we get $[V_0,A] \leq V_0$ and so $V_0$ is $A$-invariant. By the same 
argument, $[V_0,B]\leq 
V_0$ and
hence $G=\langle A, B \rangle$ stabilises $V_0$. \\
Now we show that $V = V_0$. Set $\oV =V/V_0$. Since $[V,A,A] \leq C_V(A)$, we have $[\oV,A,A]=0$. 
Moreover, because $[V,G]=V$, we also have $[\oV,G]=\oV$. Thus $\oV$ is a quadratic module for $G$, 
therefore $V=V_0$ since $V$ is squarefree.
\end{proof}

From now on, we will always assume that $V$ is in standard form. 
This is not an 
essential restriction since by \ref{standard form} every cubic module $V$ 
gives rise to a cubic module $W$ in standard form or 
to a quadratic module and since we already know the quadratic rank one groups.
\\
If $X$ is an abelian group, we say $\characteristic X =0$ if $X$ is torsion-free and 
$\characteristic X =p$ if $X$ has exponent
$p \in \mathbb{P}$.
\begin{lemma}\label{characteristic}
Let $V$ be a squarefree cubic module in standard form.
\begin{enumerate}
\item Either $V$ is an elementary-abelian $p$-group for a prime $p$ or $V$ is torsion-free and 
uniquely divisible. Thus  
$\characteristic V =0$ or $\characteristic V =p \in 
\mathbb{P}$.
\item If $A_0 \ne 1$, then $\characteristic V = \characteristic A/A_0 = \characteristic A_0$. 
If $\characteristic V 
=0$, then $A_0$ is uniquely divisible.
\end{enumerate}
\end{lemma}
\begin{proof}
Since $[V,G] =V$, we have $[V/W,G] =V/W$ for every $G$-submodule $W$ of $V$. This implies that 
either 
$V=W$ or $V/W$ is a cubic module for $G$. Suppose there is 
a prime $p$ such that $pV \ne V$. Then $X:=V/pV$ is a cubic module for $G$. Thus
$C_A(X)=1$ by \ref{faithful2}. 

The commutator map induces biadditive maps from
$A \times X/[X,A]$ to $X/C_X(A)$ and from $A \times [X,A]$ to $C_X(A)$. Let $A_1$ be the 
intersection of these kernels. Since we also have a biadditive map from $A_1 \times X $ to 
$C_X(A)$, it follows that $A/A_1$ and $A_1$ have both exponent $p$ and so $A$ has exponent at 
most $p^2$. 
Especially, this means that if $pV \ne V$, then $V =qV$ for all primes $q \ne p$.
Since $V$ has no quadratic quotient, by \ref{decompV} one sees that
$C_X(A) = (C_V(A)+pV)/pV, C_X(B) =(C_V(B)+pV)/pV$ and $Z_2(X,A) \cap Z_2(X,B) = ((Z_2(V,A)\cap 
Z_2(V,B))+pV)/pV$.
It follows that $A_1=A_0$, and so $A/A_0$ is elementary-abelian. 
Now the map $C_V(B) \times A/A_0 \times A/A_0 :(v,a,b) \mapsto [v,a,b]$ is triadditive and 
non-degenerate.
This shows that $C_V(B)$ and hence also $C_V(A)$ has exponent $p$. Moreover, the map
$(Z_2(V,A) \cap Z_2(V,B)) \times A/A_0 \to C_V(A):(v,A_0a) \mapsto [v,a]$ is biadditive and non-degenerate, 
thus 
$Z_2(V,A) \cap Z_2(V,B)$ and hence $V$ has exponent $p$. 
\\
Suppose that $V =pV$ for all primes $p$. 
For a prime $p$ set $V_p:=\{v\in V| pv =0\}$. Suppose that $V_p\ne 0$ and 
set $N:=C_G(V_p)$. If 
$N \leq Z(G)$, 
then $V_p$ is a cubic or quadratic module for $G$. Again by a commutator 
argument one sees that 
$A$ has exponent
at most $p^2$. But this implies $p^2 V =0$, a contradiction. So we have 
$G=NA$. But then we 
have
$0 \ne C_{V_p}(A) =C_{V_p}(G) \leq C_V(G) =0$, a contradiction. Hence 
$V_p=0$ for all primes 
$p$ and so 
$V$ is torsion-free and uniquely divisible. By applying the commutator 
maps, one sees that
$A_0$ and $A/A_0$ are also torsion-free. Since $A_0$ is the root group of 
a special rank $1$ 
group, $A_0$ 
is also uniquely divisible by \cite[I.(5.2)]{T1}.
 \end{proof} 
 
 \bigskip
 We will later see that if $\characteristic V =0$, then also $A/A_0$ is 
 uniquely divisible 
(see \ref{module}).\\
If the quadratic kernel is non-trivial, we have further information about 
the structure of $V$.
\begin{prop}\label{decompV2}
Suppose that $V$ is squarefree and in standard form.
If $1 \ne A_0$, then the following hold.
\begin{enumerate}
\item $C_V(A) = [V,A_0] = [V,a]$ for all $a \in A_0^*$. 
\item $[V,A]=C_V(a)$ for all $a \in A_0^*$. 
\item $V=C_V(A)\oplus C_V(G_0) \oplus C_V(B)$.
\item $Z_2(V,A) \cap Z_2(V,B) =[V,A] \cap [V,B] =C_V(G_0)$.
\item $Z_2(V,A)=[V,A] = C_V(A) \oplus C_V(G_0)$.
\end{enumerate}
\end{prop}
\begin{proof}
\begin{enumerate}
\item Let $a \in A_0^*$ and set $W:= [V,a] + [V,B]$. Since $[V,B] \leq 
Z_2(V,B)$ and
$[V,a] \leq C_V(A)$, this sum is direct. 
Now both $B$ and $\langle a \rangle$ stabilise $W$ and act trivially on $V/W$. 
Since $G =\langle a,B \rangle$ 
and 
$[V,G]=V$, we get $V=W$. 
Hence if $v \in C_V(A)$, then there are $v_1 \in [V,a] \leq C_V(A)$ and 
$v_2 \in [V,B]$ 
with $v=v_1+v_2$, thus $v_2 = v-v_1 \in C_V(A) \cap Z_2(V,B)=0$ and so 
$[V,a] =C_V(A)$. Since $[V,a] \leq [V,A_0] \leq C_V(A)$, the other 
equation follows.
\item If $a \in A_0^*$, then by definition $[V,A] \leq C_V(a)$. Suppose 
there is a $v \in 
C_V(a) \setminus [V,A]$. 
Now the proof of (a) implies \begin{align*}
V=[V,B] \oplus C_V(A) =([V,B] \oplus C_V(A))\mu(a) =\\
[V,B]\mu(a) \oplus C_V(A)\mu(a)= [V,A]\oplus C_V(B).\end{align*}
 Thus we get $W:=C_V(B) \cap C_V(a) \ne 
0$. 
But then $G = \langle B, a \rangle \leq C_G(W)$, a contradiction.

\item Set $V_0 = [V,A_0] + [V,B_0]$. Then $V_0 = [V,G_0]$ and $[V,A_0] 
\cap [V,B_0]\leq 
C_V(A)\cap C_V(B)
=C_V(G)=0$. 
If $\characteristic V \ne 2$, then by \cite[4.2(a)]{T2} we get $V = 
C_V(G_0) \oplus V_0$. 
If  $\characteristic V = 2$ and $a \in A_0^*$ and $b=\fatb(a)$, then $\langle 
a,b\rangle 
\cong S_3$ by \ref{S3}, so $t:=ab$ has order $3$. 
If $v \in C_V(t)$, then $v+[v,a] = va =vt^{-1}a=vba^2= vb = v+[v,b]$ and 
thus $$[v,a] = [v,b] \in [V,A_0] \cap 
[V,B_0] = 0.$$ 
Thus $C_V(t) = C_V(a) \cap C_V(b) =C_V(G_0)$ by (a). 
If $v \in C_V(G_0) \cap V_0 =([V,a] \oplus [V,b]) 
\cap C_V(G_0)$, then there are $v_1,v_2 \in V$ with $v=[v_1,a] +[v_2,b]$. 
Hence we have $[v_2,b]=v-[v_1,a] \in C_V(G_0) +[V,a] \leq C_V(a)$, thus 
$[v_2,b] \in C_V(a) \cap [V,b] = [V,A] \cap C_V(B) =0$. By the same 
argument, we get $[v_1,a]=0$ and thus 
$v=0$, so $C_V(G_0) \cap V_0 =0$. 
Since $[V,t] \leq [V,G_0]=V_0$, $C_V(t) =C_V(G_0)$ and $V=[V,t] \oplus 
C_V(t)$,  
we again get $V=V_0 \oplus C_V(G_0)$. \\
Thus in all cases we have $$V = C_V(G_0) \oplus V_0 = C_V(G_0) \oplus 
[V,A_0] \oplus [V,B_0]
= C_V(A) \oplus C_V(G_0) \oplus C_V(B).$$ 
\item Since $G_0 = \langle A_0, B_0 \rangle$, we get $C_V(G_0) =C_V(A_0) \cap 
C_V(B_0) = [V,A] \cap 
[V,B]$. The equation $Z_2(V,A)\cap Z_2(V,B)=C_V(G_0)$ follows by \ref{decompV} and (c).
\item 
Since $C_V(A) + C_V(G_0) \leq Z_2(V,A)$ and $V=C_V(B) \oplus Z_2(V,A)$, we 
get $C_V(A) + C_V(G_0)=[V,A]=Z_2(V,A)$ by (a) and (b). Moreover, $C_V(A)\cap C_V(G_0) 
\leq C_V(G)=0$, so 
the claim follows. 
\end{enumerate}
\end{proof}
\begin{coro}\label{mu quadrat}
Suppose that $V$ is squarefree and in standard form and $a\in A_0^*$.
\begin{enumerate}
\item $\mu(a)^2 $ inverts every element in $C_V(A) \oplus C_V(B)$ and 
centralises 
$C_V(G_0)$
\item $\mu(a)^2  \in Z(H)$.
\end{enumerate}
\end{coro}
\begin{proof}
\begin{enumerate}
\item Since $C_V(A) \oplus C_V(B)$ is a quadratic module for $G_0$, this is \cite[3.2]{T2}.
\item This follows from the fact that the decomposition in \ref{decompV2} is $H$-invariant.
\end{enumerate}
\end{proof}
From now on we assume that $V$ is squarefree and in standard form.
   Let $e \in A^*$ be fixed. 
 If $A_0 \ne 1$, then we assume $e \in A_0^*$. If $V$ has characteristic $2$ and 
 $A$ is not abelian, we choose $e$ in such a way that $e = a^2$ for an element $a \in A^*$.
 We set $\mu:=\mu(e^{-1})$. For 
 $a \in A^*$ set $h(a):= \mu \mu(a)$ and $h(1) :=0$. 
 Let $\varrho:H \cup \{0\} \to \End(C_V(A)): h \mapsto h|_{C_V(A)}$. For $a\in A$ we set 
 $\varrho(a)=\varrho(h(a))$.  

 By \cite[3.7]{T2}, we have that 
 $J:=\{\varrho(a)| a \in A_0\}$ is a special quadratic Jordan 
 division algebra 
 in $\End(C_V(A))$. Let $R$ be the subring of $\End(C_V(A))$ generated by 
 $J$ and $S$ the 
 subring of 
 $\End(C_V(A))$ generated by $\varrho(H)$. Notice that $R$ is generated by 
 $\varrho(H_0)$. Since 
 $\varrho(\mu^ 2) 
 =-1$ by \ref{mu quadrat}, for all $r \in R$ there are 
 $n\in \NN$ and elements $h_1,\ldots,h_n \in 
 H_0$ with 
 $r =\sum_{i=1}^n \varrho(h_i)$. 
 \begin{lemma}\label{rho}
 Let $v \in C_V(A)$ and $a \in A^*$. Then we have
  \begin{enumerate}
 \item $v\varrho(a) -[v\mu,a] = v\mu + [v\mu,a,\fatb(a^{-1})] \in Z_2(V,A) \cap Z_2(V,B)$.
 \item If $a \in A_0^*$, then $v\varrho(a) = [v\mu,a]$. 
 \item $\varrho({ab})=\varrho(a) + \varrho(b)$ for all $a,b \in A_0$. 
 \end{enumerate}
 \end{lemma}
 \begin{proof}
 \begin{enumerate}
 \item 
 We have $v\varrho(a)=vh(a) =v\mu \mu(a)$. 
 Since $v\mu \in C_V(A)\mu =C_V(B)$, the claim follows with
 \ref{vmu}.
 \item If $a \in A_0^*$, then $v\varrho(a)-[v\mu,a]=
 vh(a) -[v\mu,a] \in C_V(A) \cap Z_2(V,B) = 0$, since we assume that 
 $V$ is in standard form. For $a=1$, $v\varrho(1)=
 vh(1) = 0 =[v\mu,1]$, thus the 
 claim holds in every case.
 \item For $v \in C_V(A), a,b \in A_0$ we have 
 $$v\varrho({ab}) = vh({ab}) = [v\mu,ab] = 
 [v\mu,a] + [v\mu,b] = vh(a) +vh(b) = v\varrho(a) +v\varrho(b)$$ and 
 thus $\varrho({ab})= \varrho(a) + \varrho(b)$. 
 \end{enumerate}
 \end{proof}
 \begin{lemma}\label{-mu}
 Let $a\in A$ and $h\in H$.
 \begin{enumerate}
 \item 
 $\varrho(h(a)^{-\mu})= -\varrho(a^{-1})$.
 \item If $a\in A_0$, then $\varrho(h(a)^{-\mu}) =\varrho(a)$. 
 \item $h({a^h}) = h^{-\mu} h(a) h$.   
 \item $h^{-\mu} h \in H_0$.
 \end{enumerate}
 \end{lemma}
 \begin{proof}
(a)-(c) is trivially true for $a=1$, so we may assume that $a\in A^*$. 
 \begin{enumerate}
 \item 
 We have
 $$ h(a)^{\mu} = \mu \mu(a)^{\mu} = \mu(a)\mu = h(a^{-1})^{-1} \mu^2.$$ 
 Since $\varrho(\mu^2)=-1$ by \ref{mu quadrat}, we get 
 $$\varrho(h(a)^{-\mu}) =-\varrho(h(a^{-1}))=-\varrho(a^{-1}).$$
  \item If $a \in A_0^*$, then for all $v\in C_V(A)$ we have
 $$-v\varrho(a)=-[v\mu,a] =[v\mu,a^{-1}] =\varrho(a^{-1})$$ by \ref{rho}(b). 
   Thus the claim follows with part (a).
 \item 
 We have $\mu(a^h) = \mu(a)^h = h^{-1} \mu(a) h$ and thus $$h(a^h) = \mu h^{-1} \mu(a) h 
 = h^{-\mu^{-1}} \mu \mu(a) h = h^{-\mu^ {-1}} h(a) h.$$ 
 By \ref{mu quadrat})(b) $\mu^2 \in Z(H)$, so the claim follows.
 \item This follows by (c) with $a=e$.
 \end{enumerate}
 \end{proof}

  \begin{prop}\label{propf} Set $f: A \times A \to \End(C_V(A)):(a,b) \mapsto (v \mapsto 
  [v\mu,a,b])$. Then
 \begin{enumerate} 
 \item $f$ is biadditive.
 \item $f(a^h,b^h) =\varrho(h^{-\mu})f(a,b)\varrho(h)$ for all $a,b\in A$ and $h\in H$.
 \item $f(a,b^h) = f(a,b) \varrho(h)$ and $f(a^h,b) = \varrho(h^{-\mu})f(a,b) $ for all $a,b \in A$ 
 and all $h 
 \in H_0$.
 \item $C_A(V/C_V(A))$ is the left kernel of $f$ and $C_A([V,A])$ is the right kernel of $f$.
 \item $f(a,b) -f(b,a) = \varrho([a,b])$ for all $a,b \in A$. 
 \item If the characteristic of $V$ is $2$, then $f(a,a) = \varrho({a^2})$.
 \end{enumerate}
 \end{prop}
 \begin{proof}
 \begin{enumerate}
 \item This follows from \ref{trilinear}. 
 \item For all $v\in C_V(A)$ we have
 \begin{align*}
 vf(a^h,b^h) =v\mu(a^h-1)(b^h-1)=v\mu h^{-1}(a-1)(b-1)h =\\vh^{-\mu^{-1}}\mu(a-1)(b-1)h=
 vh^{-\mu}f(a,b)h =v\varrho(h^{-\mu})f(a,b)\varrho(h).
 \end{align*}
 Here we have used that $\mu^2 \in Z(H)$ by \ref{mu quadrat}(b).
  \item We first claim that 
 $$v\mu(a-1)h(b-1) =vf(a,b)$$
 for all $a,b \in A$, all $v \in C_V(A)$ and all $h \in H_0$ holds. 
 We have $v\mu(a-1) = [v\mu,a] = [v\mu,a] -vh(a) +vh(a)$. Since $[v\mu,a] -vh(a) \in Z_2(V,A) 
 \cap Z_2(V,B)\leq C_V(h)$, we  get 
 $$v\mu(a-1)h = [v\mu,a] -vh(a) +vh(a)h$$ and thus
 $$v\mu(a-1)h(b-1) = [v\mu,a](b-1) + (vh(a) h -vh(a))(b-1) = [v\mu,a,b] =vf(a,b),$$
 since $vh(a) h -vh(a) \in C_V(A)$.
Hence we have $$vf(a,b^h) =v\mu (a-1)(b^h-1)=v\mu (a-1) h^{-1}(b-1) h =v f(a,b)h$$ and so 
$f(a,b^h)=f(a,b)\varrho(h)$. 
Moreover, we get using (b)
 $$f(a^h,b)=\varrho(h^{-\mu})f(a,b^{h^{-1}}))\varrho(h)=\varrho(h^{-\mu})f(a,b)\varrho(h^{-1})\varrho(h)=
\varrho(h^{-\mu})f(a,b).$$ 
 \item Suppose that $a \in A$. Then $f(a,b) =0$ for all $b \in A$ iff $[v\mu,a] \in C_V(A)$ for 
 all 
 $v \in C_V(A)$. Therefore $f(a,b) =0 $ for all $b\in A$ iff $[C_V(B),a] \subseteq C_V(A)$. 
 Since 
 $V =Z_2(V,A) \oplus C_V(B)$ and $[Z_2(V,A),a] \leq C_V(A)$ by definition, this holds 
 iff 
 $[V,a] \subseteq C_V(A) $ and thus iff $a \in C_A(V/C_V(A))$.
 \newline Similarly, $f(b,a) =0$ for all $b \in A$ holds iff $[C_V(B),A] \subseteq C_V(a)$. 
 Since 
 $V = C_V(B) \oplus Z_2(V,A)$ and $[Z_2(V,A),A] \subseteq C_V(A)$, this holds iff $[V,A] 
 \subseteq 
 C_V(a)$ 
 and thus iff $a \in C_A([V,A]).$
 \item Since $[a,b] \in A^{\prime} \subseteq A_0 \subseteq Z(A)$ and $[v\mu, [a,b]] \in [V,A_0] =
 C_V(A)$, we get by \ref{rho}(b)
 \begin{align*} vf(a,b) -vf(b,a) = v\mu (a-1)(b-1) -v\mu(b-1)(a-1) = v\mu(ab -ba) =\\
  v\mu ba(a^{-1} b^{-1}ab -1) =v\mu b a ([a,b] -1) = v\mu ([a,b] -1) ba =\\
  [v\mu, [a,b]] ba = [v\mu, [a,b]]=v\varrho([a,b]).\end{align*}
  Thus the claim follows.
 \item For all $a \in 
 A$ and all $v \in C_V(A)$ we have
 $v\mu f(a,a) =v\mu (a-1)^2 = v\mu (a^2 -1) = [v\mu,a^2] = v\varrho(a^2).$
 \end{enumerate}
 \end{proof}

 Set $\oA = A/A_0$ and $\oa:=A_0 a$ for $a\in A$. By the previous proposition, we may consider
 $f$ as a 
 biadditive map from $\oA \times \oA$ to $\End(C_V(A))$. We will use the additive notation for 
 $\oA$, so 
 $\oa +\ob$ means 
 $\overline{ab}$ for $a,b \in A$ and $\overline{0}$ means the neutral element in $\oA$.  
 \newline The map $f$ is linked to the map $a \mapsto \varrho(a)$ from $A$ to $S$.
 \begin{prop}\label{handf} Let $a,b\in A$. Then we have
 \begin{enumerate}
 \item $f(a,b) = \varrho(ab) -\varrho(a) -\varrho(b)$.  
 \item $f(a,a)=\varrho(a)+\varrho(a^{-1})=\varrho(a)-\varrho(h(a)^{-\mu})$.
 \item $\varrho((h(ab)^{-\mu})=\varrho(h(a)^{-\mu})+\varrho(h(a)^{-\mu})-f(b,a)$.
 \end{enumerate}
 \end{prop}
 \begin{proof}
 \begin{enumerate}
 \item For $v \in C_V(A)$ we have by \ref{rho}(a)
 $$v\varrho({ab}) -[v\mu,ab] -(v\varrho(a) -[v\mu,a]) - (v\varrho(b) -[v\mu,b]) \in Z_2(V,A) \cap Z_2(V,B).$$
 Thus \begin{align*} v\varrho(ab) -v\varrho(a) -v\varrho(b) -v\mu(ab-1) +v\mu(a-1) +v\mu(b-1) 
 =\\ v\varrho({ab}) -v\varrho(a) -v\varrho(b) -v\mu(ab - a -b +1) =\\
 v\varrho({ab}) -v\varrho(a) -v\varrho(b) -vf(a,b) 
 \in Z_2(V,A) \cap Z_2(V,B).\end{align*}
 But we also have $v\varrho({ab}) -v\varrho(a) -v\varrho(b) -vf(a,b) \in C_V(A).$
 Since $C_V(A) \cap Z_2(V,B)= 0$, the claim follows. 
\item This follows by (a) and \ref{-mu}(a) with $b=a^{-1}$.
\item By (a), \ref{-mu}(a) and the biaddititivy of $f$ we have
\begin{align*}
\varrho(h(ab)^{-\mu}) =-\varrho((ab)^{-1}) =-\varrho(b^{-1}a^{-1})=\\
-\varrho(b^{-1})-\varrho(a^{-1})-
f(b^{-1},a^{-1})=\varrho(h(b)^{-\mu})+\varrho(h(a)^{-\mu})-f(b,a).
\end{align*} 
 
 \end{enumerate}
  
 \end{proof}
 
 \bigskip
 We immediately get
 \begin{coro}\label{coro handf} $f(a,b) \in S$ for all $a,b \in A$. 
 \end{coro}

 
 \bigskip
 We define  
$$\Phi: C_V(A) \times A \to Z_2(V,A) \cap Z_2(V,B): (v,a) \mapsto  
[v\mu,a] -vh(a).$$
This is well-defined by \ref{rho}. 
This important map reveals the connection between $A$ and $Z_2(V,A) \cap 
Z_2(V,B)$. 
\begin{lemma}\label{eigenschaftenphi}
 \begin{enumerate}
\item 
 $\Phi$ is a biadditive map from $C_V(A) \times A$ to $Z_2(V,A) \cap 
 Z_2(V,B)$.
 \item $\Phi(vh,a) = \Phi(v,a^{h^{-\mu}})$ for all $v \in C_V(A), a \in 
 A$ and $h \in H_0$. 
 \item $C_A(V/C_V(A))$ is the right kernel of $\Phi$. 
\end{enumerate}
\end{lemma}
\begin{proof} 
\begin{enumerate}
\item 
It is clear that $\Phi(v+w,a) =\Phi(v,a) + \Phi(w,a)$ holds for all $v,w 
\in C_V(A)$ and all 
$a \in A$. If $v \in C_V(A)$ and $a,b \in A$, 
then \begin{align*}
\Phi(v,ab) = [v\mu,ab] -vh({ab}) = -v\mu +v\mu ab -vh(a) -vh(b) -vf(a,b) 
=\\
-v\mu +(v\mu +[v\mu,a])b -vh(a) -vh(b)-vf(a,b)= \\ -v\mu + v\mu b 
+[v\mu,a] +[v\mu,a,b] -vh(a)
 -vh(b)-[v\mu,a,b] =\\
[v\mu,a] -vh(a) +[v\mu,b] -vh(b) +[v\mu,a,b] -[v\mu,a,b] =\\
[v\mu,a] -vh(a) +  [v\mu,b] -vh(b) =\Phi(v,a) + \Phi(v,b).\end{align*}
\item For $v \in C_V(A), h\in H_0$ and $a \in A$ we have $\Phi(vh,a) \in 
Z_2(V,A)\cap Z_2(V,B) \subseteq 
C_V(h^{-\mu})$. Thus we 
get by \ref{mu quadrat}(b) and \ref{-mu}(c)  \begin{align*} \Phi(vh,a)=
[vh\mu,a]-v h h(a) = ([v\mu h^{\mu},a] -vhh(a))h^{-\mu} =\\
-v\mu h^{\mu} h^{-\mu} + v\mu h^ {\mu} a h^ {-\mu}
-vhh(a) h^ {-\mu} =\\
 -v\mu + v\mu a^{h^{-\mu}} -vh({a^{h^ {-\mu}}}) = [v\mu,a^{h^{-\mu}}]-
 vh({a^{h^{-\mu}}})
   = \Phi(v,a^{h^{-\mu}}).
\end{align*}
\item Suppose $a \in A$ with $[v\mu,a] - vh(a) =0$ for all $v \in 
C_V(A)$. Then we have 
$$ [V,a] = [[V,A] \oplus C_V(B),a] \leq C_V(A) +[C_V(A)\mu,a]\leq C_V(A) 
$$
and thus $a \in C_A(V/C_V(A))$. 
If $a \in C_A(V/C_V(A))$, then 
$$[v\mu,a] - vh(a) \in C_V(A) \cap Z_2(V,A)\cap Z_2(V,B) =C_V(A) \cap Z_2(V,B) =0.$$ Thus the claim follows. 
\end{enumerate}
\end{proof}

\bigskip
The right kernel of $\Phi$ contains $A_0$ by (c). Therefore, if $A_0\ne 1$, 
we may regard $\Phi$ as a map from 
$C_V(A) \times \oA$ to $Z_2(V,A)\cap Z_2(V,B) =C_V(G_0)$.

\begin{lemma}\label{f b a mu}
 For all $a,b \in A^*$ we have $f(b\diamond \mu,a\diamond \mu) =\varrho(b)^{-1} f(b,a) \varrho(h(a^{-1})^{-1})$.
 \end{lemma}
 \begin{proof}

By \ref{h+}, \ref{-mu}(a) and \ref{handf}(a) we have
  \begin{align*}
 \varrho(b\diamond \mu(a\diamond \mu)^{-1}) =\varrho(b)^{-1}\varrho(ba^{-1})
 \varrho(h(a)^{\mu})=\\-\varrho(b)^{-1} (\varrho(b)+\varrho(a^{-1})
 +f(b,a^{-1}))\varrho(a^{-1})^{-1}=\\
 -\varrho(a^{-1})^{-1}-
 \varrho(b)^{-1}+\varrho(b)^{-1}f(b,a)\varrho(a^{-1})^{-1}.
 \end{align*}
 On the other hand, we have by \ref{mu a tau}(b), \ref{mu quadrat} and 
 \ref{-mu}(a)
 \begin{align*}
\varrho(h((a\diamond \mu)^{-1})) = -\varrho(h(a\diamond \mu)^{-\mu})=
-\varrho((\mu^2 h(a)^{-1})^{-\mu})=-\varrho(h(a)^\mu \mu^{-2}) =\\
\varrho(h(a)^\mu) = -\varrho(a^{-1})^{-1}.
\end{align*}
Thus we get again by \ref{mu a tau}(b)
  \begin{align*}
  \varrho((b\diamond \mu)(a\cdot 
  \mu)^{-1})=\varrho(b\diamond \mu)+\varrho((a\diamond \mu)^{-1})+f(b\diamond 
  \mu, (a\diamond \mu)^{-1})=\\
  \varrho(\mu^2 h(b)^{-1}) -\varrho(a^{-1})^{-1}-f(b\diamond \mu,a\diamond \mu)=
  -\varrho(b)^{-1}-\varrho(a^{-1})^{-1}-f(b\diamond \mu,a\diamond \mu).
  \end{align*}
 Thus we get
 $$f(b\diamond \mu,a\diamond \mu)=-\varrho(b)^{-1}f(b,a)\varrho(a^{-1})^{-1}.$$ 
 \end{proof}

 \begin{lemma}\label{f a mu} Suppose that $A_0 \ne 1$. For all $a,b \in A$ we have 
 \begin{enumerate}
 \item $\Phi(v,a) =\Phi(vh(a),a\diamond\mu)$.
 \item
 $f(a,b) = \varrho(a) f(a\diamond \mu,b))-f(a,b\diamond \mu)\varrho(b^{-1})$.
 \item 
 $\Phi(v,a^{\sim}) = \Phi(vh(a)^{-\mu} h(a)^{-1},a)$.
 \item
 $f(a^{\sim},b) = \varrho(h(a)^{-\mu}) \varrho(a)^{-1} f(a,b)$.
  \end{enumerate}
 \end{lemma}
  \begin{proof}
 \begin{enumerate}
 \item We have $\Phi(v,a) = [v\mu,a] -vh(a) = -v\mu-[v\mu,a,\fatb(a)^{-1}]$ by \ref{vmu}. 
 The element $\mu\in 
 G_0$ fixes $\Phi(v,a) \in Z_2(V,A)\cap Z_2(V,B)=C_V(G_0)$. Since $\fatb(a)=(a\diamond 
 \mu)^{\mu^{-1}}$ by \ref{b(a)}, we get
 \begin{align*}
  \Phi(v,a)=\Phi(v,a)\mu = -v\mu^2-[v\mu,a,\fatb(a)^{-1}]\mu = \\
v-[\Phi(v,a),\fatb(a)^{-1}]\mu-[vh(a),\fatb(a)^{-1}]\mu=\\
v-[\Phi(v,a)\mu^{-1},\fatb(a)^{-1}]\mu-[vh(a)\mu,\fatb(a)^{-\mu}] =\\ v-[\Phi(v,a),
\fatb(a)^{-\mu}]-
[vh(a),(a\diamond \mu)^{-1}]=\\  v-vf(a,(a\diamond \mu)^{-1})-\Phi(vh(a),(a\diamond 
\mu)^{-1})-vh(a) h({(a\diamond \mu)^{-1}})
=\\v(1+f(a,a\diamond \mu)+h(a) 
h(a\diamond \mu)^{-1})+\Phi(vh(a),a\diamond \mu).
 \end{align*}
 Thus we get
 $$\Phi(v,a) -\Phi(vh(a),a\diamond \mu) =v(1+f(a,a\diamond \mu) -h(a) h(a\diamond \mu)^{-1}) 
 \in C_V(G_0) \cap C_V(A) =0.$$
 Hence the claim follows.
 \item
 We have $$vf(a,b) = [\Phi(v,a),b] = [\Phi(vh(a),a\diamond \mu),b) =vh(a) f(a\diamond \mu,b)=v\varrho(a)f(a,b)$$
 for all $v \in C_V(A)$. Hence the first equation follows. For the 
 second equation we apply \ref{f b a mu} and get
$$\varrho(a)^{-1} f(a,b\diamond \mu)=
f(a\diamond \mu,b\diamond \mu) =-\varrho(a)^{-1} f(a,b)\varrho(b^{-1})^{-1}$$ 
and thus $f(a,b\diamond \mu)=-f(a,b)\varrho(b^{-1})^{-1}$.
  \item We have by applying (a) and \ref{-mu}(a)
\begin{align*} 
 \Phi(v,a^{\sim}) = \Phi(v,(a\diamond \mu^{-1})^{-1}\diamond \mu)=
 \Phi(vh((a\diamond \mu^{-1})^{-1})^{-1},(a \diamond \mu^{-1})^{-1})=\\
 -\Phi(-vh(a\diamond \mu^{-1})^\mu,a\diamond \mu^{-1})=\Phi(vh(a\diamond \mu^{-1})^\mu,
 a\diamond \mu^{-1}).
 \end{align*}
By \ref{mu a tau}(b) we have $h(a\diamond \mu^{-1}) =\mu^2 h(a)^{-1}$.
Applying (a) again and using the fact that $\varrho(\mu^2)=-1$ we get
 \begin{align*}
 \Phi(v,a^\sim) =\Phi(v\mu^2 h(a)^{-\mu},a\diamond \mu^{-1}) =
 \Phi(-vh(a)^{-\mu} h(a\diamond \mu^{-1}),a)=\\
 -\Phi(vh(a)^{-\mu}\mu^2 h(a)^{-1},a)=
 \Phi(vh(a)^{-\mu}h(a)^{-1},a).
 \end{align*}
 \item We have by (c)
 $$vf(a^{\sim},b) = [\Phi(v,a^{\sim}), b] = [\Phi(vh(a)^{-\mu} h(a)^{-1},a),b] = vh(a)^{-\mu} h(a)^{-1} f(a,b).$$ 
 Thus the claim follows.
 \end{enumerate}
 \end{proof}
 
  \begin{prop}\label{Phi h} Suppose that $A_0\ne 1$. Then for $a,b \in A^*$ and $v \in C_V(A)$ we have
 \begin{enumerate}
 \item $\Phi(v,a^{h(b)}) = \Phi(vh(b)^{-\mu},a) -\Phi(vh(b)^{-\mu} f(a,b)h(b)^{-1}, b)$.
 \item $\Phi(v,a)h(b) = \Phi(v,a) -\Phi(vf(a,b)h(b)^{-1},b)$.
 \end{enumerate}
 \end{prop}
 \begin{proof}
 \begin{enumerate}
 \item By \ref{formulas}(a) we have $((a\diamond \mu^{-1}) (b\diamond \mu^{-1})^{-1})\diamond \mu = (ab^{-1}) \diamond \mu(b) b^\sim$. 
 Thus we have by \ref{f a mu}(a) and (c) \begin{align*}
 \Phi(v, (ab^{-1})\diamond \mu(b)) =  \Phi(v,(((a\diamond \mu^{-1})
 (b\diamond \mu^{-1})^{-1})\diamond \mu) (b^\sim)^{-1} )= \\
 \Phi(vh((a\diamond \mu^{-1})(b\diamond \mu^{-1})^{-1}),a\diamond \mu^{-1}
 (b\diamond \mu^{-1})^{-1})-\Phi(vh(b)^{-\mu} h(b)^{-1},b).
 \end{align*}
 Now by \ref{h+} we have
 $h((a\diamond \mu^{-1})(b\diamond \mu^{-1})^{-1})=h(a)^{-1}h(ab^{-1})h(b)^\mu$.
 Thus we get again by \ref{f a mu}(a)
 \begin{align*}
 \Phi(v,(ab^{-1})\diamond \mu(b)) = \\
 \Phi(vh(b)^{-\mu} h({ab^{-1}})^{-1} h(a),(a\diamond \mu^{-1})(b\diamond \mu^{-1})^{-1}) -\Phi(vh(b)^{-\mu} h(b)^{-1},b)=\\
 -\Phi(vh(b)^{-\mu} h({ab^{-1}})^{-1},a)
  +\Phi(vh(b)^{-\mu} h({ab^{-1}})^{-1} h(a) h(b)^{-1},b)
  \\ -\Phi(vh(b)^{-\mu} h(b)^{-1},b).
 \end{align*}
 If we replace $a$ by $ab$ and apply \ref{handf}(a) and \ref{f a mu}(a) and 
 (b), we get
 \begin{align*}
 \Phi(v,a\diamond \mu(b)) =
  -\Phi(vh(b)^{-\mu} h(a)^{-1},a)-
 \Phi(vh(b)^{-\mu} h(a)^{-1},b)+\\ \Phi(vh(b)^{-\mu} h(a)^{-1} h({ab}) h(b)^{-1},b)-\Phi(vh(b)^{-\mu} h(b)^{-1},b) =\\
  -\Phi(vh(b)^{-\mu} h(a)^{-1},a)-
 \Phi(vh(b)^{-\mu} h(a)^{-1},b)+\\ \Phi(vh(b)^{-\mu} h(a)^{-1} h({a}) h(b)^{-1},b)+
 \Phi(vh(b)^{-\mu} h(a)^{-1} h({b}) h(b)^{-1},b)+\\
 \Phi(vh(b)^{-\mu} h(a)^{-1} f(a,b) h(b)^{-1},b)- \Phi(vh(b)^{-\mu} h(b)^{-1},b) =\\
 -\Phi(vh(b)^{-\mu}h(a)^{-1},a) +\Phi(vh(b)^{-\mu} h(a)^{-1} f(a,b) h(b)^{-1},b)= \\
 -\Phi(vh(b)^{-\mu},a\diamond \mu) +\Phi(vh(b)^{-\mu} f(a\diamond \mu,b) h(b)^{-1},b).
 \end{align*}
 If we finally replace $a$ by $a\diamond \mu$, we get by \ref{mu quadrat} 
 and \ref{propf}(c)
 \begin{align*}
 \Phi(v,a^{h(b)}) = -\Phi(vh(b)^{-\mu},a\mu^2)
 +\Phi(vh(b)^{-\mu} f(a^{ \mu^2},b)h(b)^{-1},b)=\\-\Phi(vh(b)^{-\mu}\mu^2,a)+
 \Phi(vh(b)^{-\mu}\mu^2 f(a,b)h(b)^{-1},b)=\\
 \Phi(vh(b)^{-\mu},a)-\Phi(vh(b)^{-\mu}f(a,b)h(b)^{-1},b).\end{align*}
 \item 
 We have $h(a^{h(b)})=\mu \mu(a^{h(b)}) =\mu h(b)^{-1} \mu(a)h(b) =
 h(b)^{-\mu^{-1}} h(a) h(b) =h(b)^{-\mu^{-1}} h(a) h(b)$.
 Therefore we have by part (a) and \ref{mu quadrat}(b)
  \begin{align*}
 \Phi(v,a)h(b) = (vh(a)-[v\mu,a])h(b)= vh(a)h(b) -v\mu(a-1)h(b) =\\
 vh(b)^{\mu^{-1}} h(a^{h(b)})-v\mu h(b)(a^{h(b)}-1)=\\
 vh(b)^{\mu^{-1}} h(a^{h(b)})-v h(b)^{\mu^{-1}}(a^{h(b)}-1)=\\
 \Phi(vh(b)^{\mu^{-1}},a^{h(b)}) =\Phi(vh(b)^\mu,a^{h(b)} =
5 \Phi(v,a)-\Phi(vf(a,b)h(b)^{-1},b).
  \end{align*}
 \end{enumerate}
 \end{proof}
 \begin{coro}\label{f h} 
Suppose that $A_0\ne 1$. Then for $a,b,c\in A^*$ we have 
 \begin{enumerate}
 \item $f(a^{h(b)},c) = \varrho(h(b)^{-\mu}) f(a,c) -\varrho(h(b)^{-\mu}) f(a,b) \varrho(b)^{-1} f(b,c)$. 
 \item $f(a,c^{h(b)}) = f(a,c) \varrho(b) + f(a,b) \varrho(h(b)^{\mu}) 
 f(b,c) \varrho(b) $.
 \end{enumerate}
 \end{coro}
 \begin{proof}
 \begin{enumerate}
 \item For all $v\in C_V(A)$ we have by \ref{Phi h}(a)
 \begin{align*}
 vf(a^{h(b)},c) =[\Phi(v,a^{h(b)}),c] = [\Phi(vh(b)^{-\mu},a)-\Phi(vh(b)^{-\mu} f(a,b)h(b)^{-1},c] = \\vh(b)^{-\mu} f(a,c) 
 -vh(b)^{-\mu} f(a,b) h(b)^{-1} f(b,c).\end{align*} Thus the claim follows.
 \item 
We have by \ref{mu a tau}(b) and \ref{propf}(b)
\begin{align*}
f(a,c^{h(b)}) =\varrho(h(b)^{-\mu}) f(a^{h(b)^{-1}},c)\varrho(b)=\varrho(h(b)^{-\mu}) f(a^{\mu^2 h(b\diamond \mu)},c)
\varrho(b).
\end{align*}
Since $\mu^2\in H_0\cap Z(H)$, we can apply \ref{propf}(c) and get
\begin{align*}
f(a,c^{h(b)})=\varrho(h(b)^{-\mu}) f(a^{h(b\diamond \mu)\mu^2},c)\varrho(b)=\\
\varrho(h(b)^{-\mu})\varrho(\mu^{2})
f(a^{h(b\diamond \mu)},c) \varrho(b)=-\varrho(h(b)^{-\mu})f(a^{h(b\diamond \mu)},c)\varrho(b).
\end{align*} 
Now we can apply (a) and get
\begin{align*}
f(a,c^{h(b)})= -\varrho(h(b)^{-\mu})\varrho(h(b\diamond \mu)^{-\mu})f(a,c)\varrho(b)+\\
\varrho(h(b)^{-\mu})\varrho(h(b\diamond
\mu)^{-\mu})f(a,b\diamond \mu)\varrho(h(b\diamond \mu))^{-1}f(b\diamond \mu,c)\varrho(b).
\end{align*} 
By \ref{mu a tau}(b) we have $\varrho(h(b\diamond \mu))=\varrho(\mu^2 h(b)^{-1})=-\varrho(b)^{-1}$. 
Thus we get with \ref{-mu}(b) and \ref{f a mu}(b)
\begin{align*}
f(a,c^{h(b)})=f(a,c)\varrho(b)+f(a,b\diamond \mu)\varrho(b)f(b\diamond \mu,c)\varrho(b)=\\
f(a,c)\varrho(b) -f(a,b)\varrho(b^{-1})^{-1}f(b,c)\varrho(b)
=\\
f(a,c)\varrho(b)+f(a,b)\varrho(h(b)^{\mu})f(b,c)\varrho(b).\end{align*}
 \end{enumerate}
 \end{proof}
 To ensure that $A_0$ is big enough, we prove the following lemma.
\begin{prop}\label{A_0=2}
If $A_0$ has order $2$, then $G \cong \PSU(3,2)$.
\end{prop}
\begin{proof} The map $\overline{a} \mapsto a^2$ is a quadratic map from $A/A_0$ to $A_0$. If 
$A/A_0$ has order greater 
than $4$, then there is an element $a \in A \setminus A_0$ with $a^2 =1$. Thus $h(a)^{-\mu} = -
h({a^{-1}}) = h(a)$ and so 
$h(a)^2 = h(a) h(a)^{-\mu} \in H_0 =1$. 
Let $z$ be the non-trivial element in $A_0$. Then $\varrho({az}) = 
\varrho(a) +\varrho(z) = \varrho(a)+1$ and thus 
$\varrho({az})^2 = (1+\varrho(a))^2 = 1+\varrho(a)^2 = 1+1=0$, a contradiction. 
Thus $|A:A_0| \leq 4$. If $|A:A_0|
=2$, then $A \cong Z_4$, so 
$G$ is a Frobenius group of order $20$. Since this group has not got a cubic module, this is a 
contradiction. If 
$A/A_0$ has order $4$, then $A \cong Q_8$ and $G$ is isomorphic to $\PSU(3,2)$.
\end{proof}

\begin{coro}\label{G proper} Suppose that $A_0\ne 1$ and $A_0\ne A$. If $G$ is improper, then $G\cong \PSU(3,2)$.
\end{coro}
\begin{proof}
We may assume that $V$ is squarefree and in standard form. 
If $G$ is improper, then $H\leq Z(G)$. Hence \ref{mu quadrat} implies that $\characteristic V=2$, otherwise 
$C_V(A)+C_V(B)$ would be a $G$-invariant subspace of $V$ on which $G$ acts quadratically. By \ref{proper 
Moufang set} and \ref{special elements}(b) $A^*$ contains at most one special element. Since $A_0$ is 
a special root subgroup of $A$, this implies that $|A_0|=2$, thus the claim follows by \ref{A_0=2}.
\end{proof}

 The next result might be useful for an inductive approach. 
\begin{prop}\label{rootgroups} 
\begin{enumerate}
\item $C_A(v)$ is a root subgroup of $A$ for all $v \in Z_2(V,A)\cap Z_2(V,B)$.
\item If $W $ is a subgroup of $Z_2(V,A)\cap Z_2(V,B)$, 
then $C_A(V/(W +C_V(A)))$ is a root subgroup of $A$.
\end{enumerate}
\end{prop}
\begin{proof}
\begin{enumerate}
\item If $a \in C_A(v)$ and $b = \fatb(a)^{-1}$, then $a^{\sim -} = \fata(b)^{-1} $ and
$\mu(a) =\mu a b a^{\sim-}$. We have 
\begin{align*}v\mu(a) = va b a^{\sim-} = vba^{\sim-} = 
(v +[v,b])a^{\sim-} \in (Z_2(V,A)\cap Z_2(V,B))\mu(a)=\\
Z_2(V,A^{\mu(a)}) \cap Z_2(V,B^{\mu(a)}) =Z_2(V,B)\cap Z_2(V,A)\end{align*}
and so $v +[v,b] \in Z_2(V,A) $ and hence $[v,b] \in Z_2(V,A)$. But $v \in Z_2(V,A)\cap Z_2(V,B) 
\leq Z_2(V,B) $ and so 
$[v,b] \in Z_2(V,A) \cap C_V(B) =0$. It follows $b \in C_B(v)$. Similarly, one can show that if 
$b \in C_B(v)$, then $\fata(b) \in C_A(v)$. Thus the claim follows.
\item Let $v \in C_V(A)$, $a \in C_A(V/(C_V(A) +W))$ and $b=\fatb(a)$. Then by \ref{rho}(a) we have
$$vh(a)= [v\mu,a]+v\mu+[v\mu,a,b^{-1}].$$
Now $vh(a) -[v\mu,a] \in Z_2(V,A)\cap Z_2(V,B) \cap (C_V(A) + W) = W$. 
Thus $v\mu -[v\mu,a,b] \in W$. Now 
$v\mu -[v\mu,a,b] = v\mu -[[v\mu,a] -vh(a),b] -[v h(a),b]$ and so $[vh(a),b] \in W +C_V(B)$ because 
$v\mu - [[v\mu,a] -vh(a),b] \in C_V(B)$. This implies $[C_V(A),b] \subseteq W + C_V(B)$. 
Since $V =C_V(A) \oplus [V,B]$ and $[[V,B],b] \subseteq [V,B,B]\subseteq C_V(B)$, it follows 
that $b \in C_B(V/(C_V(B) +W))$. 
A similar argument shows that if $b \in C_B(V/(W +C_V(B)))$, then $\fata(b) \in C_A(V/(W+C_V(A)))$. 
Thus 
the claim follows. 
\end{enumerate}
\end{proof}
\chapter{Construction of irreducible submodules}
In this chapter we will show that if $S$ is a skew field, then one can reduce to the case that $V$ is 
irreducible as a $G$-module.
We will not need any finiteness assumption. 
\newline
 During the whole chapter we will assume that $V$ is a squarefree cubic module for $G$ in standard form. 
Let $W$ be a $H$-submodule of $C_V(A)$. We set $X(W) := W + W\mu + \Phi(W,A)$. 
\begin{lemma}\label{X(W)} 
Suppose that $\mu|_{Z_2(V,A)\cap Z_2(V,B)} =\pm \id$. Then we have:
\begin{enumerate}
\item $X(W)$ is a $G$-submodule of $V$.
\item $X(W)$ is the direct sum of $W, W\mu$ and $\Phi(W,A)$.
\item $X(W) \cap C_V(A) = W, X(W) \cap C_V(B) = W\mu$ and $X(W) \cap Z_2(V,A)\cap Z_2(V,B) = \Phi(W,A)$.
\item If $A_0\ne 1$, then $[X(W),G] = X(W)$. 
\item $W$ is a irreducible $H$-module iff $X(W)$ is a irreducible $G$-module.
\item $V$ is a completely reducible $G$-module iff $C_V(A)$ is a completely reducible $H$-module.
\end{enumerate}
\end{lemma}
\begin{proof}
\begin{enumerate}
\item Since $G=\langle A,\mu \rangle$, we only have to show that $X(W)$ is normalised by $A$ and by $\mu$.
For $w \in W$ we have $w\mu \in W\mu \leq X(W)$ and $w\mu^2 \in W \leq X(W)$, since $W$ is $H$-invariant.
Moreover, by our assumption we have $w\mu =\pm w$ for all $w\in \Phi(W,A)$, so $X(W)$ is $\mu$-invariant. 
Suppose $a,b \in A$ and $w \in W$. Then $wa =w \in W$, 
\begin{align*} w\mu a = w\mu +[w\mu,a] =w\mu +[w\mu,a] -wh(a) +wh(a) =\\
w\mu +\Phi(w,a) +wh(a) \in W\mu + \Phi(W,A) + W = X(W)\end{align*}
and \begin{align*} \Phi(w,b)a = ([w\mu,b] - wh(b))a =[w\mu,b] a -wh(b) 
=\\ [w\mu,b] +[w\mu,b,a] -wh(b) 
 = \Phi(w,b) +wf(b,a) \in \Phi(W,A) + W \subseteq X(W)\end{align*}
 by \ref{coro handf}. 
Thus the claim follows.
\item This is clear since $W\leq C_V(A)$, $W\mu \leq C_V(B)$, 
$\Phi(W,A)\leq Z_2(V,A)\cap Z_2(V,B)$ 
and 
$V = C_V(A) \oplus (Z_2(V,A)\cap Z_2(V,B)) \oplus 
C_V(B)$.
\item This follows easily from (b).
\item We first claim that we have $[X(W),A]=W+\Phi(W,A)$. 
We have $W=[W\mu,A_0] \leq [X(W),A]$ by \ref{rho}(b) and so 
$$\Phi(w,a) = [w\mu,a] -vha \in [X(W),A] +W =
[X(W),A],$$ thus $W+\Phi(W,A) \leq [X(W),A]$. Because we have $[W\mu,A] \leq 
\Phi(W,A) +W$, $[\Phi(W,A),A] \leq W$ and 
$[W,A] =0$, we get equality. 

Now conjugating with $\mu$ shows $[X(W),B]=W\mu+\Phi(W,A)$ and thus the claim follows.
\item Suppose $W$ is a irreducible $H$-module and let $0 \ne Y$ be a $G$-submodule of $X(W)$. 
Since $[Y,A,A,A]=0$, we have $0\ne Y\cap C_V(A)\leq X(W) \cap C_V(A)=W$. Since $W$ is irreducible, we get 
$Y\cap C_V(A)=W$. Thus also $W\mu\subset Y$. 
For all $w\in W$ and $a\in A$ we have 
$\Phi(w,a)=[w\mu,a]-wh(a)\in Y$, which shows $\Phi(W,A)\cap Y$ and $Y=X(W)$. Thus $X(W)$ is irreducible.

Now if $W'$ is a proper submodule of $W$, then we have by (c) 
$X(W')\cap C_V(A) =W' \subsetneq W=X(W)\cap 
C_V(A)$ which shows that $X(W')\subsetneq X(W)$. Hence $X(W)$ is reducible. 
\item Suppose $C_V(A) = \bigoplus_{i \in I} W_i$ is a direct decomposition of $C_V(A)$ as a sum of 
irreducible $H$-modules.
We claim that $V = \bigoplus_{i \in I} X(W_i)$ is a direct decomposition of $V$ as a sum of irreducible 
$G$-modules.

Set $X:= \sum_{i \in I} X(W_i)$. Then $C_V(A), C_V(B)=C_V(A)\mu 
\subseteq X$. Therefore we have $V=X+(Z_2(V,A)\cap Z_2(V,B))$. Therefore we get $[V/X,A]=0=[V/X,B]$ and 
hence $G$ acts trivially on $V/X$. Since $V$ is in standard form and squarefree, we get $V=X$.

We still have to show that the sum is direct. Suppose this is not the case. 
Then there is an element $i \in I$ and a finite subset 
$I_0$ of $I$ with $i \not\in I_0$ and $X(W_i) \cap (\sum_{j \in I_0} X(W_j)) 
\ne 0$. Since $X(W_i) $ is a irreducible $G$-module, 
this implies $X(W_i) \subseteq \sum_{j\in I_0} X(W_j)$. But then we also 
have 
$$W_i =X(W_i) \cap C_V(A) \subseteq C_V(A) \cap (\sum_{j\in I_0} X(W_j)) = 
\sum_{j\in I_0} X(W_j),$$ as one can easily see. But this is a contradiction 
since the decomposition of $C_V(A)$ is direct.
\newline 
Suppose $V =\bigoplus_{i \in I} V_i$ such that $V_i$ is irreducible for all $i 
\in I$. Then $V_i \cap C_V(A) \ne 0$. Therefore $0 \ne X(V_i \cap C_V(A))$ and so 
$X(V_i \cap C_V(A))\cap V_i$. Since $V_i$ is irreducible, 
it follows by (e) that $X(V_i\cap C_V(A)) =V_i$ and that 
$V_i \cap C_V(A)$ is an irreducible $H$-module. 
Set $W =\sum_{i \in I} (V_i \cap C_V(A))$. Then 
$V_i =X(C_V(A) \cap V_i) \subseteq X(W)$ for all $i \in I$ and so $X(W) = 
V$. Thus $W = X(W) \cap C_V(A) = V \cap C_V(A) = C_V(A)$. 
Since $V=\bigoplus_{i \in I} V_i$ is a direct sum, the sum $\sum_{i \in I} (C_V(A)\cap V_i)$ 
is also direct. \end{enumerate} \end{proof}
\begin{remark} By \ref{mu quadrat} the condition $\mu_{Z_2(V,A)\cap Z_2(V,B)}=\pm \id$ is always satisfied 
if $A_0\ne 1$.
\end{remark}
\begin{coro} If $S$ is a skew field or a finite-dimensional central-simple 
algebra, then $V$ is a completely 
reducible $G$-module.
\end{coro}
\chapter{Cubic rank one groups with trivial kernel}
In the following chapter we will treat the case that $A_0=1$. Since $A' \leq A_0$ by 
\ref{quadratic kernel}, we get that $A$ is abelian.
\begin{lemma}\label{Z_2} 
Suppose that $V$ is squarefree. 
For all $a \in A^*$ and $v\in Z_2(V,A) \cap Z_2(V,B)$ we have 
$[v,a] =[v\mu(a),a^{\sim}]$.
\end{lemma}
\begin{proof}
We have 
\begin{align*}
v\mu(a) = v a \fatb(a)^{-1} a^{\sim-}=(v+[v,a]) \fatb(a)^{-1} a^{\sim-} = \\
(v+[v,\fatb(a)^{-1}] +[v,a]+[v,a,\fatb(a)^{-1}]) a^{\sim-} =\\
v-[v,a^{\sim}] -[v,\fatb(a)] -[v,\fatb(a),a^{\sim-}] +[v,a] +[v,a,\fatb(a)^{-1}]+\\
[v,a,\fatb(a)^{-1}, a^{\sim-}].
 \end{align*}
Hence we have
\begin{align*} v\mu(a) -v+[v,a^{\sim}] +[v,\fatb(a),a^{\sim-}]-
[v,a]-[v,a,\fatb(a)^{-1},a^{\sim-}]=\\
-[v,\fatb(a)]+[v,a,\fatb(a)^{-1}]
\in Z_2(V,A) \cap [V,B] \subseteq Z_2(V,A) \cap Z_2(V,B).\end{align*}
Thus we have $[v,\fatb(a),a^{\sim-}] -[v,a,\fatb(a)^{-1},a^{\sim-}] \in 
[Z_2(V,A),A] \subseteq C_V(A)$ and therefore 
\begin{align*} [v,a^{\sim}] -[v,a]+ [v,\fatb(a),a^{\sim-}] -[v,a,\fatb(a)^{-1},a^{\sim-}] =
\\ v-v\mu(a)-[v,\fatb(a)]+[v,a,\fatb(a)^{-1}] 
\in Z_2(V,A) \cap Z_2(V,B) \cap C_V(A) =0.\end{align*}
Hence we have $v\mu(a) =v-[v,\fatb(a)] +[v,a,\fatb(a)^{-1}]$ and  
\begin{align*} [v,a] = [v,a^{\sim}] +[v,\fatb(a),a^{\sim-}] -[v,a,\fatb(a)^{-1},a^{\sim-}] 
 =\\
[-v+[v,\fatb(a)]-[v,a,\fatb(a)^{-1}],a^{\sim-}]= [-v\mu(a),a^{\sim-}] =[v\mu(a),a^{\sim}].
\end{align*} 
\end{proof}
    \begin{theorem}\label{special} Let $G$ be a cubic rank $1$ group with abelian unipotent 
    groups $A$ and 
    $B$. Then $G$ is special.
  \end{theorem}
  \begin{proof}
Since quadratic rank one groups are special by \ref{Jordan}, we may assume that $V$ is squarefree.   
By Segev's Theorem (see \ref{Segev}) $G$ is special or $G$ is 
improper. Suppose that $G$ is not special and hence that $G$ is 
improper. Then the corresponding Moufang set for $G$ is also improper, so 
we have $H\leq Z(G)$ and $Z(G)\mu(a)=Z(G)\mu(b)$ for all $a,b\in A^*$. This also implies that 
$\mu(a)\mu(b)=\mu(b)\mu(a)$ for all $a,b\in A^*$.

Suppose that $\characteristic V =2$. 
Since $[V,A,A,A]=0$, one sees easily that $A$ has exponent at most 
$4$. 
 By \cite[5.2]{Kerby} $A$ contains at most one involution. Since $A$ is 
 abelian of exponent $4$, it 
 follows that $A$ is cyclic of order $4$ and $G \cong \mathrm{Aff}
 (\mathbb{F}_5)$, a contradiction since 
 this group has no cubic module. Hence $\characteristic V \ne 2$. 
    \\
Suppose that $t\in A$ is an involution. Ghen the minimal polynomial of $t$ 
divides both $X^2-1$ and $(X-1)^2$. Since $\characteristic V \ne 2$, this 
implies that $t$ acts trivially on $V$,
a contradiction to \ref{faithful2}. So $A$ contains no involution.
By \ref{special elements} it follows that if $a \in A^{*}$ is special, then $a$ has order $3$ 
and that $\{a,-a\}$ are the only special elements in $A$. 
    \\
Suppose that $A_0 \ne 1$. Since $A_0$ is a special root subgroup of $A$, it follows that $|A_0| 
\leq 3$ and thus $|A_0|=3$. Hence $G_0 \cong \SL(2,3)$ by \ref{mu quadrat}. 
Now if $z$ is the central involution in 
$G_0$, then $z \in H \leq Z(G)$. Therefore $z$ centralises $A$ and $B$. 
But $z$ inverts every element in 
$C_V(A)+C_V(B)$ and fixes every element in $Z_2(V,A)\cap Z_2(V,B)$ by 
\ref{mu quadrat}, so
$C_V(z) =Z_2(V,A)\cap Z_2(V,B)$ is a $G$-invariant submodule of $V$ on 
which $G$ acts quadratically, a contradiction since we assume that 
$V$ is squarefree. Hence $A_0=1$.

    By assumption there is a non-special element $a \in A^*$.  
By \ref{Z_2} we have for $[v,a]=[v\mu(a), a^{\sim}] = -[v\mu(a), a^{\sim-}]$, since 
the commutator map from $Z_2(V,A) \times A$ to $C_V(A)$ is biadditive. 
Since $\mu(a) = \mu(a^{\sim-})$ and 
   $a^{\sim -\sim} = a^{-\sim-}$ by \ref{formula sim}(f), we have \begin{align*}
   [v,a]=-[v\mu(a),a^{\sim-}]=[v\mu(a)^2, a^{\sim -\sim-} ]=
   [v\mu(a)^2, a^{-\sim--}]=\\
   [v\mu(a)^2,a^{-\sim}]=[v\mu(a)^3,a^{-\sim\sim}]=
   -[v\mu(a)^3, a]\end{align*}
   and therefore $[v\mu(a)^6,a]=[v,a]$.
   Since $\mu(a)^6 \in Z(G)$, we replace $V$ by $C_V(\mu(a)^6)$ and therefore assume that 
   $\mu(a)^6=1$. 
   Since $\mu(a)^3$ interchanges $A$ and $B$, $\mu(a)^3$ has order $2$. If $\mu(a)^2=1$, then we have
   $$[v,a]=[v\mu(a)^2 ,a^{-\sim}]=[v,a^{-\sim}].$$
   It follows that $a^{-1}a^{-\sim} \in A_0=1$ and so 
   $a =a^{-\sim}$, thus $a^{\sim} =a^{-1}$. Hence $a$ is special, a 
   contradiction. 
   Therefore $o(\mu(a))=6$. Moreover, since $\mu(a)^2 \in Z(G)$, we may assume that $C_V(\mu(a)^2)=0$. Thus 
   $\characteristic V\ne 3$ and $\mu(a)^2$ has 
    minimal polynomial $X^2 +X +1$.
Furthermore, if $a,b \in A^*$, then $\mu(a) \mu(b) \in Z(G)$ and has order at most $6$.

    Let $F$ be $\mathbb{F}_p$ if $\characteristic V =p$ and set $F=\mathbb{Q}$ if $\characteristic V =0$. 
   Let $a \in A^*$, set $b =a^2$ and $c=aa^{-\sim}$.
   Note that if $c =1$, then $a=a^{-\sim-}=a^{\sim-\sim}$ and so 
$a^{\sim} =a^{\sim-\sim\sim} =a^{\sim-}$,
so $a^{\sim}$ has order $2$, a contradiction. Therefore $b,c \ne 1$. 
    Let $E$ be the subring in $\End_G(V)$ generated by
   $Z:=Z(G) \cap \langle \{\mu(a),\mu(b), \mu(c)\} \rangle$ and $F$.
   Then $E$ is finite-dimensional over $F$, and since $E$ is generated by elements of order at 
   most $6$. Furthermore, 
   $\characteristic F \ne 2,3$ implies that $E$ is semisimple. If we replace $V$ by $Ve$ for a 
   primitive idempotent
   $e \in E$, we may assume that $E$ is a field. Thus $Z$ is cyclic of order at most $6$. 
    Suppose that $\mu(a) \mu(b)^{-1}$ has order $6$. Then we may assume that 
    $(\mu(a) \mu(b)^{-1})^3$ inverts every element in $V$. We have
    $(\mu(a) \mu(b)^{-1})^3 =    \mu(a)^3 \mu(b)^{-3}$
    and so \begin{align*}
    [v,b]=-[v(\mu(a) \mu(b)^{-1})^3,b]=-[v\mu(a)^3 
    \mu(b)^{-3},b]=[v\mu(a)^3,b]=[v\mu(a)^3,a^2]= \\
    2[v\mu(a)^3,a]=-2[v,a]=-[v,a^2]=-[v,b]
    \end{align*} for all $v \in [V,A]\cap [V,B]$.
   Since $A_0 =1$, this implies $\characteristic V =2$, a contradiction. Hence $\mu(a) 
   \mu(b)^{-1} $ has order at 
   most $3$. Since $\mu(a) =\mu(a^{\sim -})$, the same argument shows that $\mu(a) 
   \mu(c)^{-1}$ has order
   at most $3$. Thus if $x \in \{b,c\}$ we have
   $\mu(x) \in \{1,\mu(a)^2,\mu(a)^4 \} \mu(a) =\{\mu(a), \mu(a)^3,\mu(a)^{-1}\}$.
Suppose $\mu(b)=\mu(a)^3$. Then we have
$$[v,b]=-[v\mu(a)^3,b]=-[v\mu(b),b]=-[v,b^{\sim}]=[v,b^{\sim-}],$$
thus $b=b^{\sim-}$ and so $b^\sim=b^{-1}$, which implies that $b$ is 
special, a contradiction. The same argument applies for $c$.  Hence we get $\mu(b),\mu(c) =\mu(a)^{\pm 1}$.
   \\
    For $v \in C_V(A)$ the element $[v\mu,a]-v\varrho(a)$ is contained in 
    $Z_2(V,A)\cap Z_2(V,B)$. Therefore we obtain by applying \ref{Z_2} twice and \ref{formula sim}(f)
    \begin{align*}
    vf(a,a)=[v\mu,a,a]=[([v\mu,a] -v\varrho(a)),a]+[v\varrho(a),a]=\\
    [([v\mu(a),a]-v\varrho(a))\mu(a),a^
    \sim]=
    -[([v\mu(a),a]-v\varrho(a))\mu(a),a^{\sim -}]=\\
    -[([v\mu(a),a]-
    v\varrho(a))\mu(a)\mu(a^{\sim-}),a^{\sim-\sim}]=
-[([v\mu(a),a]-v\varrho(a))\mu(a)^2,a^{-\sim-}]=    
\\    
    [([v\mu(a),a] -v\varrho(a))\mu(a)^2,a^{-\sim})]=
  [[v\mu(a),a]-v\varrho(a), (a^{-\sim})^{\mu(a)^{-2}}]\mu(a)^2=\\
  [[v\mu(a),a]-v\varrho(a), a^{-\sim}]\mu(a)^2 
  =vf(a,a^{-\sim})\mu(a)^2.
  \end{align*}
    Thus we have by \ref{handf}(a)
    \begin{align*}
    v(\varrho(b) -2\varrho(a)) =v(\varrho(a^2)-2\varrho(a))=
    vf(a,a)=vf(a,a^{-\sim}))\mu(a)^2=\\
    v(\varrho(c) -\varrho(a)
    -\varrho(a^{-\sim}))\mu(a)^2=v(\varrho(c)-2\varrho(a))
    \mu(a)^2\end{align*}
    for all $v \in C_V(A)$. 
    If $\mu(c) =\mu(b) =\mu(a)$, we get
    $-v\varrho(a) =-v\varrho(a) \mu(a)^2$ for all $v \in C_V(A)$, and so $C_V(A) \leq C_V(\mu(a)^2)=0$, a contradiction. If $\mu(b) =\mu(c) =\mu(a)^{-1}$, we have $\varrho(b)=\varrho(c) =\varrho(a)\mu(a)^{-2}$ 
    and hence we get
      $$\varrho(a)(\mu(a)^{-2}-2)=\varrho(a)(\mu(a)^{-2}-2)\mu(a)^2=
      \varrho(a)(1-2\mu(a)^2).$$
This implies that $2\mu(a)^4-3\mu(a)^2-1=0$. But $\mu(a)^2$ has minimal polynomial $X^2 +X +1$, so this is a contradiction. Now suppose $\mu(b) =\mu(a) $ and $\mu(c) =\mu(a)^{-1}$. Then we have
      $$-\varrho(a) =\varrho(a)(\mu(a)^{-2} -2)\mu(a)^2 =\varrho(a)
      (1-2\mu(a)^2)$$ and thus
      $2\mu(a)^2=2$. This is again a contradiction. Finally, assume $\mu(c)=\mu(a)$ and 
      $\mu(b) =\mu(a)^{-1}$. 
      Then we have
      $$\varrho(a)(\mu(a)^4 -2)=-\varrho(a) \mu(a)^2$$ and so
      $\mu(a)^4 +\mu(a)^2 -2=0$, again a contradiction. 
      Thus our assumption that $G$ is not special yields a contradiction in 
      any case and the claim follows.
   \end{proof}
   \begin{lemma}\label{mu a 2} Suppose that $A_0=1$. Then there is a submodule $W$ of $V$ such 
   that $G$ 
   acts quadratically on $V/W$ and that $\mu(a)^2|_W=\id$ for all $a\in A^*$ and 
   \end{lemma}
   \begin{proof}
   If $A$ has exponent $2$, then $\mu(a) =\mu(a^{-1}) =\mu(a)^{-1}$ and so $\mu(a)^2=1$ 
   for all $a\in A^*$. 
   Suppose that $A$ has not exponent $2$.
  Since $G$ is special and $A$ is abelian, $\mu(a)^2 \in Z(G)$ for all $a\in A^*$ by \cite[7.5.2]{DS}. 
  By \ref{Z_2} and \ref{special} we have $[v,a] = [v\mu(a),a^{\sim}]= [v\mu(a),a^{-1}] =-[v\mu(a),a] $ for 
  all $v\in Z_2(V,A) \cap Z_2(V,B)$. By replacing $v$ with $v\mu(a)$ we get
  $[v\mu(a)^2,a] = -[v\mu(a),a] = [v,a]$. Since $\mu(a)^2 \in Z(G)$, we conclude 
  $$[v,a] =[v\mu(a)^2,a] = v\mu(a)^2 (a-1) =v(a-1)\mu(a)^2 = [v,a]\mu(a)^2.$$
  Thus we have $[Z_2(V,A) \cap Z_2(V,B),a] \leq C_V(\mu(a)^2)$. Since $\mu(a)^2 \in Z(G)$, 
  $C_V(\mu(a)^2)$ is a $G$-submodule of $V$. 
  Let $\widehat{V}=V/C_V(\mu(a)^2)$ and 
  $\widehat{A}_0 =C_A([\widehat{V},A]) \cap C_A(\widehat{V}/C_{\widehat{V}}(A))$. 
  Since $[Z_2(V,A)\cap Z_2(V,B) \leq C_V(\mu(a)^2)$, we get $a\in \widehat{A}_0$. But 
  $\widehat{A}_0$ is a $H$-invariant root subgroup of $A$, and since $A$ has not exponent $2$, 
  we get $\widehat{A}_0=A$ by \ref{root subgroups}. Thus we have 
  $[V,A,A] \leq C_V(\mu(a)^2)$ and, since $a$ was arbitrarily chosen, $[V,A,A] \leq \bigcap_{a\in 
  A^*} C_V(\mu(a)^2) =:W$. Thus the claim follows. \end{proof}
  
  \bigskip
  By replacing $V$ with $W$, if necessary, we will assume that $\mu(a)^2=1$ for all $a\in A^*$.
  \begin{prop}\label{fha} Suppose that $A_0=1$ and $\mu(a)^2=1$ for all $a\in A^*$. 
  Then we have
\begin{enumerate}  
\item 
      $f(c,a^{h(b)}) = 
      f(c,a) \varrho(b) -f(c,e)f(e,a)\varrho(b)- f(c,b)f(b,a)+f(c,b)f(a,e)$ $f(b,e)$,
\item $f(e,e)f(e,a)=2f(e,a) =f(e,a)f(e,e)$.      
      \end{enumerate}
      for all $a,b\in A^*$. 
  \end{prop}
  \begin{proof} Note that $f$ is symmetric by \ref{propf}(e). Since all $\mu$-maps are assumed to be 
  involutions, we have $h(x)^{\mu}=\mu^2 \mu(x)\mu =h(x)^{-1}$,   
  $h(x^{-1})=h(x)$ and $h(x\diamond \mu) =h(x)^{-1}$ for all $x\in A^*$ by \ref{-mu}.
\begin{enumerate}
\item
 By applying \ref{formulas}(a) and \ref{f b a mu} twice we have \begin{align*}
  f(c, (ab^{-1})\diamond \mu(b)) = f(c,(((a\diamond \mu)(b\diamond \mu)^{-1})\diamond \mu)b)
=\\
f(c\diamond \mu,((a\diamond \mu)(b\diamond \mu)^{-1}))\diamond \mu)+f(c,b)=\\
\varrho(c\diamond \mu)^{-1}f(c\diamond \mu,(a\diamond \mu)(b\diamond \mu)^{-1})\varrho((a\diamond \mu)
(b\diamond \mu)^{-1})^{-1}+f(c,b)=\\
 \varrho(c)f(c\diamond \mu,a\diamond \mu)\varrho((a\diamond \mu)(b\diamond \mu)^{-1})^{-1}
 -\\ \varrho(c)f(c\diamond \mu,b\diamond \mu)\varrho((a\diamond \mu)(b\diamond \mu)^{-1})^{-1}+ f(c,b)=\\
f(c,a)\varrho(a)^{-1} \varrho((a\diamond \mu)(b\diamond \mu)^{-1})^{-1}-\\
f(c,b)\varrho(b)^{-1}\varrho((a\diamond \mu)(b\diamond \mu)^{-1})^{-1}+f(c,b)
\end{align*}  
By \ref{h+} we have $$h((a\diamond \mu)(b\diamond \mu)^{-1}) =h(a)^{-1}h(ba^{-1})h(b)^\mu=
h(a)^{-1}h(ab^{-1})h(b)^{-1}.$$ Moreover, we have
\begin{align*} h((a\diamond \mu)(b\diamond \mu)^{-1}) =h((b\diamond \mu)(a\diamond \mu)^{-1})=\\ h(b)^{-1}h(ba^{-1})h(a)^{-1}
 =h(b)^{-1}h(ab^{-1})h(a)^{-1}.
\end{align*}
Plugging these formulas into the equation above we get
  \begin{align*}f(c,(ab^{-1})\diamond \mu(b)) = \\
  f(c,a) \varrho(ab^{-1})^{-1} \varrho(b) -f(c,b) \varrho(ab^{-1})^{-1} 
  \varrho(a) +f(c,b).\end{align*}
  Replacing $a$ by $ab$ and applying \ref{handf} and \ref{f b a mu} again we get 
  \begin{align*} 
    f(c,a\diamond \mu(b)) =f(c,ab) \varrho(a)^{-1} \varrho(b) -f(c,b) \varrho(a)^{-1} \varrho(ab) +f(c,b) = \\
        f(c,a)\varrho(a)^{-1} \varrho(b) +f(c,b) \varrho(a)^{-1} \varrho(b)-f(c,b) \varrho(a)^{-1} 
        \varrho(ab) +f(c,b) =\\
    \varrho(c) f(c\diamond \mu,a\diamond \mu) \varrho(b)- f(c,b) \varrho(a)^{-1} (\varrho(ab) -\varrho(a) 
    -\varrho(b)) =\\ 
  \varrho(c) f(c\diamond \mu,a\diamond \mu) \varrho(b) -f(c,b) \varrho(a)^{-1} f(a,b)=\\ \varrho(c) f(c\diamond \mu, 
  a\diamond \mu) \varrho(b) -
  f(c,b)f(a\diamond \mu,b\diamond \mu)\varrho(b).
    \end{align*}
  Replacing $a$ by $a\diamond \mu$ and applying \ref{f b a mu} once again we get
  \begin{align*}  f(c,a^{h(b)})  \stackrel{(1)}{=}    \varrho(c) f(c\diamond \mu, a) \varrho(b) -f(c,b) f(a,b\diamond \mu) \varrho(b)    
  \stackrel{(2)}{=}\\ \varrho(c)f(c\diamond \mu,a) \varrho(b)
    -f(c,b) \varrho(a) f(a\diamond \mu ,b) .\end{align*}
  Using (1) for $b=e$ and the fact $e\diamond \mu=e\diamond \mu(e)=-e$ by \ref{special elements}, 
  we get 
  \begin{align*} 
  f(c,a) = f(c,a^{h(e)}) = \varrho(c) f(c\diamond \mu,a)-f(c,e)f(a,e\diamond \mu) =\\
   \varrho(c) f(c\diamond \mu,a)+f(c,e)f(a,e),
  \end{align*}
thus 
$$\varrho(c) f(c\diamond \mu,a) = f(c,a) -f(c,e)f(a,e).$$
Hence we get with equation (2)
\begin{align*}
f(c,a^{h(b)}) = \\ f(c,a) \varrho(b) -f(c,e)f(e,a)\varrho(b)-f(c,b)f(a,b)+f(c,b)f(a,e) f(b,e)=\\
f(c,a) \varrho(b) -f(c,e)f(e,a)\varrho(b)-f(c,b)f(b,a)+f(c,b)f(a,e) f(b,e).\end{align*} 
\item Using equation (1) for $b=c=e$ we get\begin{align*}
f(e,a)=f(e,a^{h(e)}) =\varrho(e) f(e\diamond \mu,a)\varrho(e)-f(e,e)f(a,e\diamond \mu)\varrho(e)=\\
-f(e,a)+2f(e,e)f(e,a)\end{align*}
and hence $f(e,e)f(e,a)=2f(e,a)$.
By \ref{f b a mu} we have $\varrho(c)f(c\diamond \mu,e)=-\varrho(c)f(c\diamond \mu,e\diamond \mu)\varrho(e)=
-f(c,e)$. So equation (2) with $a=b=e$ implies\begin{align*}
f(c,e)=f(c,e^{h(e)})=\varrho(c)f(c\diamond \mu,e) -f(c,e)\varrho(e)f(e\diamond \mu,e) =\\
-f(c,e)+f(c,e)f(e,e)
\end{align*}
and hence $2f(c,e)=f(c,e)f(e,e)$.
\end{enumerate}
   \end{proof}
   
\bigskip
We are now able to prove Main Theorem 1.   
    \begin{theorem}\label{Jordan2} Let $G$ be a rank one group having a cubic module $V$. 
  If the quadratic kernel is trivial, 
  then $G/Z(G)\cong \PSL_2(J)$ for a quadratic Jordan division algebra $J$.
\end{theorem} 
  \begin{proof}
  By \ref{normal form} 
  we may suppose that $V$ is in standard form. If $V$ is not squarefree, then by \ref{Jordan} 
  $G\cong \PSL(J,R)$ for 
  a special quadratic Jordan algebra inside a ring $R$, so the claim is true in this case. 
Therefore we may suppose that $V$ is squarefree. Hence by \ref{mu a 2} we have $\mu(a)^2=1$ for all 
$a\in A^*$.

  For $a,b,c \in A$ set $a^{h(b,c)} :=a^{h(bc)} (a^{h(b)} a^{h(c)})^{-1}$ 
  with the convention 
  $a^{h(1)} =1$ for all $a\in A$. Note that $A$ is abelian and so $h(b,c)$ is an endomorphism of 
  $A$. 
 With \ref{fha} we get \begin{align*}
 f(c,e^{h(a)}) =\\
  f(c,e) \varrho(a) -f(c,e) f(e,e) \varrho(a) -
  f(c,a)f(c,e)+f(c,a)f(e,e)f(a,e) =\\
  f(c,e)\varrho(a) -2f(c,e)\varrho(a)-
f(c,a)f(a,e) +2f(c,a)f(c,e)=\\
f(c,a)f(a,e)-f(c,e)\varrho(a).
 \end{align*}
Hence we get \begin{align*}
f(c,e^{h(a,b)}) =f(c,e^{h(ab)}-e^{h(a)}-e^{h(b)}) =\\
f(c,ab)f(ab,e)-f(c,a)f(a,e)-f(c,b)f(b,e) -f(c,e)
(\varrho(ab)-\varrho(a)-\varrho(b)) =\\
f(c,a)f(b,e) +f(c,b)f(a,e) -f(c,e)f(a,b).\end{align*}
Moreover, we have \begin{align*}
f(c,a^{h(b,e)}) = f(c,a^{h(be)}(a^{h(b)}a)^{-1}) = \\
f(c,a) (\varrho({be}) -\varrho(b)-\varrho(e))-f(c,e) f(e,a)
(\varrho(be)-\varrho(b)-\varrho(e)) 
-\\ f(c,be)f(be,a)+f(c,b)f(b,a)+f(c,e)f(e,a)\\
+f(c,be)f(a,e)f(be,e) -f(c,b)f(a,e)f(b,e)-f(c,e)f(a,e)f(e,e)=\\
f(c,a)f(b,e)-f(c,e)f(e,a)f(e,b)-f(c,e)f(b,a)-f(c,b)f(e,a)+\\
f(c,b)
f(a,e)f(e,e)+f(c,e)f(a,e)f(b,e)=\\
f(c,a)f(b,e)-f(c,e)f(b,a)+f(c,b)f(e,a)
=f(c,e^{h(a,b)}).
\end{align*}
Thus we have
$$0=f(c,e^{h(a,b)}-a^{h(b,e)}) =f(e^{h({a,b})}-a^{h({b,e})},c)$$
for all $c\in A$. This implies $e^{h({a,b})}-a^{h({b,e})} \in A_0 =1$ and 
hence $e^{h({a,b})} =a^{h({b,e})}$. Thus (QJ2) holds for the Moufang set corresponding to $G$ by 
\ref{weak QJ2}. 
Therefore $G$ is the rank one group for a quadratic Jordan division algebra 
$J$ by \ref{QJ2}. Hence $G/Z(G)$ is the little projective group of $\mathbb{M}(J)$ which is just 
$\PSL_2(J)$.
 \end{proof}

 Note that by Example \ref{TKK} for every quadratic Jordan division algebra 
 there is a cubic module for 
 $\PSL_2(J)$.
\chapter{A characterisation of the adjoint module}
A. Deloro showed that if $K$ is a field of characteristic not $2$, 
$G=\SL_2(K)$, $A$ the group 
of lower triangular matrices in $G$ and $V$ is an irreducible 
cubic module for $G$ such that $C_V(a)=C_V(b)$ for all $a,b\in A^*$ holds, 
then $V$ is isomorphic to the adjoint module (see the final corollary 
of \cite{Del}). With our methods we can 
generalise this result.
 \begin{theorem}
 Suppose that $V$ is a squarefree cubic module in standard form for $G$ 
 such that $C_V(a) =C_V(A)$ for all $a\in A^*$. 
 Then $G/C_G(V) \cong \PSL_2(K)$ for a 
 commutative field $K$ of odd characteristic and $V$ is a direct sum of copies of the adjoint 
 module for $G$.
 \end{theorem} 
   \begin{proof}
We may assume that $G$ acts faithfully on $V$.  
The condition implies that $A_0=1$, so by Main Theorem 1 
$G$ is the rank one group corresponding to quadratic 
Jordan division algebra.     
     Suppose first that $\characteristic V =2$. Then $A$ is an elementary-abelian $2$-group, hence 
     for all $a\in A^*$ we have 
      $f(a,a)=\varrho(a^2)-\varrho(a) -\varrho(a) =\varrho(1)=0$, so $[v\mu,a,a]=vf(a,a)=0$ for all $v \in C_V(A)$. 
      Hence
      $[v\mu,a]\in C_V(a)=C_V(A)$. But this means that $G$ operates quadratically on $V$, a 
      contradiction. Thus $\characteristic V \ne 2$. \\
   For all $a \in A^{*}$ and all $v \in Z_2(V,A) \cap Z_2(V,B)$ we have by \ref{Z_2}
   $[v,a] =-[v\mu(a),a]$ and thus $[v(1+\mu(a)),a] =0$. Hence $v(1+\mu(a)) \in C_V(a) \cap 
   Z_2(V,A) \cap Z_2(V,B) = C_V(A) \cap Z_2(V,A) \cap Z_2(V,B) = \{0\}$. 
   Thus $v\mu(a) = -v$, It follows  
   that $v\mu(a)\mu(b)=v$ for all $a,b \in A^*$, thus 
   $vh =v$ for all $h \in H$. Hence for all $a,b \in A$, $v\in C_V(A)$ and all $h \in H$ we get
   \begin{align*} vf(a,b^h) =[v\mu,a,b^h] = [[v\mu,a]-v\varrho(a),b^h] +[v\varrho(a),b^h] =\\ [([v\mu,a]-
   v\varrho(a))h,b^h]= 
     [[v\mu,a]-v\varrho(a),b]h =vf(a,b)\varrho(h),\end{align*}
since $[v\mu,a]-v\varrho(a)\in Z_2(V,A)\cap Z_2(V,B)$ by \ref{rho}(a).    
     Thus 
\begin{equation}\label{eq5}
f(a,b^h) =f(a,b) \varrho(h)
\end{equation}     
      and hence $$f(a^h,b) =\varrho(h^{-\mu}) f(a,b^{h^{-1}}) 
     \varrho(h)=\varrho(h^{-\mu})f(a,b)$$ by \ref{propf}(b).
     On the other hand, we have
     $$ f(a,b)\varrho(h) = f(a,b^h) = f(b^h,a) =\varrho(h^{-\mu}) f(b,a)
      =\varrho(h^{-\mu}) f(a,b).$$
     Now if $h=h(c)$ for some $c \in A^{*}$, then $h^{-\mu} =h$ and so 
     $\varrho(c) \in 
     C_S(f(a,b))$. 
    Since the elements $h(c)$ generate $H$, we get $\varrho(H) \subseteq 
    C_S(f(a,b))$ for all $a,b \in A$.
    Thus we have
    $$ f(a,b) \varrho(h^{-\mu}) = \varrho(h^{-\mu}) f(a,b) =f(a,b) 
    \varrho(h)$$ and so
    $$0=f(a,b)(\varrho( h^{-\mu})-\varrho(h))=f(a,b^{h^{-\mu}} -b^h) $$
    for all $a,b \in A$ and all $h \in H$.
    This implies that $b^{h^{-\mu}} -b^h \in A_0 =1$ for all $b \in A$ and all 
    $h \in H$. 
    Therefore
    $hh^{\mu} \in H \cap C_G(A) =Z(G)$. Now since $Z_2(V,A) \cap Z_2(V,B) 
    \subseteq     C_V(H)\subseteq C_V(hh^{\mu})$, 
    it follows that $G$ acts quadratically on $V/C_V(hh^\mu)$ or 
        $C_V(hh^{\mu}) =V$. Since we assume that $V$ is squarefree, 
        the latter must hold, and as a consequence we obtain $hh^{\mu}=1$. 
        Hence $h^{\mu} =h^{-1}$ for all $h \in H$, which 
    implies that
    $H$ is abelian. Thus by \cite{DW} there is a commutative field $K$ such that
    $G$ is a central extension of $\PSL_2(K)$. This field can be constructed as follows: 
Let $(K,+)$ be a group isomorphic to $A$, $\alpha: K \to A:a \mapsto \alpha(a)$ 
an isomorphism, $1\in K$ 
the preimage of $e$ 
and $\tau $ a permutation of $X:=K \dot{\cup}\{\infty\}$ interchanging $0$ and $\infty$ 
such that $\alpha({a\tau})= \alpha({a})\diamond \mu$ for all $a\in K^{\#}$. Then 
$\mathbb{M}(K,\tau)$ is a Moufang set isomorphic to the Moufang set corresponding to $G$ 
and $\tau =\mu_e$. Set $h_a:=\tau \mu_a$ for $a\in K^{\#}$, $h_0 =0$ and 
$h_{a,b} =h_{a+b}-h_a-h_b$.  
The multiplication on $K$ is given by 
    $a \cdot b = ah_{1,b}$ for $a,b \in A$. Note that $a^2 =1^{h_a}$. 
    Set 
    $\varphi: K \to S: a \mapsto \frac{1}{2}f(e,\alpha({a}))$. Then we have using equation (\ref{eq5})
    \begin{align*}
        \varphi(a) \varphi(b) =
        \frac{1}{4} f(e,\alpha({a}))f(e,\alpha({b}))=
        \frac{1}{4}f(e,\alpha({a}))(\varrho(h(e\alpha({b}))-\varrho(e)-\varrho(\alpha({b})))=\\
   \frac{1}{4}f(e,\alpha({a})^{h(e\alpha({b}))} \alpha({a})^{-h(\alpha({b}))}
   \alpha({a})^{-h(e)})=\frac{1}{2}f(e, \alpha({\frac{1}{2}ah_{e,b}}) ) =\varphi(a\cdot b).
 \end{align*}
Since $K$ is a field, $\varphi$ is a monomorphism. Moreover, for $v\in C_V(A)$ we have by \ref{fha}(b)
$(2-f(e,e))f(e,a)=0$ for all $a\in A$ and thus 
$0=v(2-f(e,e))f(a,e)=[v(2-f(e,e))\mu,a,e]=[v(2-f(e,e))\mu-v(2-f(e,e))\varrho(a),e]$. 
Therefore we get $[v(2-f(e,e))\mu,a]-v(2-f(e,e))\varrho(a) \in C_V(e)=C_V(A)$. But by \ref{rho}(a) 
$[v(2-f(e,e))\mu,a]-v(2-f(e,e))\varrho(a)$ is also in $Z_2(V,A)\cap Z_2(V,B)$. Thus 
$[v(2-f(e,e))\mu,a]-v(2-f(e,e))\varrho(a)=0$. Therefore 
$[v(2-f(e,e))\mu,a]=v(2-f(e,e))\varrho(a)\in C_V(A)$ and $v(2-f(e,e))\mu\in C_V(B)\cap Z_2(V,A)=0$. 
This implies $v(2-f(e,e))=0$. We conclude $f(e,e)=2$ and $\varphi(1)=\frac{1}{2}f(e,e)=1$. 
Now with (\ref{eq5}) we get $ \varrho(\alpha({a})) =\frac{1}{2}f(e,e^{h(\alpha(a)}) =\varphi(\alpha(a)^2)$.
Since $S$ is generated by the elements $\varrho(a)$, we conclude that $\varphi$ is bijective. 
Thus $S \cong K$ and we may consider $C_V(A)$ as a $K$-vector space. 
    
Since $\mu|_{Z_2(V,A)\cap Z_2(V,B)} =-\id$, 
    we are allowed to apply \ref{X(W)}. One can easily see that if $L$ is a one-dimensional $K$-subspace 
    of $C_V(A)$, then $X(L)$ is isomorphic to the adjoint module for $G$. 
           Let $(b_i)_{i \in I}$ be a $K$-basis of 
    $C_V(A)$. Then by \ref{X(W)} we have $V =\bigoplus_{i \in I} X(\langle b_i \rangle)$ and the claim follows.
      \end{proof} 
 \chapter{Cubic rank one groups with non-trivial kernel}
 In this chapter we assume that $V$ is squarefree and in standard form and that $A_0\ne 1$. 
We set $\oA =A/A_0$ and $\oa=A_0a$ for $a\in A$. Recall that 
 $J:=\{\varrho(a)|\in A_0\}$ and that $R$ resp. $S$ is the subring of $\End(C_V(A))$ generated by 
 $H_0$ resp. $H$. By \ref{Jordan} $J$ is a special quadratic Jordan algebra inside $R$.

\begin{lemma}\label{normalise} For all $b \in A$, both $\varrho(b)$ and $\varrho(b) +1$ normalise 
$R$.
\end{lemma}
\begin{proof} For $b \in A, a \in A_0$ we have by \ref{-mu}(c) \begin{align*} 
\varrho(b)^{-1} \varrho(a) 
\varrho(b) =\varrho(h(b)^{-1})\varrho(e) \varrho(h(b)^{\mu}) \varrho(h(b)^{-\mu}) \varrho(a) \varrho(b) =\\
\varrho(e^{h(b)^{\mu}}) \varrho(a^{h(b)})) \in R\end{align*} and 
\begin{align*}
\varrho(b) \varrho(a) \varrho(b)^{-1} =\varrho(b) \varrho(h(b)^{-\mu}) \varrho(h(b)^{\mu}) \varrho(h(a)) 
\varrho(h(b)^{-1}) =\\ \varrho(e^{h(b)^{-\mu}}) \varrho(a^{h(b)^{-1}}) \in R.\end{align*}
Thus $\varrho(b) $ normalises $R$. This also holds for $\varrho({be}) =\varrho(b) 
+\varrho(e)=\varrho(b)+1$. 
\end{proof}  
\begin{prop}\label{R=S} Suppose that
\begin{enumerate}
\item $R$ has only finitely many maximal ideals.
\item If $x \in R$ is invertible in $S$, then $x^{-1} \in R$.
\end{enumerate} 
Then $J$ is a commutative quadratic Jordan division algebra or $R=S$.
\end{prop}
\begin{proof} Suppose that $J$ is not commutative. 
Set $Z(J):= \{x \in J| x+\B(R) \in \Cent(R/\B(R))\}$. 
If $M$ is a maximal ideal of $R$ and $x \in J$ with $x + M \in \Cent(R/M)$, then 
$[x,y]^2, [x,y,z] \in M \cap J = 0$ for all $y,z \in J$ by \ref{commutator}. 
Thus $[x,y]$ is nilpotent, and since 
$J \subseteq \Cent_R([x,y])$ and $R$ is generated by $J$, 
we get $[x,y] \in \Cent R$ for all $y \in J$. Therefore $[x,y]r[x,y]=[x,y]^2 r =0$ for all $r\in R$, 
so $[x,y] \in \B(R)$ and $x+\B(R) \in 
\Cent(R/\B(R))$. 
Therefore the preimage of $\Cent(R/M)$ in $J$ is just $Z(J)$ and does not depend on $M$. 
Note that a finite quadratic Jordan division algebra is a field and therefore 
$J/Z(J)$ is either a vector space over $\QQ$ or 
a infinite-dimensional $\FF_p$-vector space. Therefore $(J/Z(J),+)$ contains 
infinitely many cyclic subgroups. 
 
Suppose that $b \in A$ with $\varrho(b)
\not\in R$. Then $\varrho(b)^{-1}\not\in R$ by assumption (b) and $ \varrho({ba}) \not\in R$ for all $a \in 
A_0$. 
As in \ref{subrings}, we set $I_x:=xR \cap R$ for $x\in S$. By \ref{normalise},  
$I_{\varrho(ba)}$ is an ideal of $R$ for all $a\in A_0$.
Now $J/Z(J)$ contains infinitely 
many cyclic subgroups but $R$ has only finitely many maximal ideals by assumption (a). 
Therefore there 
are $a,c \in A_0$ and a maximal ideal $M$ of $R$ such that 
$I_{\varrho({ba})}, I_{\varrho({bc})} \subseteq M$ and such that $\varrho({a}), \varrho({c})$ are 
not 
in $Z(J)$ and the groups generated by $\varrho(a)$ and $\varrho(c)$ in $J/Z(J)$ are distinct. 
Thus also $\varrho({ac^{-1}}) \not\in Z(J)$. 
By replacing $b$ with $bc$ and $a$ with $c^{-1} a$, we may assume 
$I_{\varrho(b)}, I_{\varrho({ba})} \subseteq M$ and $\varrho(a) \not \in Z(J)$.
Therefore 
$$u -(\varrho({b})-1)^{-1} u (\varrho({b})-1)\in M \text{ and }
 u -(\varrho({ba})-1)^{-1} u (\varrho({ba})-1) \in 
M$$ for all
 $u \in R$ by \ref{subrings} and by \ref{normalise} with $be^{-1}$ and $bae^{-1}$ instead of 
 $b$. We get
\begin{align*}\varrho(b) u - u \varrho(b) =-u + \varrho(b) u -u\varrho(b) +u =\\ 
(\varrho({b})-1)u - u(\varrho({b})-1) \in (\varrho({b})-1) M \subseteq \varrho(b) M + M\end{align*}
and similarly 
\begin{align*}
u\varrho(b) +u\varrho(a) -\varrho(b) u -\varrho(a) u = (\varrho({ba})-1)u - u(\varrho({ba}) -1)\\
\in (\varrho({ba})-1) M \subseteq \varrho(b) M + (\varrho({a})-1) M =\varrho(b) M + M.\end{align*}
Thus we get \begin{equation}\label{eq?}
u\varrho(a) -\varrho(a) u \in (\varrho(b) M +M) \cap R =:J_b.\end{equation}
We see immediately that $J_b$ is a right ideal of $R$. Moreover, for all $u,v \in M$ and all $s 
\in R$ we have
$$s (\varrho(b) u +v) = \varrho(b) \varrho(b)^{-1} s \varrho(b) u +sv \in \varrho(b) M +M.$$
This shows that $J_b$ is also a left ideal of $R$. If $J_b = R$, then there are $u,v \in M$ 
with 
$1 = \varrho(b) u +v $. This implies $1-v = \varrho(b) u \in R \cap \varrho(b) R = I_{\varrho(b)} 
\subseteq M$ and thus 
$1 =1-v +v \in M$, a contradiction. Thus $(\varrho(b) M + M) \cap R \ne R$. Since $M \subseteq 
(\varrho(b) M + M) \cap R $, 
we get $M = (\varrho(b) M + M) \cap R $. But now we have $ u\varrho(a) - \varrho(a) u \in M$ for 
all $u \in R$ by (\ref{eq?}) and thus 
$\varrho(a) + M \in \Cent(R/M)$. But this implies $\varrho(a) \in Z(J)$, a contradiction. 
We conclude $\varrho(b) \in R$.  
\end{proof}

\begin{coro}\label{Jnoskew field} If $R=J$ and $R$ is a skew field, then $R$ is commutative.
\end{coro} 
\begin{proof} 
If $J=R$ is a non-commutative skew field, then we can apply \ref{R=S} and get $R=S$. 
But if $b \in A \setminus A_0$, there is a $a \in A_0$ with $\varrho(a) = -\varrho(b) $ and so 
$\varrho({ab}) = \varrho(a) +\varrho(b) =0$, but $ab \ne 1$, a contradiction.
 \end{proof}
 
\begin{prop}\label{A0} 
\begin{enumerate} 
\item If $\characteristic R \ne 2$, then  
$C_A([V,A]) = A_0 = C_A(V/C_V(A))$.
\item If $\characteristic R =2$, then $A_0 =C_A([V,A]) = C_A(V/C_V(A))$ or $R=J$ is a 
commutative field. 
\end{enumerate}
\end{prop}
\begin{proof} Set $A_1:=C_A([V,A])$ and $A_2:= C_A(V/C_V(A))$ and define $B_1,B_2 \leq B$ 
analogously.  
By \ref{rootgroups} $A_1$ and $A_2$ are both root subgroups of 
$A$ which contain $A_0$ and which act quadratically on $V$. Thus $A_1$ and $A_2$ are special.
\begin{enumerate}\item 
Since $A_0$ is a 
$H$-invariant subgroup of $A_0$, \ref{invariant subgroups} implies $A_1 = A_0 =A_2$ if 
$\characteristic V \ne 2$. 
\item
Suppose $\characteristic V =2$ and let $a \in A_1$. Then $f(\ob,\oa) = 0$ for all $b \in A$. If 
$a$ is not in $A_2$, then 
there is an element $b \in A$ with $f(\oa,\ob) \ne 0$. Thus $0 \ne \varrho([a,b]) = f(\oa,
\ob)$. 
Now for $c,d \in A_0$ and $h = h([a,b])^{-1} h(c) h(d)$, we get 
$$\varrho(c) \varrho(d) =f(\oa,\ob) \varrho (h) = f(\oa,\ob^h) = f(\oa,\ob^h)-f(\ob^h,\oa)=
\varrho([a,b^h]) \in J.$$
Thus $J$ is multiplicatively closed and so $R=J$. Since every element in $J^{\#}$ is a unit, 
$J$ is a skew field.
We can apply \ref{Jnoskew field} and conclude that $R=J$ is a commutative field if $A_1 \ne 
A_0$.

For $a \in A_2$ and $b \in A$ we have $f(\ob,\oa) =\varrho(h_{[b,a]})$ and thus can use the same 
argument.
\end{enumerate}
\end{proof}

\bigskip
We will now define a $R$-module structure on $\oA$. 
For $a \in A$ and $h_1, \ldots, h_n \in H_0$ set $\oa^{\sum_i h_i} := \sum_i \oa^{h_i}$. 
\begin{prop}\label{module}
\begin{enumerate}
\item The map $\cdot: \oA \times R \to \oA$ 
defined by $\oa \cdot (\sum_i \varrho(h_i)) = \oa^{\sum_i h_i}$ for $h_i \in H_0$ is well-defined 
and defines 
a $R$-module-structure on $\oA$.
\item For all $a,b\in A$ and $r\in R$ we have $f(\oa,\ob\cdot r)=f(\oa,\ob)r$. 
\end{enumerate}
 
\end{prop}
\begin{proof} For part (a) we only have to show that if $h_1, \ldots, h_n \in H_0$ with $\sum_{i=1}^n 
\varrho(h_i) =0$, then 
also $\oa^{\sum_{i=1}^n h_i} = \overline{0}$ for all $a \in A$. Suppose that $h_0, \ldots, h_n 
\in H_0$ are chosen this way. 
Then for all $a,b \in A$ we get by using \ref{propf}(c)
$$0 =f(\oa,\ob)\sum_{i=1}^n \varrho(h_i) =f(\oa,\sum_{i=1}^n \ob^{h_i}) = f(\oa, \ob^{\sum_{i=1}^n 
h_i}).$$
Thus $b^{h_1} \ldots b^{h_n} \in C_A([V,A])$. If $A_0 = C_A([V,A])$, then $\ob^{\sum_{i=1}^n 
h_i} =\overline{0}$.
If $A_0 \ne C_A([V,A])$, then $R$ is a field of characteristic $2$ by \ref{A0}. Thus $H_0$ is 
abelian and 
$\varrho(h^{-\mu})=\varrho(h)$ for all 
$h \in H_0$ by \ref{-mu}(a). Thus 
$$0=\sum_{i=1}^n \varrho(h_i) f(\ob,\oa) = \sum_{i=1}^n \varrho(h_i^{-\mu}) f(\ob,\oa) = \sum_{1=1}^n 
f(\ob^{h_i},\oa) = 
f(\ob^{\sum_{i=1}^n h_i},\oa).$$
Hence we also get $b^{h_1} \ldots b^{h_n} \in C_A(V/C_V(A))$ and so (a) follows. Part (b) follows again 
by \ref{propf}(c). 
\end{proof}
\begin{coro}
If $\characteristic R \ne 2$, then $G$ is quasi-simple and thus generated by the conjugates of 
$A_0$.
\end{coro}
\begin{proof} If $\characteristic R \ne 2$, then $[\oA,H] =\oA$ by \ref{module}. 
Moreover, $A^{\prime}$ is a 
$H$-invariant subgroup of $A_0$ and 
thus by \ref{invariant subgroups} either $A_0 =A^{\prime}$ or $A^{\prime}=1$. This shows $A \leq 
G^{\prime}$ and thus $G$ is perfect. 
Therefore $G$ is quasi-simple by \cite[I (1.10)]{T1}.  
The conjugates of $A_0$ generate a normal subgroup of $G$ and so the last claim follows.
\end{proof}
\begin{lemma}\label{abelian}
If $A$ is abelian, then $R/\B(R)$ is a commutative field of characteristic $2$.
\end{lemma}
\begin{proof}  
By \ref{G proper} $G$ is proper, thus by \ref{Segev} $G$ 
is special. Since $A_0$ is a $H$-invariant subgroup of $A$, \ref{invariant subgroups} 
implies that $A$ is an elementary-abelian 
$2$-group, so $R$ has characteristic $2$. Since $A$ is abelian, we have 
$[a,b]=1$ and thus $f(a,b) =f(b,a)$ for all $a,b \in A$ by \ref{propf}(d). 
We set  $$I:=\{r \in R| f(a,b)r =0 \text{ for all } a,b \in A\}.$$ Then $I$ is a right ideal of 
$R$. If $r \in I, 
s \in R$ and $a,b \in A$, then $f(a,b)sr = f(a,b\cdot s) r = 0$ and thus $I$ is an ideal of 
$R$. Of course $I \ne R$, 
otherwise $f(a,b)=0$ for all $a,b \in A$ and thus $G$ would act quadratically on $V$. Let $a,b 
\in A$ and 
$r,s \in R$. Then
\begin{align*} f(a,b)rs = f(a,b \cdot r)s = f(b \cdot r,a) s =f(b\cdot r,a\cdot s) =f(a\cdot s,b\cdot r) 
\\=f(a\cdot s,b)r =
f(b,a\cdot s)r =f(b,a) sr = f(a,b) sr.\end{align*}
Thus $rs-sr \in I$. Hence $R/I$ is commutative. Since $R/I$ is an envelope for $J$, $J$ is 
commutative. Thus the 
universal semi-prime envelope of $J$ is a commutative field by \ref{semi-prime}. Since 
$R/\B(R)$ is a semi-prime envelope of $J$, the 
claim follows.
\end{proof}

For $n \in \ZZ$ set $\lambda_n := h({e^n})$. Note that $\varrho(\lambda_n) = n\varrho(e) =n$ 
by \ref{rho}, that $a^{\lambda_n} =a^{n^2}$ for all 
$a \in A_0$ by \cite[4.6(6)]{DMS} and that $\lambda_{-1}=\mu\mu(e^{-1})=\mu^2$.
\begin{prop}\label{divisible} Suppose that the characteristic of $V$ is 
not $2$. Let $a \in A^*$. 
Then there is an element $\widehat{a} \in \oa$ 
such that $\widehat{a} ^n =\widehat{a}^{\lambda_n}$ for all $n \in \ZZ$ with $n$  
relatively prime to the characteristic of $A_0$. Moreover, $\widehat{a}$ is special, 
$\mu(\widehat{a})^2 =1$, $\widehat{a^n}=\widehat{a}^n$ for all $n\in \ZZ$ relatively prime to 
the characteristic 
of $V$, 
$\widehat{a^h} =\widehat{a}^h$ for all $h\in H$ and $\widehat{a}$ is  
the unique element in $\oa$ which is inverted by $\mu^2=\lambda_{-1}$.   
\end{prop}
\begin{proof} 
We have $-\oa = \oa\cdot {\lambda_{-1}} = \oa^{\lambda_{-1}}$, so 
$x :=  a^{\lambda_{-1} } a \in A_0$.  
Since the characteristic of $A_0$ does not divide $2$, there is a $y \in A_0$ with 
$y^2 =x^{-1}$. 
Set $\widehat{a}:=ya$. Then $$ \widehat{a}^{\lambda_{-1} } = 
y^{\lambda_{-1}}a^{\lambda_{-1}} 
=yxa^{-1}  =yy^{-2} a^{-1}  = y^{-1} a^{-1} =a^{-1}y^{-1}= (ay)^{-1} =\widehat{a}^{-1}.$$
We claim that $\widehat{a}$ is the unique element in $A_0 a$ which is inverted by 
$\lambda_{-1}=\mu^2$. Indeed, for $b\in A_0$ we have $(\widehat{a}b)^{\mu^2} = 
\widehat{a}^{\mu^2} b^{\mu^2} =\widehat{a}^{-1} b$ and $(\widehat{a}b)^{-1} =b^{-1} 
\widehat{a}^{-1}
=\widehat{a}^{-1} b^{-1}$, so we have equality if and only $b^2=1$. Since 
$A_0$ has no element of order $2$, the claim follows.
\\
Next we claim that $\mu(\widehat{a})^2=1$ and that 
$\varrho(\widehat{a})=\varrho(\widehat{a}^{-1})$. Since 
$\lambda_{-1} =\mu^2  \in Z(H)$ by \ref{mu quadrat},
we have \begin{align*} \mu(\widehat{a}^{-1}) =\mu(\widehat{a}^{\mu^2}) 
=\mu(\widehat{a})^{\mu^2} =
\mu^{-2}\mu(\widehat{a})\mu^2= \mu^{-3}h(\widehat{a})\mu^2=\mu^{-1}h(\widehat{a})=
\mu(\widehat{a})\end{align*}
and thus $\varrho(\widehat{a}^{-1}) =\varrho(\widehat{a})$. \\
Now we claim that $\widehat{a}$ is special. Let $\breve{a}$ be the unique element in 
$A_0 (\widehat{a}\diamond \mu)$ which is inverted by $\mu^2$. Then there is an element $b\in A_0$ 
with $\breve{a} =(\widehat{a}\diamond \mu) b$. Since we have $\mu(\widehat{a}\diamond \mu) =\mu(\widehat{a}^{-1})^{\mu} 
=\mu(\widehat{a})^{\mu} = \mu(\widehat{a}\diamond \mu)^{-1}$ by \ref{formulas}(b), we get
\begin{align*}
\varrho(\breve{a})+\varrho(b)   =\varrho(\breve{a}b) =\varrho(\widehat{a}\diamond \mu) =
\varrho((\widehat{a}\diamond \mu)^{-1}) = 
\\ \varrho((\breve{a}b)^{-1}) =\varrho(\breve{a}^{-1} b^{-1}) = \varrho(\breve{a}^{-1}) +\varrho(b^{-1}) =
\varrho(\breve{a}) -\varrho(b)
\end{align*}
and therefore $2\varrho(b)=0$. Since $\characteristic V \ne 2$, we get $\varrho(b)=0$ and thus $b=1$. 
Therefore we have $(\widehat{a}\diamond \mu)^{-1} = (\widehat{a}\diamond \mu)^{\mu^2} =(\widehat{a}^{\mu^2})\diamond \mu =
(\widehat{a}^{-1})\diamond \mu$, hence $\widehat{a}$ is special. \\
Since $\mu^2 \in Z(H)$, for $h\in H$ we have $(\widehat{a}^h)^{\mu^2} = (\widehat{a}^{\mu^2})^h 
=(\widehat{a}^{-1})^h =(\widehat{a}^h)^{-1}$, thus $\widehat{a^h}=\widehat{a}^h$. \\
For $n$ relatively prime to the characteristic of $A_0$ we 
have $$ A_0 \widehat{a^{\lambda_n}} = A_0 a^{\lambda_n}= A_0 a^n =A_0   
\widehat{a}^n.$$
Since $(\widehat{a}^n)^{\mu^2} = (\widehat{a}^{\mu^2})^n =\widehat{a}^{-n}$, we get
$\widehat{a}^{\lambda_n} =\widehat{a^{\lambda_n} }= \widehat{a}^n$. 
Now $\widehat{a}^n$ is the unique element in $A_0a^n=A_0\widehat{a}^n$ which is inverted by $\lambda_{-1}$, 
thus we have $\widehat{a^n}=\widehat{a}^n$. 
\end{proof}
\begin{lemma}\label{f h 2} Suppose that $\characteristic V \ne 2$. 
For all $a \in A \setminus A_0$ we have
\begin{enumerate}
\item $f(a,a)=2\varrho(\widehat{a})$.
\item $h(\widehat{a})^{-\mu}=\mu^2 h(\widehat{a})$ and 
$-\varrho(\widehat{a}^{-1})=\varrho(h(\widehat{a})^{-\mu})=-\varrho(\widehat{a})$.
\item $f(a^{h(\widehat{a})},c) =\varrho(\widehat{a}) f(a,c)$ and $f(c,a^{h(\widehat{a})})=-
f(c,a)\varrho(\widehat{a})$.
\end{enumerate}
\end{lemma}
\begin{proof}
\begin{enumerate}
\item We have 
$\varrho(\widehat{a}^2)=\varrho(\widehat{a}^{\lambda_2})=\varrho(\lambda_2^{-\mu})\varrho(\widehat{a})
\varrho(\lambda_2) =4\varrho(\widehat{a})$, so $f(a,a)=f(\widehat{a},\widehat{a})=
\varrho(\widehat{a}^2)-2\varrho(\widehat{a})=2\varrho(\widehat{a})$ by \ref{handf}(a).
\item We have $h(\widehat{a})^{-\mu} = (\mu(\widehat{a})^{-1}\mu^{-1})^\mu=\mu^{-1}\mu(\widehat{a})
=\mu^2 h(\widehat{a})$, hence by \ref{-mu} $-\varrho(\widehat{a}^{-1}) 
=\varrho(h(\widehat{a})^{-\mu})=\varrho(\mu^2)\varrho(\widehat{a})=-
\varrho(\widehat{a})$. 
\item This follows from (a), (b) and \ref{f h} for $b=\widehat{a}$. 
\end{enumerate}
\end{proof}

\begin{lemma}\label{root groups 2} 
Let $G_1$ be a subgroup of $G$ which contains $G_0$. Set $A_1 = A \cap G_1$ and $B_1 =B \cap 
G_1$. Then $\langle A_1,B_1
\rangle$ is a rank one group which is normal in $G_1$.
\end{lemma}
\begin{proof}
By Lemma \ref{A_0=2} we may assume that $|A_0| >2$ and thus $H_0 \ne 1$. 
Let $C_0$ be a subgroup of $G_1$ which is conjugate 
to $A_0$. Then there is an element $a \in A$ with $B_0^a =C_0$. We claim that $a \in A_1$. Set 
$P=\langle A_0,C_0 \rangle$ 
and $K =N_P(A_0) \cap N_P(C_0)$. Then $P =G_0^a $ and $K = H_0^a$. If $\characteristic V\ne 2$, set 
$h=\lambda_{-1}$ and $j=\lambda_2$, and if $\characteristic V=2$ let $h=h(b)$ for an element $b 
\in A_0 \setminus \{1,e\}$ and set $j=h(be)$. In both cases we have 
$\overline{x}\overline{x}^{-h} =\overline{x}^j$ for all $x \in A$. We have $h^a \in K \leq P \leq 
G_1$ and 
so $a^{-h} a= [h,a] =h^{-1} h^a \in G_1$. Thus $\overline{a} = \overline{a}^{hj^{-1}} \in 
A_1/A_0$. 
Since $A_0 \leq G_1$, we get $a \in G_1 \cap A =A_1$. \\
We get $A_0^{G_1} = \{A_0\} \cup B_0^{A_1}$ and by a similar argument $A_0^{G_1} = \{B_0\} \cup 
A_0^{B_1}$. 

Now if $a\in A_1^*$, then there is $b\in B_1^*$ with $A_0^b=B_0^a$. Thus $1\ne A_0^b =B_0^a \leq A^b \cap 
B^a$, which implies $A^b =B^a$. It follows $A_1^b=(A \cap G_1)^b =A^b \cap G_1^b = A^b \cap G_1 =
B^a \cap G_1^a = (B\cap G_1)^a =B_1^a$, therefore $\hat{G}_1:=\langle A_1,B_1\rangle$ is a rank one group with 
unipotent subgroups $A_1$ and $B_1$. We have $\hat{G}_1 \trianglelefteq G_1$ because $\hat{G}_1$ is 
generated by the $G_1$-conjugates of $A_1$. 
\end{proof}

We will now show that the map $\varrho(h) \mapsto \varrho(h^{-\mu})$ extends uniquely to an 
anti-automorphism $*$ of $R$. Note that if $H_0$ is abelian (and hence $R$ is commutative), 
then  
by \ref{-mu}(b) $\varrho(h^{-\mu})=\varrho(h)$ for all $h\in H_0$ (since $H_0$ is generated by the elements $h(a)$ with $a \in A_0^{\#}$), so in this case the claim holds for $*$ the identity. 
 \begin{prop}\label{antiauto} Suppose that $A$ is not abelian or $\B(R)=0$. 
There is a unique involutory anti-automorphism $*$ of $R$ with 
$\varrho(b)^* = \varrho(h(b)^{-\mu})$ for all $b \in A_0$. 
\end{prop}
\begin{proof}
The map $h \mapsto h^{-\mu}$ is an anti-automorphism of $H_0$. Since $R$ is generated by 
$\varrho(H_0)$, 
there is at most one possibility to extend this map to an anti-automorphism of $R$. We have to 
show that 
if $h_1, \ldots, h_n \in H_0$ with $\sum_{i=1}^n \varrho(h_i) =0$, then also $\sum_{i=1}^n 
\varrho(h_i^{-\mu}) =0$. 
In this case, the map $*$ with $(\sum_{i=1} \varrho(h_i))^* =\sum_{i=1}^n \varrho(h_i^{-\mu})$ for 
$h_1, \ldots, h_n \in H_0$ is a well-defined anti-automorphism of $R$. 
\newline Suppose that $h_1, \ldots, h_n \in H_0$ with $\sum_{i=1}^n \varrho(h_i) =0$. 
By the remark above and by \ref{abelian} we may assume that $A$ is not abelian. Thus 
by \ref{propf}(f) and by \ref{f h 2}(a) there is an element $a \in A$ such 
that $f(a,a)$ is a unit in $R$.
Therefore we get by \ref{propf}(c)
$$0= f(\overline{0},\oa) =f(\oa \cdot \sum_{i=1}^n\varrho(h_i), \oa) =\sum_{i=1}^n f(\oa \cdot 
\varrho(h_i), \oa) =
\sum_{i=1}^n \varrho(h_i^{-\mu} ) f(\oa,\oa) .$$
This implies $\sum_{i=1}^n \varrho(h_i^{-\mu}) =0$, as desired. Since $\mu^2 \in Z(H_0)$, we get 
$\varrho(a)^{**} =
\varrho(a)$ for all $a \in H_0$. Thus $*$ is an involution.
\end{proof}
\begin{lemma}\label{h in R} Suppose that $R=S$ and $\characteristic R\ne 2$. Then we have
$\varrho(\widehat{a})^*=-\varrho(\widehat{a})$ for all $a\in A \setminus A_0$.
\end{lemma}
\begin{proof}
We have by \ref{f h 2}(c)
\begin{align*} f(\oa^{h(\widehat{a})},\ob)=\varrho(\widehat{a}) f(\oa,\ob)=f(\oa\cdot\varrho(\widehat{a})^*,\ob)
\end{align*}
for all $b\in A_0$. Thus we have $\oa^{h(\widehat{a})}=\oa\cdot \varrho(\widehat{a})^*$. Again by \ref{f h 2}(c) 
we get
\begin{align*}
f(\oa,\oa)\varrho(\widehat{a})^*=f(\oa,\oa\cdot \varrho(\widehat{a})^*)=
f(\oa,\oa^{h(\widehat{a})})=-f(\oa,\oa)\varrho(\widehat{a}).
\end{align*} 
Since $f(\oa,\oa)=2\varrho(\widehat{a})$ is a unit in $R$, we conclude that $\varrho(\widehat{a})^*=
-\varrho(\widehat{a})$ holds.
\end{proof}
\begin{prop}\label{f R linear} 
\begin{enumerate}
\item The map $f: \oA \times \oA \to S$ is $*$-sesquilinear. 
\item If $J$ is not a commutative field of characteristic $2$, then $f$ is non-degenerate.
\item If $J$ is commutative, then the map $g: \oA \times \oA \to R: (\oa,\ob) \mapsto f(a,b)-f(b,a) =\varrho({[a,b]})$ 
is an alternating $R$-bilinear map. 
\end{enumerate}
\end{prop}
\begin{proof}
\begin{enumerate}
\item This follows by \ref{propf}(a) and (c).
\item This follows by \ref{A0}.
\item If $J$ is commutative, then $*$ is the identity, so $f$ is $R$-bilinear. By \ref{propf}(e), 
$g(\oa,\ob) = f(a,b) -f(b,a)$, so $g$ is also $R$-bilinear. It is clear that $g$ is 
alternating.
\end{enumerate}
\end{proof}
\begin{lemma}\label{ample} If $\characteristic R=2$ and $A$ is not abelian, then $J$ is an ample subspace 
in $\mathcal{H}(R,*)$.
\end{lemma}
\begin{proof} Since $A$ is not abelian, then by the choice of $e$ 
there is an element $a \in A$ with $a^2 = e$. Thus $f(\oa,\oa) =\varrho(e) = 1$ by \ref{propf}(f).
If $r \in R$ and $b \in A$ with $\ob = \oa \cdot r$ we get $r^* r = f(\oa \cdot r, \oa \cdot r) = 
f(\ob,\ob) = \varrho(b^2) \in J$. 
This also implies $$r+r^* = (r+1)^* (r+1) + r^* r +1 \in J.$$ 
Furthermore, for all $a,b \in A_0$ we have 
$\varrho(b)^*\varrho(a)\varrho(a)=\varrho(h(a^{h(b)}))\in J$ 
by \ref{-mu}(c) and 
\ref{antiauto}. Since $J$ generates $R$ as a ring, the claim follows by \ref{amplifier}.
\end{proof}  
\begin{theorem}\label{Rskew field}
 \begin{enumerate}
\item If $\characteristic R \ne 2$, then $R$ is a skew field.
\item If $\characteristic R = 2$, then $R/\B(R)$ is a skew field.
\end{enumerate}
\end{theorem}
\begin{proof}
\begin{enumerate}
\item If $r \in R$ and $a \in A \setminus A_0$ with $\oa \cdot r =\overline{0}$, then $0 =f(\oa,\oa \cdot 
r) = f(\oa,\oa) r $. But $f(\oa,\oa)$ is a unit in $R$ by \ref{f h 2}, 
thus $r=0$. So if $r \ne 0$, then $\oa \cdot r \ne 0$ and thus 
$f(\oa \cdot r, \oa \cdot r)$ is again a unit. We get $r^* f(\oa, \oa) r = f(\oa \cdot r, \oa \cdot r)$ and 
$r f(\oa,\oa) r^* = f(\oa \cdot r^*, \oa \cdot r^*)$. Hence $r$ is invertible. 
\item If $A$ is abelian, then the claim follows by \ref{abelian}. 
If $A$ is not abelian, then $J$ is an ample subspace in $\mathcal{H}(R,*)$. By \cite[Theorem 2.1.8]{H}, 
$R/\B(R)$ is either a skew field, or $R/\B(R)$ is commutative and $*$ induces the identity on $R/\B(R)$, 
or $R/\B(R)$ is the direct product of a skew field and its opposite 
and $*$ induces the exchange involution on $R/\B(R)$, or $R/\B(R) \cong \mathrm{Mat}_2(F)$ for a 
commutative field $F$ 
and $*$ induces the standard involution on $R/\B(R)$, which is given by 
$$\left( \begin{array}{cc} x & y \\ z & w \end{array} \right)^ * = \left( \begin{array}{cc} w & -y \\ -z & 
x \end{array} \right)=\left( \begin{array}{cc} w & y \\ z & x \end{array} \right)$$ 
for $x,y,z,w \in F$. 
If $R/\B(R)$ is commutative, $J$ must be commutative as well. Since $R/\B(R)$ is the semi-prime envelope of 
$J$, 
$R/\B(R)$ must be a field. \\
Suppose that $R/\B(R)\cong \mathrm{Mat}_2(F)$. 
Let $\oJ$ be the image of $J$ in $\oR:=R/\B(R)$. Since $J$ contains all traces, 
$\Cent\oR \subseteq \oJ$. Note that $\oJ$ cannot be $\Cent \oR$ because $J$ generates $R$. 
Therefore there are 
$x,y \in F $ with  
$$0 \ne \left( \begin{array}{cc} 0 & x \\ y & 0 \end{array} \right) \in \oJ.$$
Since every element in $J^{\#}$ is invertible, neither $x$ nor $y$ can be $0$. 
But since $J$ is ample in $R$, we get
$$ \left( \begin{array}{cc} 0 & 0 \\ 1 & 0  \end{array} \right) \left( \begin{array}{cc} 0 & x \\ y & 0 
\end{array} \right) \left( \begin{array}{cc} 0 & 0 \\ 1 & 0 \end{array} \right) =
\left( \begin{array}{cc} 0 & 0 \\ x & 0 \end{array} \right) \in \oJ,$$ 
a contradiction, since every element in $J^{\#}$ is invertible. Thus this case cannot hold.
\\
Suppose that $R/\B(R) =K \times K^o$ for a skew field $K$. We denote the anti-automorphism induced by $*$ 
on $K \times K^o$ also by $*$. There is an automorphism $\phi$ of $K$ such that 
$(x,y)^* = (y^{\phi^{-1}},x^{\phi})$ for $x \in K,y \in K^o$. Thus the image of
$J$ is the set $\{(x,x^{\phi})|  x \in K\}$. If $b \in A \setminus A_0$, 
then $\overline{\varrho(b)}$ is mapped on $(x,y)$ with $x \in K, y \in K^o$. But then there is an 
element $a \in A_0$ 
such that $\varrho(a)$ is mapped on $(x,x^{\phi})$ and so $\varrho({ba})$ is mapped on $ (0,x^{\phi} +y)$, a contradiction, 
since $\varrho({ba})$ is invertible.
\end{enumerate}
\end{proof}

We summarise our results:
\begin{coro}\label{cases} If that $A_0 \ne 1$, then one of the following cases holds:
\begin{enumerate}
\item $J \subseteq R$ is commutative, thus $R/\mathfrak{B}(R)$ is a field. If $\characteristic 
R \ne 2$, then $R=J$ is a field, $R\subseteq \Cent S$ and $G_0 \cong \SL_2(R)$.
\item $R/\mathfrak{B}(R)$ is a skew field with involution $*$. 
If $\characteristic R =2$, then  
$(R/\mathfrak{B}(R),\oJ,*)$ is an involutory set, where $\oJ$ denotes the image of $J$ in 
$R/\B(R)$. If $\characteristic R \ne 2$, then 
$\mathfrak{B}(R) =0$.
\end{enumerate} 
\end{coro}
\begin{proof}
This follows by \ref{R=S}, \ref{ample} and \ref{Rskew field}.
\end{proof}
In the most cases $\B(R)$ is actually $0$ and so $R$ is a skew field. 
\begin{prop} If $J$ is not commutative and $\B(R) \ne 0$, then either $ R/\B(R)$ is a biquaternion algebra, 
$*$ induces a symplectic involution on $R /\B(R)$ and $J \cong K_*=\{a+a^ *| a \in K\}$, or $R/\B(R)$ is a 
quaternion algebra.
\end{prop}
\begin{proof} Set $K:=R/\B(R)$ and let $K_0$ be the image of $J$ in $K$. Then there is a Jordan isomorphism 
$\phi: K_0 \to J$. If $Z_{48}(K_0) \ne 0$, then $\phi$ can be extended to an associative homomorphism 
${\widehat \phi}: K \to R$ by the $Z$-algebra Theorem in \cite{Mc}. Since $J \subseteq \im({\widehat \phi})$ 
and $J$ generates $R$, we have 
$R = \im(\phi) \cong K$. Now \cite[2.6.5]{Kn} tells us that $Z_{48}(K_0)=0$ implies that either 
$K$ is a biquaternion algebra, $*$ a symplectic involution and $K_0=K_*$, or $K$ is a quaternion algebra.
\end{proof}

\begin{prop} Suppose that $R=S$ and that $A$ is not abelian. If $v \in 
C_V(A)$ with $\Phi(v,a) 
=0$ for all 
$a\in A$, then $v=0$.
\end{prop}
\begin{proof} Let $W =vR$. Then $W$ is $H$-invariant. If $r 
\in R$, then 
$r =\sum_{i=1}^n \varrho(h_i)$ with $h_i\in H_0$. Thus for all $a \in A$ we 
get by \ref{eigenschaftenphi}(b)
$$\Phi(vr,a) =\Phi(\sum_{i=1}^n v h_i, a) = \sum_{i=1}^n \Phi(vh_i,a) = 
\sum_{i=1}^n 
\Phi(v,a^{h_i^{-\mu}}) =0.$$
Thus $ W \oplus W\mu =X(W)$ is a quadratic $G$-module, so $A$ must be abelian, a contradiction to our 
assumption.
\end{proof}

Suppose that $\characteristic R =2$ and 
$J$ is not commutative. Then $R=S$ and $R/\B(R)$ is a skew field. Moreover, 
$A$ is not abelian, 
so there is $c \in A$ with $c^2 =e$ and thus 
$f(c,c) =1$. Let $W \leq C_V(A)$ be a finitely generated $R$-module. Note that
$C_V(G) =0$ implies that $C_{X(W)}(G) =0$ and that by \ref{X(W)}(d) we have $X(W)=[X(W),G]$.

Define 
$$A_1 :=\{a \in A| f(\oa,\ob) \in \B(R) \text{ for all } b \in A\},\ H_1:= 
\{\mu \mu(a)| a \in 
A_1\}$$ and 
$$\oX(W) := X(W)/(X(W\B(R))+\Phi(W,A_1)).$$ Then we have:
\begin{prop} If $W \ne 0$, then $0 \ne \oX(W)$ is a cubic module for $G$ 
with 
$A_1 =C_A([\oX(W),A]) \cap C_A(\oX(W)/C_{\oX(W)}(A))$, 
$\oX(W) =[\oX(W),G]$ and $C_{\oX(W)}(G)=0$. If $R^ {\prime}$ is the subring 
of $\End(\oX(W))$ 
generated by 
$H_1$, then $R^{\prime} = R/\B(R)$. 
\end{prop} 
\begin{proof} First note that $[\Phi(W,A_1),A] \leq W\B(R)$ and $\Phi(W,A_1) 
\leq C_V(\mu)$; this implies 
that $X(W\B(R)) + \Phi(W,A_1)$ is a $G$-submodule of $X(W)$.  
If $W \ne 0$, then $W \ne W\B(R)$ by the Nakayama Lemma (\cite[13.12]{I}, note that $\B(R) 
\subseteq J(R)$). 
Thus we also get $X(W) \ne X(W\B(R))+\Phi(W,A_1)$. Now $X(W) =[X(W),G]$ and 
thus 
$[\oX(W),G] =\oX(W)$. Thus $[\oX(W),G,G,G] \ne 0$ and $[\oX(W),A,A,A] =0$, 
so $\oX(W)$ is a cubic module 
for $G$.
\newline Suppose that $z \in C_{\oX(W)}(G)$. Since $K:=R/\B(R)$ is a skew 
field,
 $W/W\B(R)$ is a finite-dimensional vector space over $K$. Let $\{{\ow}_1, 
 \ldots,{\ow}_n\}$ 
be a $K$-basis of $\oW$ and let $w_i$ be a preimages of $\ow_i$ in 
$W$ for $i=1,\ldots,n$. Then 
$W= w_1 R + \ldots + w_n R +W\B(R)$ and so using Nakayama's lemma again we get 
$W =w_1 R + \ldots +w_n R$. 
Thus there are $x,y \in W, a_1, \ldots a_n \in A$ 
with $$z = x +\Phi(w_1,a_1) +\ldots + \Phi(w_n,a_n) + y\mu + X(W\B(R)) + 
\Phi(W,A_1).$$ We get 
$$0=[z,e] =[y\mu,e] + X(W\B(R)) + \Phi(W,A_1)$$ and thus $y=[y\mu,e] \in 
W\B(R)$. 
If we apply $\mu$, we get $x \in W\B(R)$. For all $a \in A$ we have 
$$0= [z,a] = w_1 f(\oa_1,\oa) + \ldots + w_n f(\oa_n,\oa) + X(W\B(R)) + 
\Phi(W,A_1)$$ 
and therefore
$$w_1 f(\oa_1,\oa) + \ldots + w_n f(\oa_n,\oa) \in W\B(R).$$
Since $\ow_1,\ldots,\ow_n$ are $K$-linearly independent, this implies 
$f(\oa_1,\oa),\ldots, f(\oa_n,\oa) \in \B(R)$. Thus $a_1,\ldots,a_n \in A_1$ 
and so $z=0$.
\newline 
By \ref{X(W)}(b) one sees that $C_A(\oX(W))$ is just the image of $W$ in 
$\oX(W)$. Thus 
$C_A(\oX(W)/C_{\oX(W)}(A)) =A_1$. Since $J$ is not commutative, we also have 
$A_1 =C_A([\oX(W),A])$ by 
\ref{A0}.

Since $H_1 \leq H$, the ring $R^{\prime}$ is contained in the image of $R=S$ 
in $\End(\oX(W))$. 
Now $\B(R)$ annihilates $\oX(W)$ by construction. Thus we get $R^{\prime} 
=R/\B(R)$ which is a skew field. 
\end{proof}

This proposition shows that we can assume that $R$ is a skew field if $J$ is not commutative. 
\chapter{Cubic rank one groups with hermitian quadratic kernel}
In this chapter we will assume that $J$ is not commutative and that $\B(R) =0$. Thus 
$R$ is a skew field with involution $*$. We have $J \subseteq \mathcal{H}(R,*)$. If $\characteristic R =2$, then $J$ is ample in 
$R$. 
\begin{lemma}\label{anisotropic} \begin{enumerate}
\item If $a \in A$, then $\varrho(a)\in J$ iff $a \in A_0$.
\item If $\characteristic R \ne 2$, then $\varrho(a) \in \mathcal{H}(R,*)$ iff $a \in A_0$. 
\end{enumerate}
\end{lemma}
\begin{proof}
\begin{enumerate}
\item Let $a \in A$ with $\varrho(a) \in J$. Then there is $b \in A_0$ with $\varrho(b) =-\varrho(a)$ and thus 
$\varrho({ab}) = \varrho(a) +\varrho(b) =0$. Thus $ab=1$ and $a =b^{-1} \in A_0$. 
\item If $a\in A \setminus A_0$, then there is $b\in A_0$ with $a=\widehat{a}b$. 
By \ref{h in R} we get 
$$\varrho(a)^*=\varrho(\widehat{a}b)^*=\varrho(\widehat{a})^*+\varrho(b)^*=-\varrho(\widehat{a})+\varrho(b).$$
Thus the claim follows.
\end{enumerate}
\end{proof}
\begin{lemma} We have $\varrho(h^{-\mu}) =\varrho(h)^*$ for all $h \in H$.
\end{lemma}
\begin{proof} We consider the action of $G_0 H$ on $X:=C_V(A) \oplus C_V(B)$. We regard $R$ as a subring of 
$E:=\End(C_V(A))$. 
The map $\mu$ induces an automorphism between $C_V(A)$ and $C_V(B)$. Thus we can regard $X$ as the direct 
sum 
of two standard $E$-modules. 
Therefore we get a homomorphism $\xi$ from $G_0 H$ in the group of all invertible $2 \times 2$-matrices 
over $R$, such that 
the image of $a \in A_0$ is 
$$ \left( \begin{array}{cc} 1 & 0 \\ \varrho(a) & 1 \end{array} \right),$$
while $\mu$ is mapped to 
$$ \left( \begin{array}{cc} 0 & 1 \\ -1 & 0 \end{array} \right)$$ and 
$h(a)$ is mapped to 
$$ \left( \begin{array}{cc} \varrho(a) & 0 \\ 0 & \varrho(a)^{-1} \end{array} \right) =
 \left( \begin{array}{cc} \varrho(a) & 0 \\ 0 & \varrho(a)^{-*} \end{array} \right).$$
For $b \in A$ there is an element $y \in R$ such that the image of $h(b)$ is 
$$ \left( \begin{array}{cc} \varrho(b) & 0 \\ 0 & y \end{array} \right).$$ 
We get 
$$ \xi(h(b)^{-1} a h(b)) = \left( \begin{array}{cc} 1 & 0 \\ y^{-1} \varrho(a) \varrho(b) & 1  \end{array} 
\right).
$$
Thus we have $$y^{-1} \varrho(a) \varrho(b) = (y^{-1} \varrho(a) \varrho(b))^* = \varrho(b)^* \varrho(a) 
y^{-*}$$ and hence
$$\varrho(a) = y\varrho(b)^* \varrho(a)  y^{-*} \varrho(b)^{-1} = (y\varrho(b)^*) \varrho(a) 
(y\varrho(b)^*)^{-*}$$
for all $a \in A_0$. If we take $a=e$, we get $(y\varrho(b)^*)^{-1} = (y\varrho(b)^*)^{-*}$ and thus 
$(y\varrho(b)^*)^* = y\varrho(b)^*$. Therefore we have $(y\varrho(b)^*) \varrho(a) (y\varrho(b)^*)^{-1} = 
\varrho(a)$ for all 
$a \in A_0$. Since $R$ is generated by $\varrho(H_0)$, we conclude $y\varrho(b)^* \in k:=\mathcal{H}(\Cent 
R,*)$. 
Thus for 
all $b \in A^*$ there is a $\omega_b \in k$ with 
$$\xi(h(b)) = \left( \begin{array}{cc} \varrho(b) & 0 \\ 0  & \varrho(b)^{-*} \omega_b \end{array} \right).
$$
Of course, $\omega_a =1$ for all $a \in A_0$. 
By \ref{-mu}(a) we have $\varrho(h(b)^{-\mu}) = -\varrho(b^{-1})$ 
and thus
\begin{align*}
\left( \begin{array}{cc} -\varrho({b^{-1}}) & 0 \\ 0 & -\varrho({b^{-1}})^{-*} \omega_{b^{-1}} 
\end{array} \right) = \\
\xi(\mu^{-2}) \xi(h({b^{-1}})) = \xi(h(b)^{-\mu}) 
= \xi(\mu )^{-1} \xi(h(b))^{-1} \xi(\mu) = \\
 \left( \begin{array}{cc} 0 & -1 \\ 1 & 0 \end{array} \right)  \left( \begin{array}{cc} \varrho(b)^{-1} & 0 
 \\ 0 
& \varrho(b)^* \omega_b^{-1} \end{array} \right)  \left( \begin{array}{cc} 0 & 1 \\ -1 & 0 \end{array} 
\right)
=  \left( \begin{array}{cc} \varrho(b)^* \omega_b^{-1} & 0 \\ 0 & \varrho(b)^{-1} \end{array} \right).
\end{align*}
Hence $\varrho(b)^* = -\omega_b \varrho({b^{-1}})= \omega_b \varrho(h(b)^{-\mu})$ for all $b \in A$. 
For $a \in A_0$ we get \begin{align*}
\varrho(a) -\omega_b \varrho({b^{-1}}) = \varrho(a)^* +\varrho(b)^* =\varrho({ba})^* = 
\\-\omega_{ab} \varrho({a^{-1} b^{-1}}) =-\omega_{ab} (\varrho({a^{-1}}) +\varrho({b^{-1}}))
 = \omega_{ab} \varrho(a) - \omega_{ab} \varrho({b^{-1}}).\end{align*}
We get 
$$(1-\omega_{ab}) \varrho(a) = (\omega_b -\omega_{ab}) \varrho({b^{-1}}).$$
If $1-\omega_{ab} \ne 0$, then also $\omega_b -\omega_{ab} \ne 0$ and so 
$$\varrho({b^{-1}})=(\omega_b -\omega_{ab})^{-1}(1-\omega_{ab}) \varrho(a).$$ 
This forces $\varrho({b^{-1}})^* =\varrho({b^ {-1}})$. If $\characteristic R\ne 2$, 
we obtain a contradiction to 
\ref{anisotropic}. 
If $\characteristic R =2$, then $J$ contains all traces by \ref{ample}, thus there is $a \in A_0$ and $r 
\in R$ with $\varrho(a) = r +r^*$. 
Thus 
\begin{align*}\varrho({b^{-1}}) = (\omega_b -\omega_{ab})^ {-1}(1-\omega_{ab}) (r+r^*) = \\
(\omega_b -\omega_{ab})^ {-1}(1-\omega_{ab}) r + ((\omega_b -\omega_{ab})^ {-1}(1-\omega_{ab} ) r)^ * 
\in J\end{align*} by \ref{ample}, again a contradiction to \ref{anisotropic}. So $\omega_b=1$ for all $b 
\in A^*$. Thus 
$\varrho(b)^* = -\varrho({b^{-1}}) = \varrho(h(b)^{-\mu})$ for all $b \in A$. 
 \end{proof}
\begin{lemma} $f$ is $*$-skew-hermitian. 
\end{lemma}
\begin{proof}
For all $a,b \in A$ we have $f(\oa,\ob) = \varrho({ab}) -\varrho(a) -\varrho(b)$ and so 
\begin{align*} f(\oa,\ob)^* =
\varrho( {ab})^* -\varrho(a)^* -\varrho(b)^* = 
-\varrho((ab)^{-1}) +\varrho({a^{-1}}) +\varrho({b^{-1}})= \\ 
 -(\varrho({b^{-1} a^{-1} }) -\varrho({b^{-1}}) -\varrho({a^{-1}}) )= -f(-\ob,-\oa) = 
 f(\ob,-\oa) = -f(\ob,
 \oa).\end{align*}
\end{proof}
\begin{lemma} If $\characteristic R \ne 2$, then $J=\mathcal{H}(R,*)$.
\end{lemma} 
\begin{proof} For all $a,b \in A$ we have $$f(\oa,\ob) +f(\oa,\ob)^* = f(\oa,\ob) -f(\ob,\oa) 
=\varrho({[a,b]}) \in J.$$
Since $A$ is not abelian, there are $a,b \in A$ with $[a,b] \ne 1$ and thus $\varrho([a,b]) 
\ne 0$ and $f(a,b)\ne 0$. Since $R$ is a skew field by \ref{Rskew field}, we obtain 
$f(\oa ,\ob \cdot f(\oa,\ob)^{-1} r) = r$ and hence $r+r^* \in J$ for all $\in R$. Since 
$\characteristic 
R \ne 2$, 
we have $\mathcal{H}(R,*) =\{r+r^*| r \in R\}$ and thus the claim follows.
\end{proof}
\begin{prop} The the map $\pi: A/A_0 \to R/J: a+A_0 \mapsto \varrho(a)+ J$ is a pseudo-
quadratic form with 
associated 
skew-hermitian form
$f$. \end{prop}
\begin{proof} We first note that $\pi$ is well-defined since 
$$\varrho({ab}) = \varrho(a) +\varrho(b) \equiv \varrho(a)  \mod J$$ for all $a \in A, b 
\in A_0$. Since 
$r+r^* \in J$ for all $r \in R$ in any characteristic, 
we have $r^* \equiv r^* -(r+r^*) \equiv -r  \mod  J$. 
If $r \in R$, then there is a natural number $n$ and elements $h_1, \ldots, h_n \in H_0$ with 
$r =\sum_{i=1}^n  \varrho(h_i)$. Thus if $a \in A$, then we get by \ref{propf} and \ref{handf}
\begin{align*} 
\pi(\oa \cdot r ) = \pi( \oa \cdot h_1 + \ldots +\oa  \cdot h_n) = \varrho({a^{h_1} \ldots 
a^{h_n}})) + J 
=\\
\sum_{i=1}^n \varrho({a^{h_i}}) +\sum_{i < j} f(\oa^{h_i},\oa^{h_j}) +J=\\
\sum_{i=1}^n \varrho(h_i^{-\mu} h(a) h_i) + 
\sum_{i < j} \varrho(h_i^{-\mu}) f(\oa,\oa) \varrho(h_j)+J =\\
\sum_{i=1}^n \varrho(h_i)^* \varrho(a) \varrho(h_i) + \sum_{i < j} \varrho(h_i)^* f(a,a) 
\varrho(h_j)+J.\end{align*}
We have  
\begin{align*} 
0 =\varrho(1) =\varrho({aa^{-1}}) = \varrho(a) +\varrho({a^{-1}}) +f(\oa,\overline{a^{-1}})=\\
\varrho(a) +\varrho(-h(a)^{-\mu}) +f(\oa,-\oa) = \varrho(a) -\varrho(a)^* -f(\oa,\oa)\end{align*} 
and thus 
$$f(\oa,\oa) =\varrho(a) -\varrho(a)^*.$$
Thus we get 
\begin{align*}
\varrho(h_i)^* f(\oa,\oa) \varrho(h_j) = \varrho(h_i)^* (\varrho(a) -\varrho(a)^*) \varrho(h_j) =\\
\varrho(h_i)^* \varrho(a) \varrho(h_j) - (\varrho(h_j)^* \varrho(a) \varrho(h_i))^*.\end{align*}
Since $$-(\varrho(h_j)^* \varrho(a)  \varrho(h_i))^* \equiv \varrho(h_j)^* \varrho(a) \varrho(h_i)  \mod 
J,$$ 
we get $$\varrho(h_i)^* f(\oa,\oa) \varrho(h_j) \equiv \varrho(h_i)^* \varrho(a) \varrho(h_j) + 
\varrho(h_j)^* \varrho(a) \varrho(h_i) \mod  J.$$ 
Therefore 
\begin{align*}
\pi(\oa \cdot r) \equiv \ \sum_{i=1}^n \varrho(h_i)^* \varrho(a) \varrho(h_i)^* + \sum_{i\ne j} 
\varrho(h_i)^* \varrho(a) \varrho(h_j) 
\equiv 
\sum_{i,j =1}^n \varrho(h_i)^* \varrho(b) \varrho(h_j)\\
 \equiv \ (\varrho(h_1) + \ldots + \varrho(h_n))^* \varrho(a) (\varrho(h_1) + \ldots 
\varrho(h_n)) \equiv \ r^* \pi(\oa) r  \mod J.\end{align*}
The equation $$\pi(\oa + \ob) \equiv \pi(\oa) +\pi(\ob) + f(\oa,\ob)  \mod  J$$ 
for all $a,b\in A$ already follows from \ref{handf}(a).
\end{proof}
\chapter{The construction of a pseudoquadratic space}
In this chapter we continue to assume that $J$ is not commutative, so we may assume that 
$R$ is a skew field with involution 
$*$. After possibly replacing $V$ by $X(W)$ for a one-dimensional $R$-subspace $W$ of 
$C_V(A)$, we may additionally assume that $\dim_R C_V(A) =1$. This implies that $G$ acts irreduciblely on $V$. 
We have defined an anisotropic pseudo-quadratic form on the $R$-vector space $\oA$. We will make $V$ to a 
vector space over $R$ and translate this form to a form of $V$.
\newline
From now on, let $0 \ne v \in C_V(A)$ be fixed.
\begin{lemma} The map $\Phi(v,.):\oA \to \Phi(V,A)$ is an isomorphism with 
$$\Phi(v r,\oa) =\Phi(v,\oa r^*)$$
for all $\oa \in \oA$ and all $r \in R$. 
\end{lemma}
\begin{proof} If $w \in C_V(A), a \in A$ and $h \in H_0$, then by \ref{eigenschaftenphi}
\begin{align*}\Phi(w\varrho(h),\oa) = \Phi(wh,\oa) = \Phi(w,\oa^{h^{-\mu}}) =\Phi(w,\oa\cdot 
\varrho(h^{-\mu})) =
\Phi(w,\oa \cdot \varrho(h)^*).\end{align*}
Now there are $h_1, \ldots, h_n \in H_0$ with $w =\sum_{i=1}^n v\varrho(h_i)$. Then for all $\oa \in \oA$ 
we get by \ref{Phi h}
\begin{align*}\Phi(w,\oa) = \Phi(v\sum_{i=1}^n \varrho(h_i),\oa) =\sum_{i=1}^n \Phi(v\varrho(h_i),\oa) =\\
\sum_{i=1}^n \Phi(v,\oa \cdot \varrho(h_i)^*) =
\phi(v,\sum_{i=1}^n \oa \cdot \varrho(h_i)^*).\end{align*}
Thus $\Phi(v,A) =\Phi(V,A)$. This computation also shows 
the last equation. 
\end{proof}

For every $w \in V$ there are $r,s \in R$ and $\oa \in \oA$ with $w = vr + \Phi(v,\oa) + vs\mu$. 
Thus we can identify $V$ with $R \times \oA \times R$ and set $(r,\oa,s) := vr + \Phi(v,\oa) +vs\mu$.

\begin{lemma} Let $r,s \in R$ and $a,b \in A$. Then:
\begin{enumerate} \item 
$(r,\oa,s)\mu = (-s,\oa,r)$.
\item $(r,\oa,s)b = (r +f(\oa,\ob) +s\varrho(b), \oa +\ob \cdot s^*, s).$
\end{enumerate}
\end{lemma}
\begin{proof} 
\begin{enumerate}
\item By \ref{mu quadrat} we get
\begin{align*}(r,\oa,s)\mu = (vr + \Phi(v,\oa) + vs\mu)\mu =\\ vs\mu^2 + \Phi(v,\oa) +vr\mu = 
-vs +\Phi(v,\oa) +vr\mu = (-s,\oa,r).\end{align*}
\item Applying \ref{eigenschaftenphi} we get
\begin{align*}(r,\oa,s)b =(vr + [v\mu,a] -v\varrho(a) +vs\mu) b = \\
vr -v\varrho(a) +[v\mu,a] +[v\mu,a,b] +vs\mu +\Phi(vs,\ob) +vs\varrho(b) = \\
vr + vf(\oa,\ob) +vs\varrho(b) +\Phi(v,\ob \cdot s^*)+\Phi(v,\oa) +vs\mu = \\
(r + f(\oa,\ob) + s\varrho(b), \oa + \ob \cdot s^*, s).\end{align*}
\end{enumerate}
\end{proof}

We will now define a $R$-vector space structure on $V$: For $ \oa \in \oA$ and $r,s, \lambda \in R$ set 
$(r,\oa,s) \circ \lambda := (\lambda^* r, \oa \cdot \lambda, \lambda^* s)$. 
Then we have
\begin{prop} \begin{enumerate}
\item $\circ$ defines a scalar multiplication of $R$ on $V$.
\item $\circ$ commutes with the action of $G$ on $V$.
\end{enumerate}
\end{prop}
\begin{proof}
\begin{enumerate}
\item This can be easily verified.
\item We only have to show that $(w \circ \lambda) \mu = w\mu \circ \lambda$ and 
$(w\circ \lambda) b = wb \circ \lambda$ for all $b \in A$ and all $\lambda \in R$ hold.
If $r,s \in R, \oa \in \oA$, then 
\begin{align*}((r,\oa,s)\circ \lambda) \mu = (\lambda^* r, \oa \cdot \lambda,  \lambda^* s)\mu =\\
(-\lambda^*s,\oa \cdot \lambda, \lambda^ * r) 
= (-s, \oa,r) \circ \lambda = ((r,\oa,s) \mu) \circ \lambda.\end{align*}
Moreover, 
\begin{align*} ((r,\oa,s) \circ \lambda) b  = (\lambda^ * r, \oa \cdot \lambda, \lambda^* s) b\\
= (\lambda^* r +f(\oa \cdot \lambda, \ob) + \lambda^* s \varrho(b) , \oa \cdot \lambda + \ob \cdot 
(\lambda^* s)^*, 
\lambda^* s) =\\
(\lambda^* (r + f(\oa,\ob) +s\varrho(b)), (\oa +\ob \cdot s^*)\cdot \lambda, \lambda^* s) \\
=(r+f(a,b)+s\varrho(b),\oa + \ob \cdot s^*,s) \circ \lambda =((r,\oa,s)b) \circ \lambda.\end{align*}
\end{enumerate}
\end{proof}

For $r,s \in R$ and $\oa \in \oA$ set $[r,\oa,s] := (r^*,\oa,s^*)$. Then we have 
$$[r,\oa,s] \circ \lambda = (\lambda^ * r^*,\oa \cdot \lambda, \lambda^* s^*) =
[r\lambda, \oa \cdot \lambda, s\lambda]$$
for all $\lambda \in R$,
$$[r,\oa,s]\mu = (r^*,\oa,s^*)\mu = (-s^*,\oa,r^*) =[-s,\oa,r]$$ and 
\begin{align*} [r,\oa,s]b = (r^*,\oa,s^*)b = (r^* +f(\oa,\ob) +s^* \varrho(b), \oa + \ob \cdot s, s^*) = \\
[r -f(\ob,\oa)+\varrho(b)^* s, \oa + \ob \cdot s, s]\end{align*}
for all $b \in A$.

We are now able to prove our second main theorem: 
\begin{theorem}
Let $\pi:V\to R/J: \pi([r,\oa,s]) = s^* r + \varrho(a) + J$. Then
$\pi$ is a pseudo-quadratic form of Witt index $1$ and $G =\SU(\pi)$.
\end{theorem} 
\begin{proof}
Since $\oa \to \varrho(a) + J$ is an anisotropic form of $\oA$, we get that the map
${\pi:V \to R/J}$ is a pseudo-quadratic form of Witt index $1$ (see \ref{pqf}). 
Using the notation of \ref{pqf}, one can see that $a =\alpha_{(\oa,\varrho(a))}$.
If $b \in B$, then $b= a^{\mu}$ for an element $a \in A$, and thus 
$b=\alpha_{(\oa,\varrho(a))}^{\mu} =\beta_{(x,t)}$ with $x \in V$ and $t\in K$.
Thus the claim follows.
\end{proof}

We sum up our results in our main theorem.
\begin{theorem}\label{main theorem} Let $G$ be a rank one group with unipotent subgroups $A$ 
and $B$. Suppose $V$ is a cubic 
module for $G$ with $[V,G] =V$ and $C_V(G)=0$. If $A_0 \ne 1$ and the Jordan division algebra 
$J$ defined by 
$\varrho(H_0)$ is not commutative, then $J=\mathcal{H}_0(K,*)$ for a skew field $K$ with involution $*$ such 
that 
$\langle J \rangle =K$, a $K$-vector space $M$ and a pseudo-quadratic form $\pi:M \to K/J$ of Witt index 
$1$ such that 
$G \cong \SU(\pi)$. Moreover, there are $G$-submodules $U,W $ of $V$ with $U \leq W$ and $W/U 
\cong M$ as 
$G$-module. $V$ is a direct sum of $G$-submodules isomorphic to $M$ unless $\characteristic K 
=2$ and one of the 
following hold: \begin{enumerate}
\item $K$ is a quaternion algebra.
\item $K$ is a biquaternion algebra, $*$ a symplectic involution on $K$ and $J=K_* $.
\end{enumerate}
\end{theorem}
\begin{coro} Suppose that $k$, $G$, $M$ and $\Sigma$ are as in \ref{algebraic groups}. Then either $\langle A,B \rangle$ for 
$A,B \in \Sigma$ 
distinct is isomorphic to $\SL_2(F)$ for an extension field $F$ of $k$ or there is a skew field $K 
\subseteq \End_G(M)$, 
an involution $*$ of $K$ with $\langle K_{*} \rangle =K$, an irreducible submodule $M_0$ of $M$ and a $*$-skew-
hermitian 
form $f: M_0 \times M_0 \to K$ of Witt index $1$ such that $G \cong \SU(M_0,f)$. Moreover, $M$ is a direct 
sum of copies of $M_0$.
\end{coro}
\begin{proof} If $A$ and $B$ are two distinct elements in $\Sigma$ and $U(A),U(B)$ are defined 
as in \ref{algebraic groups}, 
then $G=\langle U(A), U(B) \rangle$ and $G$ acts cubically on $M$ with $U(A)_0 =A$ and $U(B)_0 =B$. By 
our main theorem, either 
$\langle A,B \rangle \cong \SL_2(F)$ for an extension field $F$ of $k$, or there is a skew field $K$ with 
involution $*$ 
such that $K_*$ generates $K$, an irreducible submodule $M_0$ of $M$ and a pseudo-quadratic 
form 
${\pi: M_0 \to K/K_*}$ with $G=\SU(\pi)$. Since $\characteristic k \ne 2$, $\pi$ is uniquely 
determined by its corresponding 
skew-hermitian form $f$. 
Since $\characteristic k \ne 2$, $M$ is the direct sum of copies of $M_0$.
\end{proof} 

\begin{remark}\rm
\begin{enumerate}
\item 
This result is similar to the following: If $\Delta$ is the generalised 
quadrangle 
corresponding to an 
involutory set $(K,K_0,\sigma)$ with $\sigma \ne 1$ and $\langle K_0 \rangle = K$ 
and if $\Gamma$ is an 
extension of 
$\Delta$ (in the sense of (21.5) in \cite{TW}), then by (21.11) of the 
same book there is a 
$K$-vector space 
$L_0$ and an anisotropic pseudo-quadratic form $\pi:L_0 \to K/K_0$ such 
that $\Gamma$ is the 
generalised 
quadrangle corresponding to $\pi$. 
\item If $G$ can be generated by three conjugates of $A_0$, then 
in characteristic not $2$ this result already follows from Theorem 2 
of \cite{T3}. 
Note that if $K=\HH$ and $*$ is defined by $r^* = -i r^{\sigma} i$, where $\sigma$ is the 
standard involution of 
$\HH$ and $i^2 =-1$, then any hermitian form proportional to a $*$-skew-hermitian form is 
$\sigma$-hermitian. 
But $K_{\sigma} =\Cent K$ does not generate $K$. Thus we cannot replace "skew-hermitian" by 
"hermitian" in our corollary. 
\end{enumerate}
\end{remark}
\chapter{Cubic rank one groups with commutative quadratic kernel}
In this chapter we investigate the case that $A_0 \ne 1$ and $R$ is a commutative field 
 and therefore $R =\Cent S$ by \ref{cases}. We will mainly deal with the case that 
 $\characteristic R \ne 2$. Later we will also have to exclude the case $\characteristic R 
 =3$. Our first result tells us that a cubic rank one group is 'locally' a unitary group.
 \begin{theorem}\label{hermitian}
  Suppose that $\characteristic R \ne 2$. 
 If $a \in A \setminus A_0$ and $G_1=\langle G_0 \cup \{a\}\rangle$, then 
 $G_1\cong \SU_3(K)$ for a quadratic extension field $K/R$ and $[V,G_1]$ is a direct sum of 
 $G_1$-standard modules.
 \end{theorem}
 \begin{proof} Set $A_1 =A \cap G_1$ and $B_1 =B \cap G_1$. Then $G_1$ is a rank one group with 
 unipotent subgroups $A_1$ and $B_1$ by \ref{root groups 2}.  
 Without loss of generality we can assume that $a=\widehat{a}$ and thus $a$ is 
 special by \ref{divisible}. 
 Set $\omega =\varrho(a)$. Then by \ref{-mu}(a) and \ref{divisible} we have  
 $\omega=\varrho(a) =\varrho(a^{-1}))=-\varrho(h(a)^{-\mu})$ and so 
 $\omega^2 =-\varrho(h(a)^{-\mu}) \varrho(e) \varrho(a) = -\varrho(e^{h(a)}) \in R^*$ by \ref{-mu}(c). 
 If $\omega\in R$, then there is $b\in A_0$ with $\varrho(b)=-\omega$ and thus 
 $\varrho(ab)= \varrho(a)+\varrho(b)=\omega-\omega=0$, a contradiction. Thus 
 $K :=R(\omega)$ is a quadratic extension of $R$. We denote the generator of the Galois group by $*$.

 Let $ v\in C_V(A)$ and set 
 $W =vK$ and $X = W \oplus \Phi(W,a) \oplus W\mu$. Then $X$ is $G_0$-invariant. Moreover, 
 for $w\in W$ we have $[w,a] =0 \in X, [w\mu,a]=\Phi(w,a)+w\varrho(a) \in X$ 
 and by \ref{f h 2} $[\Phi(w,a),a] =wf(a,a)=w (\varrho(a)-\varrho(a)^*)=
 w \cdot 2\omega \in vK\subseteq X$. 
 This shows that $X$ is also $a$-invariant. Thus $X$ is $G_1$-invariant. Now $X$ is $K$-vector 
 space via $(w_1 +\Phi(w_2,a)+w_3\mu) \cdot  r :=w_1 \cdot r+ \Phi(w_2\cdot r,a)+(w_3 \cdot 
 r)\mu$ for $w_1,w_2,w_3 \in W$ and $r \in K$. This is well-defined since 
 if $\Phi(w,a) =\Phi(u,a)$ for $u,w\in W$, then we have $\Phi(w-u,a)=0$ and thus 
 $0=[\Phi(w-u,a),a] = (w-u)f(a,a) =(w-u) \cdot 2\omega$, so $w-u=0$. 
 One easily sees that $X$ becomes a $K$-vector space this way. Since $\mathrm{dim}_K W=1$, we 
 have $\mathrm{dim}_K X =3$. We simply write $(r_1,r_2,r_3)$ for $v\cdot r_1 +\Phi(v\cdot 
 r_2,a) +(v \cdot r_3)\mu$. 
 
 For $b\in A_0$ we have 
 \begin{align*} 
 (r_1,r_2,r_3) b = (v\cdot r_1)b +\Phi(v\cdot r_2,a)b +(v\cdot r_3) \mu b=\\
 v\cdot r_1 +\Phi(v\cdot r_2,a) +(v\cdot r_3) \mu  +[(v\cdot r_3)\mu,b] = \\
 (r_1,r_2,r_3) +(v\cdot r_3)\varrho(b)= (r_1+r_3 \varrho(b),r_2,r_3).
 \end{align*}
 
Moreover, we have 
\begin{align*}
(r_1,r_2,r_3) a = (v\cdot r_1)a +\Phi(v\cdot r_2,a)a +(v\cdot r_3) \mu a =\\
v\cdot r_1 +\Phi(v\cdot r_2,a) +[\Phi(v\cdot r_2,a),a] +(v\cdot r_3) \mu a+
\Phi(v\cdot r_3),a) +(v\cdot r_3) \varrho(a)  \\
=(r_1,r_2,r_3) +v\cdot r_2 f(a,a) +(0,r_3,0) +v\cdot r_3 \omega= 
(r_1 +2r_2\omega +r_3 \omega,r_2 +r_3,r_3)
\end{align*}
and 
\begin{align*}
(r_1,r_2,r_3)\mu = (v\cdot r_1 +\Phi(v\cdot r_2,a) + (v\cdot r_3) \mu)\mu =\\
v\cdot r_3 \mu^2 +\Phi(v\cdot r_2,a) +(v\cdot r_1)\mu = (-r_3,r_2,r_1).
\end{align*}

 We define $\langle .,.\rangle: X\times X \to R $ by $\langle (r_1,r_2,r_3),(s_1,s_2,s_3) 
 \rangle =  r_1^* s_3 - r_3^* s_1 +2\omega r_2^* s_2$. Then $\langle .,.\rangle$ is 
 $*$-antihermitian $K$-linear form on $X$, and one easily computes that it is invariant 
 under $A_0$, $a$ and $\mu$. Thus $\langle .,. \rangle$ is $G_1$-invariant. In fact, 
 one easily sees that $A_0 \cup \{a,\mu\}$ generates $\PSU(X,\langle.,.\rangle)$ and thus
 $G \cong \PSU_3(K)$.\\
By \ref{X(W)} we conclude that $V$ is a direct product of modules isomorphic to $X$.    
 \end{proof}
 \begin{coro}\label{centre} If $\characteristic R \ne 2$, then $Z(A) =A_0$.
 \end{coro}
 \begin{proof}
 If $R$ is not commutative, this follows by \ref{main theorem}. If $R$ is a field, 
 then for $a \in A \setminus A_0$ and $G_1 =\langle G_0,a\rangle$ we have 
 $Z(A) \cap G_1 =A_0$ by \ref{hermitian} and thus $Z(A) =A_0$.
 \end{proof}
 \begin{coro}\label{symplectic} Suppose that $R$ is a field. Then
 $g:\oA \times \oA\to R: g(\oa,\ob) =f(\oa,\ob)-f(\ob,\oa)$ is an alternating, 
 non-degenerate, $R$-bilinear map.
 \end{coro}
 \begin{proof}
 In \ref{f R linear} we have already proved that $g$ is alternating and $R$-bilinear and that 
 $g(\oa,\ob) =\varrho([a,b]) \in R$ for all $a,b\in A$. This formula implies that 
 $\mathrm{rad}(g) =Z(A)/A_0=\{0\}$ by \ref{centre}, so $g$ is non-degenerate.   
 \end{proof}
  \begin{lemma}\label{adjoint} For $\oa,\ob \in \oA$ and $h \in H$, we have $g( a,b^h)
  = g( a^{h^{-\mu}},b)$.
 \end{lemma}
 \begin{proof} For all $v\in C_V(A)$ we have 
 \begin{align*} 
 vf(a,b^h) = v\mu(a-1)(b^h-1) = v\mu h^{-1} (a^{h^{-1}}-1) (b-1) h = \\
 vh^{-\mu^{-1}} \mu (a^{h^{-1}}-1)(b-1) h=vh^{-\mu} f(a^{h^{-1}},b) h 
 \end{align*}
 and \begin{align*} 
 vf(b^h, a) = v\mu(b^h -1)(a-1) = v\mu h^{-1} (b-1) (a^{h^{-1}} -1) h = \\
 vh^{-\mu^{-1}} \mu (b-1) (a^{h^{-1}}-1) h = vh^{-\mu} f(b,a^{h^{-1}}) h.
 \end{align*}
Here we have used that $\mu^2\in Z(H)$ by \ref{mu quadrat}(b). Thus we get 
 \begin{align*} 
 g( a,b^h) = f(a,b^h) -f(b^h,a) = \varrho(h^{-\mu}) f(a^{h^{-1}},b) \varrho(h) -
 \varrho(h^{-\mu}) f(b,a^{h^{-1}}) \varrho(h) = \\
 \varrho(h^{-\mu}) g( a^{h^{-1}},b) \varrho(h).\end{align*}
 Since $g( a^{h^{-1}},b), \varrho(h^{-\mu})\varrho(h) =\varrho(h)\varrho( h^{-\mu} )\in R \subseteq \Cent 
 S$, we get by \ref{-mu}
 $$g( a,b^h) = \varrho(hh^{-\mu}) g( a^{h^{-1}},b)= g( a^{h^{-1} hh^{-\mu}},
 b) = g( a^{h^{-\mu}},b).$$
 \end{proof}

  \begin{lemma}\label{rational}
 For all $a,b,c,d \in A$ we have 
 \begin{align*} 
 f(a,b)f(d,c)+f(c,d)f(b,a)+f(a,c)f(d,b)+\\
 f(b,d)f(c,a)+f(a,d)f(c,b)+f(b,c)f(d,a) \in R.\end{align*}
 \end{lemma}
 \begin{proof} 
 Let $a,c\in A\setminus A_0$ with $a=\widehat{a}$ and $c=\widehat{c}$. By \ref{-mu}(c) we have 
 $\varrho(ac) \varrho(h(ac)^{-\mu}) \in R$. Now we get using \ref{propf}, \ref{handf} and 
 \ref{f h 2}
 \begin{align*}
  \varrho({ac})\varrho(h({ac})^{-\mu}) =(\varrho(a) +\varrho(c) + f(a,c))(\varrho(h(a)^{-\mu})
   +\varrho(h(c)^{-\mu}) -f(c,a))=\\
 \varrho(a) \varrho(h(a)^{-\mu}) +\varrho(c) \varrho(h(c)^{-\mu}) +
 \varrho(a) \varrho(h(c)^{-\mu}) + \varrho(c) \varrho(h(a)^{-\mu}) +\\
 f(a,c)\varrho(h(a)^{-\mu}) + f(a,c)\varrho(h(c)^{-\mu}) -\varrho(a) f(c,a)-\varrho(c) f(c,a)-f(a,c)f(c,a)=\\
 \varrho(a) \varrho(h(a)^{-\mu}) +\varrho(c) \varrho(h(c)^{-\mu}) +\varrho(a) \varrho(h(c)^{-\mu})\\ 
 +\varrho(c) \varrho(h(a)^{-\mu}) 
 +g(a,c)\varrho(h(a)^{-\mu}) 
 + f(c,a)\varrho(h(a)^{-\mu}) -
 f(a,c^{h(c)})\\
 -\varrho(a) g(c,a)-\varrho(a) f(a,c)+f(c^{h(c)},a)-f(a,c)f(c,a)=\\
  \varrho(a) \varrho(h(a)^{-\mu}) +\varrho(c) \varrho(h(c)^{-\mu}) +\varrho(a) \varrho(h(c)^{-\mu}) +
  \\
  \varrho(c) \varrho(h(a)^{-\mu}) +g(a,c)  (\underbrace{\varrho(a)+\varrho(h(a)^{-\mu})}_0)+\\
  g(a,c^{h(c)})+f(c,a^{h(a)})-f(a^{h(a)},c)-f(a,c)f(c,a)=\\
\varrho(a) \varrho(h(a)^{-\mu}) + \varrho(c)\varrho(h(c)^{-\mu})-f(a,c)f(c,a) +   \\
\varrho(a)\varrho(h(c)^{-\mu}) +\varrho(c)\varrho(h(a)^{-\mu}) +
g(c,a^{h(a)}) +g(a,c^{h(c)}).
 \end{align*}
   Now $$\varrho(a) \varrho(h(a)^{-\mu})+ \varrho(c) \varrho(h(c)^{-\mu}) +g(a,c^{h(c)})+g(c,a^{h(a)})\in 
   R.$$ 
   Thus 
\begin{equation}\label{eq6}
   \varrho(a) \varrho(h(c)^{-\mu}) +\varrho(c) \varrho(h(a)^{-\mu}) -f(a,c)f(c,a)\in R.
\end{equation}
For $x,y\in A_0$ we get by \ref{handf}(c)
\begin{align*}
\varrho(ax)\varrho((h(cy)^{-\mu})+\varrho(cy)\varrho(h(ax)^{-\mu})-f(ax,cy)f(cy,ax) =\\
(\varrho(a)+\varrho(x))(\varrho(h(c)^{-\mu})+\varrho(h(y)^{-\mu})-f(y,c))+ \\
(\varrho(c)+\varrho(y))(\varrho(h(a)^{-\mu})+\varrho(h(x)^{-\mu})-f(x,a))-f(a,c)f(c,a)=\\
(\varrho(a)+\varrho(x))(\varrho(h(c)^{-\mu})+\varrho(y))+\\
(\varrho(c)+\varrho(y))(\varrho(h(a)^{-\mu})+
\varrho(x))-f(a,c)f(c,a)=\\
\varrho(a) \varrho(h(c)^{-\mu}) +\varrho(c) \varrho(h(a)^{-\mu}) -f(a,c)f(c,a)+\\
\varrho(x)(\varrho(c)+\varrho(h(c))^{-\mu})+\varrho(y)(\varrho(a)+\varrho(h(a)^{-\mu})=\\
\varrho(a) \varrho(h(c)^{-\mu}) +\varrho(c) \varrho(h(a)^{-\mu}) -f(a,c)f(c,a)\in R.
\end{align*}
Thus equation (\ref{eq6}) remains valid for arbitrary $a,c\in A\setminus A_0$. Moreover, 
equation 
(\ref{eq6}) is 
also true if $a\in A_0$ or $c\in A_0$.

Linearising equation (\ref{eq6}) in $a$ yields
   $$f(a,b)\varrho(h(c)^{-\mu}) -\varrho(c) f(b,a)-f(a,c)f(c,b)-f(b,c)f(c,a)\in R.$$
   If we linearise this expression in $c$, we get
   \begin{align*}-f(a,b)f(d,c)-f(c,d)f(b,a)-f(a,c)f(d,b)-\\f(a,d)f(c,b)-f(b,c)f(d,a)-f(b,d)f(c,a) \in R,
   \end{align*}
as desired.
   \end{proof}

   In the following we define $\langle.,.\rangle: \oA \times \oA \to R$ by 
   $\langle \oa,\ob \rangle =\frac{1}{2} g(\oa,\ob)$ for $\oa,\ob \in \oA$. 
\begin{lemma}\label{Noch ein Lemma}
Suppose that $a\in A \setminus A_0$. Then we have 
$f(a^{h(\widehat{a})},a) =-f(a,a^{h(\widehat{a})}) = \langle
 a^{h(\widehat{a})},a\rangle=2\varrho(a)^2$.
\end{lemma}
\begin{proof} We may suppose that $a=\widehat{a}$. 
We have by \ref{f h 2}(c)
$$ f(a^{h(a)},a) =-\varrho(h(a)^{-\mu}) f(a,a) =2\varrho(a)^2 $$ and
$$f(a,a^{h(a)}) = -f(a,a) \varrho(a) =-2\varrho(a)^2,$$
thus 
$$\langle a^{h(a)},a\rangle = \frac{1}{2}(f(a^{h(a)},a)-f(a,a^{h(a)}) \rangle =2\varrho(a)^2.$$  
\end{proof}
   \begin{lemma}\label{4 form} Suppose that $\characteristic A_0\ne 2$. 
   Let $Q: \oA \times \oA \times \oA \times \oA \to R$ with \begin{align*}
      Q(a,b,c,d)=\frac{1}{2}\big(f(a,b)f(d,c)+f(c,d)f(b,a)+f(a,c)f(d,b)+f(b,d)f(c,a)+\\
      f(a,d)f(c,b)+f(b,c)f(d,a)\big)+ 
      \langle a,b\rangle \langle c,d\rangle+\langle a,c\rangle\langle b,d\rangle+\langle a,d\rangle \langle b,c\rangle.\end{align*}
   Then 
   \begin{enumerate}
   \item $Q$ is a symmetric, $R$-quadrilinear form.
   \item For all $a \in \oA^*$ and all $b \in \oA$ we have 
   $Q(a,a,a,b)=\langle 6a^{h({\widehat{a}})},b\rangle$.
   \item If $\characteristic R \ne 2,3$, then for all $a,b,c \in \oA$ there is a unique element 
   $\{a,b,c\}\in \oA$ with $\langle d,\{a,b,c\}\rangle= 
   Q(a,b,c,d)$ for all $d \in \oA$. It is 
   $\{a,a,a\}=-6 a^{h(\widehat{a})}$ for all $a \in \oA^*$.
   \end{enumerate}
   \end{lemma}
   \begin{proof}
   \begin{enumerate}
   \item By \ref{rational} it follows that the image of $Q$ is in $R$. It is clear that 
   $Q$ is $R$-quadrilinear. One has 
   \begin{align*}
  Q(a,b,c,d)-Q(b,a,c,d)=\\
  \frac{1}{2} f(a,b)f(d,c)+\frac{1}{2} f(c,d)f(b,a)+\\
  \langle a,b\rangle \langle c,d\rangle -\frac{1}{2}f(b,a)f(d,c)-\frac{1}{2}f(c,d)f(a,b)-
   \langle b,a\rangle \langle c,d\rangle =\\
   \langle a,b\rangle f(d,c)+f(c,d)\langle b,a\rangle +2\langle a,b\rangle \langle c,d\rangle=\\
   2\langle a,b\rangle \langle d,c\rangle-2\langle a,b\rangle \langle c,d\rangle =0. 
   \end{align*}
   One similarly computes $Q(a,b,c,d)=Q(a,c,b,d)=Q(a,b,d,c)$. Hence $Q$ is symmetric.
   \item
   We have by \ref{f h 2}
   \begin{align*}
     Q(a,a,a,b)=
     \frac{3}{2}\left( f(a,b)f(a,a)+f(a,a)f(b,a)\right)=\\
\frac{3}{2}\left( 2g(a,b)\varrho(\widehat{a}) +2f(b,a)\varrho(\widehat{a}) +2\varrho(\widehat{a})g(b,a)+
2\varrho(\widehat{a})f(a,b)\right) =\\
3\left(-f(b,a^{h(\widehat{a})})+f(a^{h(\widehat{a})},b)\right) =6\langle a^{h(\widehat{a})},b\rangle.       
   \end{align*}
   \item The existence follows by linearising (b), the uniqueness follows from the 
   non-degeneracy of $g$. 
      \end{enumerate}
   \end{proof}
         \begin{theorem}\label{Freudenthal triple system}
   If $Q$ and $\{.,.,.\}$ are as in \ref{4 form}, then $(\oA, \{.,.,.\},\langle.,.\rangle)$ is 
   an anisotropic Freundenthal triple 
   system.
   \end{theorem}
      \begin{proof} By \ref{Noch ein Lemma} we have $\omega:=f(a^{h(\widehat{a})},a)=-
      f(a,a^{h(\widehat{a})})={1\over 2}
    g(a^{h(\widehat{a})},a)=\langle a^{h(\widehat{a})},a\rangle=2\varrho(a)^2 $.
    For all $a,b,c \in \oA $ we have by \ref{f h 2}
    \begin{align*}
       \langle c,\{a,a^{h(\widehat{a})},b\}\rangle=Q(a,a^{h(\widehat{a})},b,c)= 
        Q(b,c,a,a^{h(\widehat{a})})=\\   
\frac{1}{2}f(b,c)f(a^{h(\widehat{a})},a)+\frac{1}{2}f(a,a^{h(\widehat{a})})f(c,b)+
\frac{1}{2}f(b,a)f(a^{h(\widehat{a})},c)  +\\
   \frac{1}{2}f(c,a^{h(\widehat{a})})f(a,b)+\frac{1}{2}f(b,a^{h(\widehat{a})})f(a,c)+ 
   \frac{1}{2}f(c,a)f(a^{h(\widehat{a})},b)+\\
   \langle b,c\rangle \langle a,a^{h(\widehat{a})}\rangle+\langle b,a\rangle \langle c,a^{h(\widehat{a})}\rangle+
   \langle b,a^{h(\widehat{a})}\rangle \langle c,a\rangle =\\ f(b,c) \omega -\omega f(c,b)+\\
 \frac{1}{2}(f(b,a)\varrho(\widehat{a})f(a,c) -
f(c,a)\varrho(\widehat{a})f(a,b) -f(b,a)\varrho(\widehat{a})f(a,c)+f(c,a)\varrho(\widehat{a})f(a,b))\\
-2\omega \langle b,c \rangle+\langle b,a\rangle \langle c,a^{h(\widehat{a})}\rangle +\langle 
b,a^{h(\widehat{a})}\rangle 
\langle c,a \rangle=\\    
   2\omega \langle b,c\rangle-2\omega \langle b,c\rangle +\langle c,\langle b,a\rangle a^{h(\widehat{a})}\rangle+ 
   \langle c,\langle b,a^{h(\widehat{a})}\rangle a\rangle=
    \langle c, \langle b,a\rangle a^{h(\widehat{a})}  +\langle b,a^{h(\widehat{a})}\rangle a \rangle.
     \end{align*} Since $\langle .,. \rangle$ is non-degenerate, we get 
     $$\{a,a^{h(\widehat{a})},b\}=\langle b,a\rangle a^{h(\widehat{a})} +\langle b,a^{h(\widehat{a})}
     \rangle a$$
     and hence multiplication with $-6$ on both sides yields  
   $$\{a,\{a,a,a\},b\}= \langle b,a\rangle  \{a,a,a\} +\langle b,\{a,a,a\}\rangle a.$$
   Thus $(\oA, \{.,.,.\}, \langle .,.\rangle)$ is a Freudenthal triple system. Moreover we have 
   $Q(a) =\langle a,\{a,a,a\}\rangle =-6 \langle a,a^{h(\widehat{a}} \rangle
   =6 \omega \ne 0$, thus $Q(a)\ne 0$ for all $a \in \oA^*$.
   \end{proof}

From now on we assume that $\characteristic R \ne 2,3$. We show that $G$ is the rank one 
group corresponding to $\mathfrak{F}$ as defined in \ref{Freudenthal Moufang}.
We will identify $A$ with $\oA \times R$ in the following. 
For $t \in A_0$ we set
$(0,\varrho(t))=t$. For $a \in A \setminus A_0$, set 
$t=a\widehat{a}^{-1}\in A_0$ and $(\bar{a},\varrho(t))=a$. Then we have
\begin{lemma} For all $\bar{a}_1,\bar{a}_2 \in \oA$ and $\lambda_1,\lambda_2 \in R$ 
we have
$$(\oa_1,\lambda_1)(\oa_2,\lambda_2) = (\bar{a}_1+\bar{a}_2, \lambda_1+\lambda_2+\langle 
\oa_1,\oa_2\rangle).$$
\end{lemma}
\begin{proof} Since $A_0$ is uniquely $2$-divisible, there is
$z \in A_0$ with $z^{-2} =[\widehat{a}_1,\widehat{a}_2]$. Then we have 
\begin{align*}
( \widehat{a}_1 \widehat{a}_2 z)^{\mu^2} =  \widehat{a}_1^{-1} \widehat{a}_2^{-1} z 
=[\widehat{a}_1,\widehat{a}_2] \widehat{a}_2^{-1} \widehat{a}_1^{-1} z =\\
z^{-2} \widehat{a}_2^{-1} \widehat{a}_1^{-1} z =\\
z^{-1} \widehat{a}_2^{-1} \widehat{a}_1^{-1} = (\widehat{a}_1\widehat{a}_2 z)^{-1} .
\end{align*}
This shows that $\widehat{a}_1 \widehat{a}_2 z$ is inverted by $\mu^2$, hence 
$\widehat{\widehat{a}_1 \widehat{a}_2}=\widehat{a}_1\widehat{a}_2 z$ by \ref{divisible}. 
Suppose that $z_1,z_2 \in A_0$ with 
$\varrho(z_1)=\lambda_1$ and $\varrho(z_2)=\lambda_2$. Since $\varrho([a_1,a_2])=
g(a_1,a_2)=2\langle a_1,a_2\rangle$, we get \begin{align*}
(\bar{a}_1,\lambda_1)(\bar{a}_2,\lambda_2) = \widehat{a}_1 z_1 \widehat{a}_2z_2 =\\
\widehat{a}_1\widehat{a}_2 z_1 z_2 =\widehat{a}_1\widehat{a}_2 z z^{-1}   z_1 z_2 =
(\bar{a}_1+\bar{a}_2, \lambda_1 +\lambda_2 +\langle \bar{a}_1,\bar{a}_2 \rangle).
\end{align*}
\end{proof}

We are now ready to prove Main Theorem 3.
\begin{theorem} 
Let $a \in A \setminus A_0$ and $\lambda\in R^*$. Then we have
\begin{enumerate}
\item $(0,\lambda)\diamond \mu = (0,-\lambda^{-1})$.
\item $(\oa,0)\diamond \mu =\left(\frac{-2}{Q(\oa)} \{\oa,\oa,\oa\},0\right)$.
\item $(\oa,\lambda)\diamond \mu =\left(\frac{1}{12\lambda^2 -Q(\oa)}(2\{\oa,\oa,\oa\} +12\lambda \oa),
\frac{12\lambda}{Q(\oa) -12\lambda^2}\right)$.
\end{enumerate}
Thus the Moufang set corresponding to $G$ is isomorphic to $\mathbb{M}(\mathfrak{F})$ as 
defined in \ref{Freudenthal Moufang}.
\end{theorem}
\begin{proof}
First note that \ref{Noch ein Lemma} implies that $Q(a) =\langle \{a,a,a\},a \rangle = 6 \langle 
a^{h(\widehat{a})},a \rangle =12 \varrho(\widehat{a})^2$ for all $a \in A \setminus A_0$. 
Since $\mathfrak{F}$ is anisotropic and therefore not reduced, $Q(\oa) =12\varrho(\widehat{a})$ 
is not a square in $R^*$ by \ref{Ferrar}.
\begin{enumerate}
\item Let $t \in A_0$ with $\varrho(t) =\lambda$. Then 
$h({t\diamond \mu}) = \mu^{-2} h(t)^{-1}$ by \ref{mu a tau}(b). 
Since $\varrho(\mu^{-2}) =-1$, we get $\varrho(t\diamond \mu) =\varrho(\mu^{-2}) \varrho(t)^{-1} = 
(-1) \varrho(t)^{-1} =-\lambda^{-1}$, hence the claim follows.
\item Since $\widehat{a}$ is special and $\widehat{a^h}=\widehat{a}^h$ 
for all $h\in H$ by \ref{divisible}, we get by \ref{special elements in MS}
\begin{align*}
\widehat{a}\diamond \mu = \widehat{a}^{-1}\diamond \mu(\widehat{a})^{-1}\diamond \mu =
\widehat{a}^{-h(\widehat{a})^{-1} \mu^2}=\widehat{a}^{-\mu^2 h(\widehat{a})^{-1}}
=\widehat{a}^{h(\widehat{a})^{-1}}= 
\widehat{a}^{h(\widehat{a}) h(\widehat{a})^{-2}} 
 \end{align*}
thus $(\oa,0)\diamond \mu = (\oa^{h(\widehat{a})} \cdot \varrho(\widehat{a})^{-2},0) =
(\frac{-1}{6}\{\oa,\oa,\oa\} \cdot \frac{12}{Q(\oa)},0) =(-\frac{2}{Q(\oa)}\{\oa,\oa,\oa\},0)$.
\item Let $t \in A_0$ with $\varrho(t) =\lambda$ and $b\in A$. Then we have by \ref{f a mu}
\begin{align*}f(\widehat{a},b) = f(\widehat{a}t,b) =
\varrho(\widehat {a}t)f((\widehat{a}t)\diamond \mu, b) = \\
(\varrho(\widehat{a})+\varrho(t))((\widehat{a}t)\diamond\mu,b)=
(\lambda +\varrho(\widehat{a})) f((at)\diamond\mu,b),
\end{align*}
thus by \ref{f h 2}
\begin{align*}
f( (\widehat{a}t)\diamond \mu,b) = (\lambda +\varrho(\widehat{a}))^{-1} f(\widehat{a},b) = 
\frac{1}{\lambda^2-\varrho(\widehat{a})^2} (\lambda -\varrho(\widehat{a}))f(a,b)=\\
\frac{1}{\lambda^2-\varrho(\widehat{a})^2} (\lambda f(a,b) -f(a^{h(\widehat{a})},b))=\\
f(\oa \cdot \frac{\lambda}{\lambda^2-\varrho(\widehat{a})^2} 
+ \frac{1}{6} \{\oa,\oa,\oa\} \cdot \frac{1}{\lambda^2 -\varrho(\widehat{a})^2},\ob)=\\
f(\oa \cdot \frac{\lambda}{\lambda^2 - \frac{1}{12} Q(\oa)} +\frac{1}{6} \{\oa,\oa,\oa\}
\cdot \frac{1}{\lambda^2 -\frac{1}{12}Q(\oa)},\ob)=\\
f((\oa\cdot 12\lambda  + 2\{\oa,\oa,\oa\}) \cdot \frac{1}{12 \lambda^2-Q(\oa)},\ob).  
\end{align*}
This shows that $$\overline{(\widehat{a}t)\diamond \mu}  = (\oa\cdot 12\lambda  + 2\{\oa,\oa,\oa\}) \cdot 
\frac{1}{Q(12 \lambda^2-Q(\oa)}.$$
Now we have again by \ref{mu a tau}(b) and \ref{mu quadrat} \begin{align*}
\varrho({(\widehat{a}t)\diamond \mu}) = -\varrho(\widehat{a}t)^{-1} = -(\varrho(\widehat{a})+\lambda)^{-1} =\\
\frac{1}{\lambda^2-\varrho(\widehat{a})^2}(\varrho(\widehat{a})-\lambda).
\end{align*}
Moreover, we have by \ref{formula sim}(b) and (d) and \ref{divisible}
\begin{align*}
\varrho(((\widehat{a}t)\diamond\mu)^{-1}) =\varrho(((\widehat{a}t)\diamond\mu)^\sim) =
\varrho((\widehat{a}t)^{-1}\diamond\mu)= \varrho((\widehat{a}^{-1}t^{-1})\diamond\mu)=\\
\varrho(\big(\widehat{a^{-1}}t^{-1}\big)\diamond\mu)= 
\frac{1}{(-\lambda)^2-\varrho\big(\widehat{a^{-1}}\big)^2}
(\varrho\big(\widehat{a^{-1}}\big)+\lambda)=\\
\frac{1}{\lambda^2-\varrho(\widehat{a})^2}(\varrho(\widehat{a})+\lambda).
\end{align*}
Set $x=(\widehat{a}t)\diamond \mu$. 
There is $s\in A_0$ such that $\widehat{x}=xs=((\widehat{a}t)\diamond \mu)s$. 
Since $\varrho(\widehat{x})=\varrho(\widehat{x}^{-1})$ by \ref{divisible}, we get by \ref{mu a tau} and \ref{-mu}
\begin{align*} 
\frac{1}{\lambda^2-\varrho(\widehat{a}^2)}(\varrho(\widehat{a})-\lambda)+\varrho(s) =\varrho(\widehat{x}) =\\
\varrho(\widehat{x}^{-1}) =\varrho((\widehat{a}t)\diamond \mu)^{-1}s^{-1})=
\frac{1}{\lambda^2-\varrho(\widehat{a})^2}(\varrho(\widehat{a})+\lambda) -\varrho(s)
\end{align*}
thus $$\varrho(s) =\frac{\lambda}{\lambda^2-\varrho(\widehat{a})^2}=
\frac{12\lambda}{12\lambda^2 -Q(\oa)}.$$ Thus the claim follows.
\end{enumerate} 
\end{proof}
\begin{coro} Suppose that the hypothesis of \ref{quadratic pairs} hold and that additionally 
$\characteristic K\ne 3$. Then one of the following holds:
\begin{enumerate}
\item There is a field or skew field $L$ containg $K$, a finite-dimensional $L$-vector space $V$ and a 
$*$-skew-hermitian form $f:V\times V \to L$ of Witt index $1$ such that $G\cong \SU(f)$.
\item There is a field extension $K/R$ such that $G/Z(G)$ 
is isomorphic to a semisimple algebraic group 
defined over $R$ of type ${}^{3}\mathrm{D}^9_{4,1},{}^{6}\mathrm{D}^9_{4,1} 
{}^{2}\mathrm{E}^{35}_6, 
\mathrm{E}^{66}_{7,1}$ or $\mathrm{E}^{133}_{8,1}$.
\end{enumerate}
\end{coro}
\begin{proof}
This follows from Main Theorem 2 and 3 along with \ref{Freudenthal Moufang}.
\end{proof}
\begin{remark}\rm
\begin{enumerate}
\item If $G$ is an algebraic group of type $\mathrm{E}_8$ defined over an algebraically closed 
field $K$ of characteristic not $2$, then by \cite[Theorem 1]{PS} 
there is no homomorphism $\varphi:G\to \GL_n(K)$ 
of algebraic groups such that the image of $G$ is generated by quadratic elements. Therefore there is 
evidence that the case $\mathrm{E}^{133}_{8,1}$ does not occur. However, we do not assume that $G$ is 
Zariski closed in $\GL(V)$, and therefore we cannot exclude this case.
\item Again by \cite[Theorem 1]{PS} 
several groups appearing in the list above have in fact more than one cubic module, and we do not 
know yet how we can distinguish them with our methods. 
\end{enumerate}
\end{remark} 
\backmatter

\end{document}